\documentclass[11pt]{article}
\usepackage{amsmath,amsthm,amssymb,
color,epic,eepic,graphics,
cite}
\renewcommand{\baselinestretch}{1.2}

\def\baselinestretch{1.4}
\newlength{\minitwocolumn}
\newcommand{\N}{{\mathbb N}} 
\newcommand{\Z}{{\mathbb Z}} 
\newcommand{\R}{{\mathbb R}} 
\newcommand{\C}{{\mathbb C}} 
\newcommand{\FF}{{\mathbb F}} 

\newcommand{\cC}{{\cal C}}
\newcommand{\cD}{{\cal D}}
\newcommand{\cA}{{\cal A}}
\newcommand{\cB}{{\cal B}}
\newcommand{\cF}{{\cal F}}
\newcommand{\cK}{{\cal K}}
\newcommand{\cU}{{\cal U}}
\newcommand{\cV}{{\cal V}}

\newcommand{\cR}{{\cal R}}

\newcommand{\cQ}{{\cal Q}}
\newcommand{\cM}{{\cal M}}
\newcommand{\cN}{{\cal N}}
\newcommand{\cH}{{\cal H}}
\newcommand{\cE}{{\cal E}}

\newcommand{\bR}{{\overline{R}}}
\newcommand{\hL}{\widehat{L}}
\newcommand{\hR}{\widehat{R}}

\newcommand{\trL}{{{}^tL}}
\newcommand{\trho}{{{\rho}}}
\newcommand{\tb}{{{b}}}
\newcommand{\tbb}{{{\bar{b}}}}
\newcommand{\tc}{{{c}}}
\newcommand{\bb}{{{\bar{b}}}}

\newcommand{\cL}{{\cal L}}

\renewcommand{\H}{{\cal H}}
\newcommand{\la}{\lambda}

\newcommand{\al}{\alpha}

\newcommand{\vep}{\varepsilon}

\newcommand{\bep}{\bar{\epsilon}}

\newcommand{\hb}{\widehat{b}}
\newcommand{\hc}{\widehat{c}}
\newcommand{\s}{\sigma}
\newcommand{\hf}{\widehat{f}}

\newcommand{\hV}{\widehat{V}}

\newcommand{\sgn}{{\rm sgn}}
\newcommand{\hsgn}{\widehat{\rm sgn}}

\newcommand{\edet}{q\mbox{-}{\rm det}}

\newcommand{\noi}{{\noindent}}
\newcommand{\nn}{{\nonumber}}
\newcommand{\bea}{\begin{eqnarray}}
\newcommand{\ena}{\end{eqnarray}}
\newcommand{\beit}{\begin{itemize}}
\newcommand{\enit}{\end{itemize}}

\newcommand{\be}{\begin{eqnarray*}}
\newcommand{\en}{\end{eqnarray*}}
\newcommand{\lb}[1]{\label{#1}}

\newcommand{\ds}[1]{{\displaystyle #1 }}

\newcommand{\End}{{\rm End}}

\newcommand{\id}{{\rm id}}
\newcommand{\Ad}{{\rm Ad}}

\newcommand{\wt}{{\rm wt}}
\newcommand{\Ker}{{\rm Ker}\ }

\newcommand{\te}{{\Theta_p}}
\newcommand{\tes}{{\Theta}_{p^*}}
\newcommand{\qdet}{{q}\mbox{\rm -det}}

\def\infq4p#1{{(#1;q^4,p)_\infty}}

\newcommand{\hd}{\widehat{d}}
\newcommand{\hi}{\widehat{i}}
\newcommand{\hj}{\widehat{j}}
\newcommand{\hk}{\widehat{k}}
\newcommand{\hl}{\widehat{l}}
\newcommand{\hrho}{\widehat{\rho}}

\newcommand{\tot}{\widetilde{\otimes}}

\newcommand{\mmatrix}[1]{\begin{matrix} #1 \end{matrix}}
\newcommand{\mat}[1]{\left(\mmatrix{#1}\right)}

\font\teneufm=eufm10
\font\seveneufm=eufm7
\font\fiveeufm=eufm5
\newfam\eufmfam
\textfont\eufmfam=\teneufm
\scriptfont\eufmfam=\seveneufm
\scriptscriptfont\eufmfam=\fiveeufm

\let\goth\frak
\newcommand{\slth}{\widehat{\goth{sl}}_2}
\newcommand{\slt}{\goth{sl}_2}

\newcommand{\slnh}{\widehat{\goth{sl}}_N}
\newcommand{\sln}{\goth{sl}_N}
\newcommand{\g}{\goth{g}}
\newcommand{\gl}{\goth{gl}}
\newcommand{\glth}{\widehat{\goth{gl}}_2}
\newcommand{\gh}{\widehat{\goth{g}}}

\newcommand{\Aqp}{{\cal A}_{q,p}}

\newcommand{\Bqla}{{{\cal B}_{q,\lambda}}}

\newcommand{\gS}{\goth{S}}
\newcommand{\ghbig}{\widehat{\mbox{\fourteeneufm g}}}  

\newcommand{\glnhbig}{\widehat{\mbox{\fourteeneufm gl}}_N}
\newcommand{\glnh}{\widehat{\goth{gl}}_N}
\newcommand{\gln}{{\goth{gl}}_N}
\newcommand{\h}{H}
\newcommand{\hh}{\goth{h}}

\font\fourteeneufm=eufm10 scaled\magstep2    
   
\makeatletter
\@addtoreset{equation}{section}
\makeatother

\newtheorem{thm}{Theorem}[section]
\newtheorem{prop}[thm]{Proposition}
\newtheorem{lem}[thm]{Lemma}
\newtheorem{cor}[thm]{Corollary}

\newtheorem{df}{Definition}[section]
\newtheorem{dfn}[thm]{Definition}

\begin{document}

\begin{center}
{\Large  }
\end{center}


\begin{center}
{\Large \bf Elliptic Quantum Groups $U_{q,p}(\glnh)$ and $E_{q,p}(\glnh)$}\\[10mm]
{\large Hitoshi Konno\footnote{E-mail: hkonno0@kaiyodai.ac.jp}}\\[5mm]
{\it Department of Mathematics, Tokyo University of Marine Science and Technology,}\\ 
{\it  Etchujima,Tokyo 135-8533, Japan}\\

\end{center}

\begin{center}
Dedicated to Professor Masatoshi Noumi on his 60th birthday. 
\end{center}

\begin{abstract}
We reformulate a central extension of Felder's elliptic quantum group 
in the FRST formulation as a topological algebra $E_{q,p}(\glnh)$ over the ring of formal power series in $p$. 
We then discuss the isomorphism between $E_{q,p}(\glnh)$ and the elliptic algebra $U_{q,p}(\glnh)$ 
of the Drinfeld realization. An evaluation $H$-algebra  homomorphism from $U_{q,p}(\glnh)$ to a dynamical extension of 
the quantum affine algebra $U_{q}(\glnh)$ resolves the problem into the one discussed by Ding and Frenkel in the trigonometric case. 
We also provide some useful formulas for the elliptic quantum determinants.  
\end{abstract}

\section{Introduction}
An elliptic quantum algebra is an associative algebra related to an elliptic solution to the Yang-Baxter equation (YBE) or the dynamical Yang-Baxter equation (DYBE). Equipped with a co-algebra structure it is called  the elliptic quantum group (EQG). 
Depending on YBE or DYBE, the corresponding EQG is called the vertex type or the face type, respectively\cite{Fro,JKOStg}. 
Through this paper  we use the terminology DYBE\cite{Felder} as an equation equivalent to  the face type YBE or the star triangle equation (see for example \cite{JMO}).

Let  $\g$ and $\gh$ denote a simple Lie algebra and an (untwisted) affine Lie algebra, respectively.
In known quantum groups, such as the Yangian $Y(\g)$ (or its double $\cD Y(\g)$) associated to the rational solutions  to the YBE 
and the affine quantum group $U_q(\gh)$ associated to the trigonometric solutions there are some  different formulations depending on the types of the generators. In particular for $U_q(\gh)$ they are the Drinfeld-Jimbo  formulation\cite{Drinfeld,Jimbo} in terms of an analogue of the Chevalley generators, Drinfeld\rq{}s new realization\cite{Dr} whose generators, called the Drinfeld generators, are  natural analogues of those in the  the loop algebras $\g[t,t^{-1}]$, and the 
Faddeev-Reshetikhin-Semenov-Tian-Shansky-Takhtadjan (FRST) formulation \cite{FRT,RS} in terms of the $L$ operators satisfying the $RLL$ relations 
associated with the $R$  matrix, a solution to the YBE. 
The isomorphisms among these three have been discussed by several authors\cite{Dr,Beck,Jing98,DF,HM,Iohara}.   
 
Correspondingly there are three different formulations of EQGs: $\Aqp(\slnh)$ and $\Bqla(\gh)$\cite{JKOStg} in terms of the Chevalley type generators, $U_{q,p}(\gh)$\cite{Konno,JKOS} and $E_{\tau,\eta}(\g)$\cite{ER,EF} in terms of the Drinfeld generators and 
 $\Aqp(\slth)$\cite{FIJKMY} and $E_{\tau,\eta}(\gln)$\cite{EV, Felder,FV,KNR} in terms of the $L$ operators. 
Here only $\Aqp(\slnh)$ is the vertex type EQG, which is related to Baxter-Belavin\rq{}s elliptic $R$ matrix\cite{Baxter,Belavin}.  The others are the face type which are related to the elliptic solutions to the face type YBE, for example \cite{JMO}.
These have their own co-algebra structures: the quasi-Hopf algebra structure\cite{DrinfeldQH} for $\Aqp(\slnh)$, $\Bqla(\gh)$\cite{JKOStg} and $E_{\tau,\eta}(\slt)$\cite{EF}, and the Hopf algebroid structure\cite{EV,KR} for $E_{\tau,\eta}(\gln)$\cite{FV,KNR,Hartwig} and $U_{q,p}(\gh)$\cite{Konno09}.

As like the cases in $Y(\g), \cD Y(\g)$ and $U_q(\gh)$, each formulation 
has both advantages and disadvantages. The quasi-Hopf algebra formulations   
$\Aqp(\slnh)$ and $\Bqla(\gh)$\cite{JKOStg} are suitable for studying formal algebraic structures such as the universal elliptic dynamical $R$ matrices, the universal form of the dynamical $RLL$ relations etc., 
 but it is hard to derive concrete representations due to the complexity of the quasi-Hopf twist operation. 

The Drinfeld realization $U_{q,p}(\gh)$ is suitable for studying both finite 
and infinite dimensional representations\cite{JKOS,KK03,KK04,KKW,Konno09,Konno08,KO} due to the nature 
of the Drinfeld generators. Recent developments include  
a characterization of the finite dimensional representations in terms of a theta function analogue\cite{Konno09} of the Drinfeld polynomials\cite{Dr,CP}
and a clarification of the quantum $Z$-algebra structures of the infinite dimensional representations\cite{FKO}.  
 An application to the algebraic analysis of the solvable lattice models\cite{JMBook} 
also have made a great success\cite{JKOS,JKOPS,KKW,Konno08,CKSW}. See also rather older works \cite{LP,AJMP,LasPu} whose results, in particular the vertex operators and the screening operators, are able to be reformulated by the representation theory of $U_{q,p}(\slnh)$\cite{KK03,KKW}. In addition  there are deep relationships between $U_{q,p}(\gh)$ and the deformed $W(\g)$ algebras: the generating functions of the Drinfeld generators ( the elliptic currents ) $e_j(z)$ and $f_j(z)$ of $U_{q,p}(\gh)$ are identified with the screening currents of the deformed $W(\g)$ algebras of the coset type\cite{Konno,JKOS,KK04W,FKO}.

The FRST formulation is suitable for studying finite dimensional representations by a fusion procedure or by taking a coproduct. In this way finite dimensional representations of $E_{\tau,\eta}(\gln)$ have been studied well \cite{FV,KNR,KN} (see also \cite{JKMO}) and applied to the elliptic Ruijsenaars models\cite{FV97, FV96NP},  the elliptic hypergeometric series\cite{KNR,KN,Rosengren09}, the partition function of the solvable lattice model\cite{Rosengren11,PRS} and the elliptic Gaudin model\cite{RST}. 

On the other hand in order to formulate infinite dimensional representations of $E_{\tau,\eta}(\gln)$ one needs it's central extension. 
There are two different proposals provided by \cite{EF} and \cite{JKOStg,JKOS}, respectively.  
Accordingly $E_{\tau,\eta}(\slt)$ in \cite{EF} and $U_{q,p}(\gh)$ in \cite{JKOS,FKO,Konno} have been the two proposals for their Drinfeld realizations. 
However the 
isomorphism between $E_{\tau,\eta}(\gln)$ in the FRST formulation and 
neither of these two Drinfeld realizations 
has been discussed precisely. 

The aim of this paper is to establish the isomorphism between $U_{q,p}(\glnh)$ and a central extension of $E_{\tau,\eta}(\gln)$ in the FRST formulation 
as a Hopf algebroid. For this purpose, we first reformulate $E_{\tau,\eta}(\gln)$ as a topological algebra over the ring of formal power series in $p$ and at the same time we give a central extension of it according to the argument in \cite{JKOStg,JKOS}. 
We denote the resultant algebra by $E_{q,p}(\glnh)$, where the generators 
are clear and their defining relations are well defined in the $p$-adic topology as in $\Aqp(\slth)$\cite{FIJKMY} and $U_{q,p}(\gh)$\cite{FKO}. 
Secondly we discuss dynamical representations of $U_{q,p}(\glnh)$.  We  especially introduce an evaluation  $H$-algebra homomorphism from 
$U_{q,p}(\glnh)$ to a dynamical extension of the quantum affine algebra $U_{q}(\glnh)$. This allows us to obtain the 
 dynamical representations ( of both finite and infinite dimensional ) from any representations of $U_{q}(\glnh)$. 
As a result   the problem resolves itself into the one discussed by  Ding and Frenkel in the trigonometric case \cite{DF}.

A part of the results have been reported in the workshops ``Recent advances in quantum integrable systems 2012\rq{}\rq{}, 10-14 September 2012, Angers, France   
and ``Elliptic Integrable Systems and Hypergeometric Functions\rq{}\rq{}, 
15-19 July 2013, Lorentz Center, Leiden, the Netherlands.

This paper is organized as follows.
In Section 2 preparing notations and conventions we introduce 
the elliptic dynamical $R$ matrix. In Section 3 we define $U_{q,p}(\glnh)$ and $E_{q,p}(\glnh)$ as  topological algebras over the ring of formal power series in $p$. We also give the trigonometric  $(p=0)$ counter parts of them. 
In section 4 we show that both $U_{q,p}(\glnh)$ and $E_{q,p}(\glnh)$ are $H$-algebras (Proposition \ref{UqpHalg} and \ref{EqpHalg}). 
Then we introduce an $H$-Hopf algebroid structure to them.
In Section 5 we introduce dynamical representations of $U_{q,p}(\gh)$ and  give a construction of the  evaluation dynamical representations from any representations of $U_q(\glnh)$.  In Section 6 we discuss an isomorphism between $U_{q,p}(\glnh)$ and $E_{q,p}(\glnh)$ as an $H$-Hopf algebroid. Our arguments mainly follow those by Ding and Frenkel in the trigonometric case\cite{DF} with some additional formulas for the lower rank subalgebras of $E_{q,p}(\glnh)$, which make the induction process more transparent. In particular by making use of the evaluation dynamical representations in Sec. 5 our proof on the injectivity resolves itself into the 
results in \cite{DF}. Appendix A contains a definition of the quantum affine algebra $U_q(\glnh)$ which we  use in Sec.5. 
In Appendix B we list the formulas necessary for discussing the evaluation dynamical representations. 
In Appendix C we summarize the formulas which identify a combination of the Gauss components of the $L$ operator with the elliptic currents of $U_{q,p}(\glnh)$. In Appendix D we summarize some formulas on adding `fractional powers in $z$' which clarify a connection between $U_{q,p}(\gh)$ in the current paper and the previous one in \cite{Konno,JKOS,KK03,KK04,Konno09}.    
Appendix E contains some formulas for the elliptic quantum determinants.

\section{The $R$-matrices}\lb{defellRmat}

 Let $\{ \epsilon_j\ (1\leq j\leq N)\}$ be 
the orthonormal basis in $\R^N$ with the inner product
$( \epsilon_j, \epsilon_k )=\delta_{j,k}$.
Setting 
$\displaystyle{\bar{\epsilon}_j=\epsilon_j-\epsilon,~
\epsilon=\frac{1}{N}\sum_{j=1}^N \epsilon_j,
}$
we define the weight lattice ${\cal P}$ of $A_{N-1}$ type by 
$\ds{
{\cal P}=\sum_{j=1}^N {\mathbb{Z}}\ \bar{\epsilon}_j.
}$
Let $I=\{1,2,\cdots,N-1\}$.   We set 
$\ds{
\alpha_j=\bar{\epsilon}_j-\bar{\epsilon}_{j+1}, \ \bar{\Lambda}_j=\bep_1+\cdots +\bep_{j}\ (j\in I)
}$
and define $\cQ=\Z \al_1\oplus \cdots \oplus \Z \alpha_{N-1}$ and $\bar{\hh}^*=\C \bar{\Lambda}_1\oplus \cdots \oplus \C \bar{\Lambda}_{N-1}$. 
We also define elements $h_{\bep_j}\ (1\leq j\leq N)$ in the dual space $\bar{\hh}$ by  $<\bep_i,h_{\bep_j}>=(\bep_i,{\bep_j})=
\delta_{j,k}-\frac{1}{N}$. Setting 
$h_j=h_{\bar{\epsilon}_j}-h_{\bar{\epsilon}_{j+1}}\ (j\in I)$ we have $<\bar{\Lambda}_i,h_{j}>=\delta_{i,j}$ so that   
$\bar{\hh}=\C h_1\oplus \cdots \oplus \C h_{N-1}$.
 For $\alpha=\sum_j a_j \bep_j \in \bar{\hh}^*$, we define $h_\alpha \in \bar{\hh}$  by
 $h_\alpha=\sum_j a_j h_{\bep_j}$ and $h_0=0$. 
 We also need two more elements $c$ and $\Lambda_0$ satisfying 
$<\Lambda_0,c>=1, <\Lambda_0,h_j>=0=<\bar{\Lambda}_j,c>\ (1\leq j\leq N)$. 
 We regard $\bar{\hh}\oplus  \bar{\hh}^*$ as a Heisenberg 
algebra  by
\bea
&&~[h_{\bar{\epsilon}_j},\bar{\epsilon}_k]
=( \bar{\epsilon}_j,\bar{\epsilon}_k ),\qquad [h_{\bar{\epsilon}_j},h_{\bar{\epsilon}_k}]=0=
[\bar{\epsilon}_j,\bar{\epsilon}_k].\lb{HA1}
\ena
 
We also introduce the dynamical parameters $P_{\bep_j}\ (j\in I)$ and their duals $Q_{\bep_j}$. They are the Heisenberg algebra defined by
\begin{eqnarray}
&&[P_{\bar{\epsilon}_j}, Q_{\bar{\epsilon}_k}]=
( \bar{\epsilon}_j, \bar{\epsilon}_k ), \qquad 
[P_{\bar{\epsilon}_j}, P_{\bar{\epsilon}_k}]=0=
[Q_{\bar{\epsilon}_j}, Q_{\bar{\epsilon}_k}].\lb{HA2}
\ena
We set 
$P_\alpha=\sum_j a_j P_{\bep_j}$ for $\alpha=\sum_j a_j \bep_j$ and $P_0=0$ etc.
In particular we set $P_{j}=P_{\alpha_j}=P_{\bep_j}-P_{\bep_{j+1}}$ and 
$Q_{j}=Q_{\alpha_j}=Q_{\bep_j}-Q_{\bep_{j+1}}$. 


For the abelian group 
$\cR_Q= \sum_i\Z Q_{\bep_i}$, 
we denote by $\C[\cR_Q]$ the group algebra over $\C$ of $\cR_Q$. 
We denote by $e^{Q_\al}$ the element of $\C[\cR_Q]$ corresponding to $Q_\al\in \cR_Q$. 
These $e^{Q_\al}$ satisfy $e^{Q_\al} e^{Q_\beta}=e^{Q_{\al}+Q_{\beta}}$
 and $(e^{Q_\al})^{-1}=e^{-Q_\al}$. 
In particular, $e^0=1$ is the identity element. 

Now let us introduce a dynamical extension of $\bar{\hh}$ and $\bar{\hh}^*$:  $H=\bar{\hh}\oplus \sum_j \C P_{\bep_j}+\C c= \sum_{j}\C(P+h)_{\bep_j}+\sum_j\C P_{\bep_j}+\C c$  and   $H^*=\bar{\hh}^*+\sum_j \C Q_{\bep_j}+\C \Lambda_0$.  Through this paper we often use the abbreviation $(P+h)_{\bep_j}$ for $P_{\bep_j}+h_{\bep_j}$. We have the paring:   
$<Q_\al,P_\beta>=(\al,\beta)=<\al,h_{\beta}>\ \al,\beta\in \bar{\hh}^*,\ <\Lambda_0,c>=1$, and the others vanish.  
Let $\cM_{H^*}$ be the field of meromorphic functions on $H^*$.  We denote by $\hf=f(P+h,P)$ an element of $\cM_{H^*}$, where $P+h=\sum_ja_j(P+h)_{\bep_j}, P=\sum_jb_j P_{\bep_j}\in H$.  The function $\hf$ is evaluated at $\mu\in H^*$ as  
$\hf(\mu)=f(<\mu,P+h>, <\mu,P>)$ etc..  
Hereafter we set $\FF=\cM_{H^*}$.   

Let ${\hbar}$ and $p$ be indeterminates. We set $q=e^{\hbar}$. 
Through this paper we also use $p^*=pq^{-2c}$. 

The following notations are often used.
\be 
&&[n]_q=\frac{q^n-q^{-n}}{q-q^{-1}},\qquad \Theta_p(z)=(z;p)_{\infty}(p/z;p)_\infty(p;p)_\infty,\\
&&(x;q_1,q_2,\cdots,q_k)_\infty=\prod_{n_1,n_2,\cdots,n_k=0}^\infty(1-x q_1^{n_1} q_2^{n_2}\cdots q_k^{n_k}),\qquad \{z\}=(z;p,q^{2N})_\infty,\\
&&(x_1,x_2,\cdots,x_l;q_1,q_2,\cdots,q_k)_\infty=\prod_{i=1}^l(x_i;q_1,q_2,\cdots,q_k)_\infty.
\en

\subsection{The elliptic dynamical $R$-matrices}\lb{edR}

 Let $V=\oplus_{i=1}^N\C v_i$, 
$E_{i,j}v_k=\delta_{j,k}v_i$. 
We consider the following elliptic dynamical $R$-matrices $R^{\pm}(z,s)\in \End (V\otimes V)$ of type $A_{N-1}^{(1)}$. For $s\in H$, 
\bea
R^{\pm}(z,s)&=&{\trho}^{\pm}(z){\bR}(z,s),
\lb{tRmat}\\
{\bR}(z,s)
&=&
\sum_{j=1}^{N}
E_{jj}\otimes E_{jj}+
\sum_{1 \leq j<l \leq N}
\biggl(
\tb_{}(z,s_{j,l})
E_{jj}
\otimes E_{ll}+
{\tbb}_{}(z)
E_{ll}\otimes E_{jj}
\nonumber\\
&&\qquad\qquad +
\tc_{}(z,s_{j,l})
E_{jl}\otimes E_{lj}+\bar{{\tc}}_{}(z,s_{j,l})E_{lj}\otimes E_{jl}
\biggr),\lb{tbRmat}
\end{eqnarray}
where $\displaystyle{s_{j,l}=
s_{\bep_j}-s_{\bep_l}
\quad (1\leq j<l \leq N)}$ and 
\begin{eqnarray}
&&\trho^+(z)=q^{-\frac{N-1}{N}}
\frac{\{q^2z\}
\{q^{2N-2}z\}
\{p/z\} \{pq^{2N}/z\}
}{
\{z\}
\{q^{2N}z\}
\{pq^2/z\} \{pq^{2N-2}/z\}
},\label{def:rhopz}\\
&&\trho^-(z)=q^{\frac{N-1}{N}}
\frac{\{pq^2z\}
\{pq^{2N-2}z\}
\{1/z\} \{q^{2N}/z\}
}{
\{pz\}
\{pq^{2N}z\}
\{q^2/z\} \{q^{2N-2}/z\}
},\label{def:rhomz}\\
&&\tb_{}(z,s)=q
\frac{\Theta_p(q^2q^{2s})\Theta_p(q^{-2}q^{2s})\Theta_p(z)
}{
\Theta_p(q^{2s})^2\Theta_p(q^2z)
},\\
&&{\tbb}_{}(z)=q
\frac{\Theta_p(z)
}{\Theta_p(q^2z)
},\\
&&\tc_{}(z,s)=
\frac{\Theta_p(q^2)\Theta_p(q^{2s}z)
}{
\Theta_p(q^{2s})\Theta_p(q^2z)
},\qquad \bar{c}(z,s)=c(z,-s).\lb{rmatcompz}
\end{eqnarray}
We also denote by $R^{\pm *}(z,s)$ the $R$ matrices obtained from $R^{\pm}(z,s)$ by replacing 
$p$ with $p^*$. 
Note that 
\bea
&&\trho^+(z p)=q^{-2\frac{N-1}{N}}\trho^-(z). \lb{trhopm}
\ena
In particular, 
\bea
R^-(z,s)^{-1}=PR^+(z^{-1},s)P.\lb{tRpRm}
\ena
Furthermore if we set 
\bea
&&\rho(z)=\frac{\trho^{+*}(z)}{\trho^+(z)}\lb{rho}
\ena
where $\trho^{\pm*}(z)=\trho^\pm(z)|_{p\to p^*}$, we have 
\bea
\rho(z)^{-1}=\rho(z^{-1}),\qquad \rho(z)=\frac{\trho^{-*}(z)}{\trho^-(z)}.\lb{trhorhos} 
\ena

\begin{prop}\cite{Felder}
The $R^+(z,s)$ satisfies the following dynamical Yang-Baxter equation.
\bea
&&R^{+(12)}(z_1/z_2,P+\pi_V(h)^{(3)})R^{+(13)}(z_1,P)R^{+(23)}(z_2,P+\pi_V(h)^{(1)})\nn\\
&&\qquad=R^{+(23)}(z_2,P)R^{+(13)}(z_1,P+\pi_V(h)^{(2)})R^{+(12)}(z_1/z_2,P),\lb{dybe}
\ena
where $\pi_V(h)^{(1)}=\pi_V(h)\otimes 1\otimes 1$ and 
$\pi_V(h)^{(1)}_{j,l}=\pi_V(h_{\bep_j})^{(1)}-\pi_V(h_{\bep_l})^{(1)}$ with $\pi_V(h_{\bep_j})=E_{jj}-\frac{1}{N}{\rm I}$ etc.. 
Here 
${\rm I}$ denotes the $N\times N$ unit matrix.
\end{prop}

\noindent
{\it Remark 1.} The elliptic  dynamical $R$-matrix \eqref{tRmat} is gauge equivalent to the $A^{(1)}_{N-1}$ type face weight obtained by Jimbo, Miwa and Okado\cite{JMO}. 
See Appendix \ref{JMOsW}.

\noindent
{\it Remark 2.}\lb{propR} The $R$-matrix preserves the weights
\bea
&&[R^\pm(z,s), \pi_V(h)\otimes \pi_V(h)]=0 \qquad \forall h\in \bar{\hh}.\lb{weightp}
\ena

Now let us set 
\begin{eqnarray}
&&\rho_0(z)=q^{-\frac{N-1}{N}}\frac{(q^2z;q^{2N})_\infty(q^{2N-2}z;q^{2N})_\infty}{(z;q^{2N})_\infty(q^{2N}z;q^{2N})_\infty},\\
&&\alpha(z)=\frac{\{pq^2z\}
\{pq^{2N-2}z\}
\{p/z\} \{pq^{2N}/z\}
}{
\{pz\}
\{pq^{2N}z\}
\{pq^2/z\} \{pq^{2N-2}/z\}
}\\
&&\Xi_p(z)=(pz;p)_\infty(p/z;p)_\infty.
\ena
Then we have 
\bea
&&\trho^\pm(z)\tb_{}(z,s)=q^{\pm1}\rho_0(z^{\pm1})^{\pm1}
\frac{(1-(q^{2s}q^{-2})^{\pm1})(1-(q^{2s}q^{2})^{\pm1})
}{
(1-q^{\pm 2s})^2 }
\frac{1-z^{\pm1}}{1-(q^2z)^{\pm1}}\nn\\
&&\qquad\qquad\qquad\qquad\qquad\times 
\alpha(z)\frac{\Xi_p(q^{2s}q^{-2})\Xi_p(q^{2s}q^{2})}{
\Xi_p(q^{2s})^2}\frac{\Xi_p(z)}{\Xi_p(q^{2}z)},\lb{compb}\\
&&\trho^\pm(z){\tbb}_{}(z)=q^{\pm1}\rho_0(z^{\pm1})^{\pm1}\frac{1-z^{\pm1}}{1-(q^2z)^{\pm1}}\alpha(z)\frac{\Xi_p(z)}{\Xi_p(q^{2}z)},\\
&&\trho^\pm(z)\tc_{}(z,s)=
\rho_0(z^{\pm1})^{\pm1}
\frac{1-q^{\pm 2}}{1-q^{\pm 2s} }
\frac{1-(q^{2s}z)^{\pm1}}{1-(q^2z)^{\pm1}}\alpha(z)\frac{\Xi_p(q^{2})\Xi_p(q^{2s}z)}{
\Xi_p(q^{2s})\Xi_p(q^{2}z)}.\lb{compc}
\end{eqnarray}
 In \cite{FIJKMY} a similar expression was obtained for 
Baxter\rq{}s elliptic $R$ matrix.  
In $R^\pm(z,s)$, we specify the factors $\frac{1-z^{\pm1}}{1-(q^2z)^{\pm1}}$ and $\frac{1-(q^{2s}z)^{\pm1}}{1-(q^2z)^{\pm1}}$  in \eqref{compb}-\eqref{compc} to be power series in $z^{\pm1}$. We then 
treat $R^\pm(z,s)$ as formal Laurent series in $z$ 
\bea
 &&R^\pm(z,s)=\sum_{n\in\Z}R^\pm(s)_nz^n
\ena
whose coefficients are in the ring $\FF[[p]]$ of formal power series in $p$. 
Note that $\alpha(z)$, $\frac{\Xi(z)}{\Xi(q^2z)}$ and $\frac{\Xi(q^{2s_{ij}} z)}{\Xi(q^2z)}$ are 
well defined formal Laurent series in $z$ with coefficients in $\FF[[p]]$. Then the matrices $R^\pm(z,s)$ satisfy 
\bea
&&R^\pm(s)_n\equiv 0\qquad {\rm mod}\ p^{{\rm max}(\mp n,0)}\FF[[p]]\qquad \forall n\in \Z.
\ena
In particular, at $p=0$ $R^+(z,s)$ (reps. $R^-(z,s)$) contains only non-negative (reps. non-positive) 
powers in $z$.  Explicitly we have $R^{\pm}_0(z,s)\equiv R^{\pm}(z,s)\biggl|_{p=0}$, 
\bea
R^{\pm}_0(z,s)&=&{\trho}_0^{\pm}(z^\pm)^{\pm 1}{\bR}^\pm_0(z,s),
\lb{trigDR}\\
{\bR}^\pm_0(z,s)
&=&
\sum_{j=1}^{N}
E_{jj}\otimes E_{jj}+
\sum_{1 \leq j<l \leq N}
\biggl(
\tb_{0}^\pm(z,s_{j,l})
E_{jj}
\otimes E_{ll}+
{\tbb}^\pm_{0}(z)
E_{ll}\otimes E_{jj}
\nonumber\\
&&\qquad\qquad +
\tc^\pm_{0}(z,s_{j,l})
E_{jl}\otimes E_{lj}+\bar{{\tc}}^\pm_{0}(z,s_{j,l})E_{lj}\otimes E_{jl}
\biggr),\nn
\ena
where  
\be
&&\rho_0^+(z)=\rho_0(z),\qquad \rho^-_0(z)=\rho_0(z^{-1})^{-1},\\
&&\tb^\pm_{0}(z,s)=
\frac{(1-q^{2(s-1)})(1-q^{2(s+1)})
}{
(1-q^{2s} )^2}q^{\pm1}
\frac{1-z^{\pm1}}{1-(q^2z)^{\pm1}},\lb{tcompb}\\
&&{\tbb}^\pm_{0}(z)=q^{\pm1}\frac{1-z^{\pm1}}{1-(q^2z)^{\pm1}},\\
&&\tc^\pm_{0}(z,s)=
\frac{1-q^{2s}z}{1-q^{2s} }
\frac{1-q^{\pm 2}}{1-(q^2z)^{\pm1}}\times \left\{\mmatrix{1\qquad \mbox{for} \ +\cr z^{-1}\qquad \mbox{for} \ -\cr}\right.,\lb{tcompc}\\
&&\bar{\tc}^\pm_{0}(z,s)=
\frac{1-q^{2s}z^{-1}}{1-q^{2s} }
\frac{1-q^{\pm 2}}{1-(q^2z)^{\pm1}}\times \left\{\mmatrix{z\qquad \mbox{for} \ +\cr 1\qquad \mbox{for} \ -\cr}\right..
\en
Hence one can regard the matrix element $R_0(z,s)_{ij}^{kl}$ as a formal power series in 
the (multiplicative) dynamical variables $q^{2s_{i,j}}$. The 0-th order term in $q^{2s_{i,j}}$ coincides with the corresponding component of 
the standard trigonometric $R$ matrix 
\bea 
R_0(z)&=&\rho_0(z)\bar{R}_0(z),\lb{trigR}\\
\bar{R}_0(z)&=&
\sum_{j=1}^{N}
E_{jj}\otimes E_{jj}+
\sum_{1 \leq j<l \leq N}
\biggl(
q\frac{1-z}{1-q^2z}(
E_{jj}
\otimes E_{ll}+
E_{ll}\otimes E_{jj})
\nonumber\\
&&\qquad\qquad +
\frac{1-q^2}{1-q^2z}
E_{jl}\otimes E_{lj}+z\frac{1-q^2}{1-q^2z}E_{lj}\otimes E_{jl}
\biggr),\nn
\ena
Note that one can parametrize the elliptic dynamical $R$ matrices associated with the other types of affine Lie algebras, at least $\gh=B_N^{(1)}, C_N^{(1)}, D_N^{(1)}$\cite{JMO,Konno06}, in a similar way to \eqref{compb}-\eqref{compc} so that they have the same property at $p=0$.

\section{The   Elliptic Quantum Algebras $U_{q,p}(\glnhbig)$ and $E_{q,p}(\glnhbig)$}
In this section we define two elliptic algebras $U_{q,p}(\glnh)$ and $E_{q,p}(\glnh)$ as topological algebras 
over $\FF[[p]]$.  

\subsection{$U_{q,p}(\glnh)$}

\begin{dfn}\lb{defUqpgl}
The elliptic algebra $U_{q,p}(\glnh)$ is a topological algebra over $\FF[[p]]$ generated by 
$e_{j,m}, f_{j,m}, k_{l,m}$,  
$(1\leq j\leq N-1, 1\leq l\leq N, m\in \Z)$, ${\hd}$ and the central element $q^{\pm c/2}$. 
We set 
\bea
&&e_j(z)=\sum_{m\in \Z} e_{j,m}z^{-m}
,\quad f_j(z)=\sum_{m\in \Z} f_{j,m}z^{-m},\lb{defenfngl}\\
&&k_l^+(z)=\sum_{m\in \Z_{\geq 0}} k_{l,-m}z^{m}+\sum_{m\in \Z_{> 0}} k_{l,m}p^mz^{-m},\\
&& k^-_l(z)=q^{2h_{\bep_l}}k^+_l(zp^*q^c).
\lb{defkngl}
\ena 
The  defining relations are as follows.  For $g(P), g(P+h)\in 
\FF$, 
\bea
&&g({P+h})e_j(z)=e_j(z)g({P+h}),\quad g({P})e_j(z)=e_j(z)g(P-<Q_{\al_j},P>),\lb{gegl}\\
&&g({P+h})f_j(z)=f_j(z)g(P+h-<{\al_j},P+h>),\quad g({P})f_j(z)=f_j(z)g(P),\lb{gfgl}\\
&&g({P})k^+_l(z)=k^+_l(z)g(P-<Q_{\bep_l},P>),\quad
g({P+h})k^+_l(z)=k^+_l(z)g(P+h-<Q_{\bep_l},P>),\nn\\
&&\lb{gkgl}
\\
&&[\hd, g(P+h)]=0=[\hd, g(P)],
\quad\lb{dggl}\\
&& [\hd, k^+_l(z)]=-z\frac{\partial}{\partial z}k^+_l(z),\quad [\hd, e_j(z)]=-z\frac{\partial}{\partial z}e_j(z), \quad  [\hd, f_j(z)]=-z\frac{\partial}{\partial z}f_j(z),\quad \lb{dedfdkgl}\\
&&\rho^+_+(z_2/z_1)k^+_{l}(z_1)k^+_{l}(z_2)=\rho^+_+(z_1/z_2)k^+_{l}(z_2)k^+_{l}(z_1),
\qquad  (1\leq l\leq N),\lb{kjkj}\\
&&\rho^+_+(z_2/z_1)\frac{(p^*z_2/z_1;p^*)_\infty(pq^2z_2/z_1;p)_\infty}{(p^*q^2z_2/z_1;p^*)_\infty(pz_2/z_1;p)_\infty}k^+_{j}(z_1)k^+_{l}(z_2)\nn\\
&&\qquad ={\rho}^+_+(z_1/z_2)\frac{(q^{-2}z_1/z_2;p^*)_\infty(z_1/z_z;p)_\infty}{(z_1/z_2;p^*)_\infty(q^{-2}z_1/z_2;p)_\infty} k^+_{l}(z_2)k^+_{j}(z_1)
  \quad (1\leq j<l\leq N),\lb{kjkl}
\ena
\bea
&&\frac{(p^*q^{c+2-j}z_2/z_1;p^*)_\infty}{(p^*q^{c-j}z_2/z_1;p^*)_\infty}k_j^{+}(z_1)e_j(z_2)=q^{-1}\frac{(q^{-c+j}z_1/z_2;p^*)_\infty}{(q^{-c-2+j}z_1/z_2;p^*)_\infty}e_j(z_2)k_j^{+}(z_1),\lb{kjej}
\\
&&\frac{(p^*q^{c-2-j}z_2/z_1;p^*)_\infty}{(p^*q^{c-j}z_2/z_1;p^*)_\infty}k_{j+1}^{+}(z_1)e_j(z_2)=q\frac{(q^{-c+j}z_1/z_2;p^*)_\infty}{(q^{-c+2+j}z_1/z_2;p^*)_\infty}e_j(z_2) k_{j+1}^{+}(z_1)
,
\lb{kjp1ej}\\
&&k_l^{+}(z_1)e_j(z_2)k_l^{+}(z_1)^{-1}=e_j(z_2)\qquad\qquad (l\not=j,j+1),\lb{ejkl}
\\
&&\frac{(pq^{-j}z_2/z_1;p)_\infty}{(pq^{2-j}z_2/z_1;p)_\infty}k_j^{+}(z_1)f_j(z_2)=q\frac{(q^{-2+j}z_1/z_2;p)_\infty}{(q^{j}z_1/z_2;p)_\infty}f_j(z_2)k_j^{+}(z_1),
\\
&&\frac{(pq^{-j}z_2/z_1;p)_\infty}{(pq^{-2-j}z_2/z_1;p)_\infty}k_{j+1}^{+}(z_1)f_j(z_2)=q^{-1}\frac{(q^{2+j}z_1/z_2;p)_\infty}{(q^{j}z_1/z_2;p)_\infty}f_j(z_2)k_{j+1}^{+}(z_1)
,\lb{fjkjp1}\\
&&k_l^{+}(z_1)f_j(z_2)k_l^{+}(z_1)^{-1}=f_j(z_2)\qquad\qquad (l\not=j,j+1),
\lb{fjkl}
\ena
\bea
&&
z_1 \frac{(q^{2}z_2/z_1;p^*)_\infty}{(p^*q^{-2}z_2/z_1;p^*)_\infty}e_j(z_1)e_j(z_2)=
-z_2 \frac{(q^{2}z_1/z_2;p^*)_\infty}{(p^*q^{-2}z_1/z_2;p^*)_\infty}e_j(z_2)e_j(z_1),\lb{ejejgl}\\
&&z_1 \frac{(q^{-1}z_2/z_1;p^*)_\infty}{(p^*qz_2/z_1;p^*)_\infty}e_j(z_1)e_{j+1}(z_2)=
-z_2 \frac{(q^{-1}z_1/z_2;p^*)_\infty}{(p^*qz_1/z_2;p^*)_\infty}e_{j+1}(z_2)e_j(z_1),\lb{ejejp1gl}\\
&&e_j(z_1)e_l(z_2)=e_l(z_2)e_j(z_1)\qquad\qquad (|j-l|>1)\\
&&
z_1 \frac{(q^{-2}z_2/z_1;p)_\infty}{(pq^{2}z_2/z_1;p)_\infty}f_j(z_1)f_j(z_2)=
-z_2 \frac{(q^{-2}z_1/z_2;p)_\infty}{(pq^{2}z_1/z_2;p)_\infty}f_j(z_2)f_j(z_1),\lb{fjfjgl}\\
&&
z_1 \frac{(qz_2/z_1;p)_\infty}{(pq^{-1}z_2/z_1;p)_\infty}f_j(z_1)f_{j+1}(z_2)=
-z_2 \frac{(qz_1/z_2;p)_\infty}{(pq^{-1}z_1/z_2;p)_\infty}f_{j+1}(z_2)f_j(z_1),\lb{fjfjp1gl}\\
&&f_j(z_1)f_l(z_2)=f_l(z_2)f_j(z_1)\qquad\qquad (|j-l|>1)\lb{fjfl}
\\
&&[e_i(z_1),f_j(z_2)]=\frac{\delta_{i,j}\kappa}{q-q^{-1}}
\left(\delta(
q^{-c}
z_1/z_2)
k_j^-(
q^{-\frac{c}{2}}
z_1)k_{j+1}^-(
q^{-\frac{c}{2}}
z_1)^{-1}\right.\nn\\
&&\left.\qquad\qquad\qquad\qquad\qquad\qquad\qquad -
\delta(
q^c
z_1/z_2)
k_j^+(
q^{-\frac{c}{2}}
z_2)k_{j+1}^+(
q^{-\frac{c}{2}}
z_2)^{-1}
\right),\lb{eifj}
\ena
\bea
&&
\frac{(p^*q^{2}{z_{2}}/{z_{1}}; p^*)_{\infty}}
{(p^*q^{-2}{z_{2}}/{z_{1}}; p^*)_{\infty}}
\left\{\frac{(p^*q^{-1}{z_{1}}/{w}; p^*)_{\infty}}
{(p^*qz_{1}/{w}; p^*)_{\infty}}
\frac{(p^*q^{-1}{z_{2}}/w; p^*)_{\infty}}
{(p^*q{z_{2}}/w; p^*)_{\infty}}
e_j(w)e_i(z_1)e_i(z_2)\right.
\nn\\
&&\left.\qquad\qquad\qquad-[2]_q\frac{(p^*q^{-1}w/{z_{1}}; p^*)_{\infty}}
{(p^*qw/z_{1}; p^*)_{\infty}}
\frac{(p^*q^{-1}{z_{2}}/w; p^*)_{\infty}}
{(p^*q{z_{2}}/w; p^*)_{\infty}}
e_i(z_1)e_j(w)e_i(z_2)\right.
\nn\\
&&\left.
\qquad+\frac{(p^*q^{-1}w/{z_{1}}; p^*)_{\infty}}
{(p^*qw/z_{1}; p^*)_{\infty}}
\frac{(p^*q^{-1}w/{z_{2}}; p^*)_{\infty}}
{(p^*qw/{z_{2}}; p^*)_{\infty}}
e_i(z_1)e_i(z_2)e_j(w)\right\}+(z_1\leftrightarrow z_2)=0,\label{serreegl}
\ena
\bea
&&
\frac{(pq^{-2}{z_{2}}/{z_{1}}; p)_{\infty}}
{(pq^{2}{z_{2}}/{z_{1}}; p)_{\infty}}
\left\{\frac{(pq^{}{z_{1}}/{w}; p)_{\infty}}
{(pq^{-1}z_{1}/{w}; p)_{\infty}}
\frac{(pq^{}{z_{2}}/w; p)_{\infty}}
{(pq^{-1}{z_{2}}/w; p)_{\infty}}
f_j(w)f_i(z_1)f_i(z_2)\right.
\nn\\
&&\left.\qquad\qquad\qquad-[2]_q\frac{(pq^{}w/{z_{1}}; p)_{\infty}}
{(pq^{-1}w/z_{1}; p)_{\infty}}
\frac{(pq^{}{z_{2}}/w; p)_{\infty}}
{(pq^{-1}{z_{2}}/w; p)_{\infty}}
f_i(z_1)f_j(w)f_i(z_2)\right.
\nn\\
&&\left.
\qquad+\frac{(pq^{}w/{z_{1}}; p)_{\infty}}
{(pq^{-1}w/z_{1}; p)_{\infty}}
\frac{(pq^{}w/{z_{2}}; p)_{\infty}}
{(pq^{-1}w/{z_{2}}; p)_{\infty}}
f_i(z_1)f_i(z_2)f_j(w)\right\}+(z_1\leftrightarrow z_2)=0
\quad   |i-j|=1.\nn\\&&
\label{serrefgl}
\end{eqnarray}  
where $\delta(z)=\sum_{n\in \Z}z^n$, $\rho(z)$ is given in \eqref{rho}, 
\bea
&&\rho^+_+(z)= \frac{\{q^2z\}^*\{q^{-2}q^{2N}z\}^*\{z\}\{q^{2N}z\}}{\{z\}^*\{q^{2N}z\}^*\{q^2z\}\{q^{-2}q^{2N}z\}},
\ena
and
$\kappa$ is given by
\bea
&&\kappa=\frac{(p;p)_\infty(p^*q^2;p^*)_\infty}{(p^*;p^*)_\infty(pq^2;p)_\infty}.
\ena   
We call $e_j(z), f_j(z), k^\pm_l(z)$ the elliptic currents. 
We also denote by ${U}'_{q,p}(\glnh)$ the subalgebra obtained by removing ${\hd}$. 
\end{dfn}
We treat these relations 
 as formal Laurent series in $z, w$ and $z_j$'s. 
In  each term of \eqref{kjkl}-\eqref{fjfl} and \eqref{serreegl}-\eqref{serrefgl}, the expansion direction of the structure function given by a ratio of  infinite products is 
chosen according to the order of the accompanied product of the
 elliptic currents. For example, in the l.h.s of \eqref{ejejgl}, $\frac{(q^{2}z_2/z_1;p^*)_\infty}{(p^*q^{-2}z_2/z_1;p^*)_\infty}$ should be  expanded in $z_2/z_1$, 
 whereas in the r.h.s $\frac{(q^{2}z_1/z_2;p^*)_\infty}{(p^*q^{-2}z_1/z_2;p^*)_\infty}$ should be expanded in $z_1/z_2$.   
All the coefficients in $z_j$'s are well defined in the $p$-adic topology.  


For a practical use, we remark that  in the sense of analytic continuation \eqref{kjkj}-\eqref{fjkjp1} and \eqref{ejejgl}-\eqref{fjfjp1gl} 
can be rewritten as follows.  
\bea
&&k^+_{l}(z_1)k^+_{l}(z_2)={\rho}(z_1/z_2)k^+_{l}(z_2)k^+_{l}(z_1),
\qquad  (1\leq l\leq N),\lb{kjkja}\\
&&k^+_{j}(z_1)k^+_{l}(z_2)={\rho}(z_1/z_2)\frac{\Theta_{p^*}(q^{-2}z_1/z_2)\Theta_{p}(z_1/z_2)}{\Theta_{p^*}(z_1/z_2)\Theta_{p}(q^{-2}z_1/z_2)} k^+_{l}(z_2)k^+_{j}(z_1)  \quad (1\leq j<l\leq N),\lb{kjkla}\\
&&k_j^{+}(z_1)e_j(z_2)k_j^{+}(z_1)^{-1}=q^{-1}\frac{\Theta_{p^*}(q^{-c+j}z_1/z_2)}{\Theta_{p^*}(q^{-c-2+j}z_1/z_2)}e_j(z_2),\lb{kjeja}
\\
&&k_{j+1}^{+}(z_1)e_j(z_2)k_{j+1}^{+}(z_1)^{-1}=q^{}\frac{\Theta_{p^*}(q^{-c+j}z_1/z_2)}{\Theta_{p^*}(q^{-c+2+j}z_1/z_2)}e_j(z_2) \qquad (1\leq j\leq N-1),\lb{kjp1eja}
\\
&&k_j^{+}(z_1)f_j(z_2)k_j^{+}(z_1)^{-1}=q^{}\frac{\Theta_{p}(q^{-2+j}z_1/z_2)}{\Theta_{p}(q^{j}z_1/z_2)}f_j(z_2),\lb{kjfja}
\\
&&k_{j+1}^{+}(z_1)f_j(z_2)k_{j+1}^{+}(z_1)^{-1}=q^{-1}\frac{\Theta_{p}(q^{2+j}z_1/z_2)}{\Theta_{p}(q^{j}z_1/z_2)}f_j(z_2) \qquad (1\leq j\leq N-1),
\lb{kjp1fja}\\
&&e_j(z_1)e_j(z_2)=
-\frac{z_2}{z_1 } \frac{\Theta_{p^*}(q^{2}z_1/z_2)}{\Theta_{p^*}(q^{2}z_2/z_1)}e_j(z_2)e_j(z_1),\lb{eethetajj}\\
&&e_j(z_1)e_{j+1}(z_2)=
-\frac{z_2}{z_1 } \frac{\Theta_{p^*}(q^{-1}z_1/z_2)}{\Theta_{p^*}(q^{-1}z_2/z_1)}e_{j+1}(z_2)e_j(z_1),\lb{eethetajjp1}\\
&&f_j(z_1)f_j(z_2)=
-\frac{z_2}{z_1 } \frac{\Theta_{p^*}(q^{-2}z_1/z_2)}{\Theta_{p}(q^{-2}z_2/z_1)}f_j(z_2)f_j(z_1),\lb{ffthetajj}\\
&&f_j(z_1)f_{j+1}(z_2)=
-\frac{z_2}{z_1 } \frac{\Theta_{p^*}(qz_1/z_2)}{\Theta_{p}(qz_2/z_1)}f_{j+1}(z_2)f_j(z_1)
.\lb{ffthetajjp1}
\ena

\begin{prop}\lb{CenterK}
Let us set
\be
&&K(z)=k^+_1(z)k^+_2(q^{-2}z)\cdots k^+_N(q^{-2(N-1)}z).
\en
Then $K(z)$ belongs to the center of $U_{q,p}\rq{}(\glnh)$. 
\end{prop}
\noindent
{\it Proof.} Direct calculation using \eqref{gkgl},\eqref{ejkl},\eqref{fjkl},  \eqref{kjkja}-\eqref{kjp1fja} shows that $K(z)$ commutes with 
$\FF$ and all of the elliptic currents of $U_{q,p}(\glnh)$. In particular, 
$[K(z),k^+_l(w)]=0\ (1\leq l \leq N)$ follows from the  
 identity
\be
&&\prod_{j=1}^N\rho(q^{-2(j-1)}z)=\frac{\Theta_{p^*}(z)\Theta_{p}(q^{-(N-1)}z)}{\Theta_{p^*}(q^{-(N-1)}z)\Theta_{p}(z)}.
\en
\qed

\noindent
{\it Remark.}\ In Appendix \ref{ellqdet} we identify $K(z)$ with the $q$-determinant of the $L$-operator. 

The elliptic algebra $U'_{q,p}(\slnh)$ is identified with the quotient algebra 
$U'_{q,p}(\glnh)/<K(z)-1>$.
More explicitly, one can realize $k^+_l(z) \ (1\leq l\leq N)$ satisfying \eqref{gkgl}, 
\eqref{kjkj}-\eqref{kjkl} and $K(z)=1$ as follows. 
  Let  $A=(a_{ij})_{i,j\in I\cup\{0\}}$ be the $A^{(1)}_{N-1}$ type generalized Cartan matrix. 
Let $\al_{i,m}\ (i\in I, m\in \Z_{\not=0})$ be the Heisenberg algebra satisfying 
\bea
&&[\al_{i,m},\al_{j,n}]=\delta_{m+n,0}\frac{[a_{ij}m]_q
[cm]_q
}{m}
\frac{1-p^m}{1-p^{*m}}
q^{-cm}.\lb{ellboson}
\ena
Let us consider the following $\cE^{+l}_m\ (1\leq l\leq N, m\in \Z_{\not=0})$,  
which we call the elliptic bosons of the orthonormal basis type\cite{FKO}.
\be
&&\cE^{+l}_m=\frac{q^{lm}}{(q-q^{-1})[m]_q^2[Nm]_q}\left(-q^{-Nm}\sum_{k=1}^{l-1}[km]_q\al_{k,m}+\sum_{k=l}^{N-1}[(N-k)m]_q\al_{k,m}\right)\quad (1\leq l\leq N-1),\\
&&\cE^{+N}_m=-\frac{1}{(q-q^{-1})[m]_q^2[Nm]_q}\sum_{k=1}^{N}[km]_q\al_{k,m}.
\en
They satisfy
\bea
&&[\cE^{+ l}_m,\cE^{+ l}_n]
=\delta_{m+n,0}\frac{[cm]_q[(N-1) m]_q }{m(q-q^{-1})^2[m]_q^3 [N m]_q}\frac{1-p^m}{1-p^{*m}}q^{-cm}, \\
&&[\cE^{+j}_m,\cE^{+ l}_n]
=-\delta_{m+n,0}q^{({\rm sgn}(l-j)N -l+j)m}
\frac{[cm]_q}{m(q-q^{-1})^2[m]_q^2[Nm]_q}\frac{1-p^m}{1-p^{*m}}q^{-cm},\\
&&[\al_{i,m}, \cE^{+ l}_n]= \delta_{m+n,0}\frac{[cm]_q}{m(q^m-q^{-m})}\frac{1-p^m}{1-p^{*m}}q^{-cm}(q^{- m}\delta_{i,l}-\delta_{i,l-1}).
\ena
Let  $K^+_{\bep_j}$ satisfies for $g(P), g(P+h)\in \FF$
\bea
&&g({P})K^+_{\bep_j}=K^+_{\bep_j}g(P-<Q_{\bep_j},P>),\\
&&\ g({P+h})K^+_{\bep_j}=K^+_{\bep_j}g(P+h-<Q_{\bep_j},P>).\lb{gKpm}
\ena
Then the following $k^+_j(z)$ satisfy the desired relations.
\bea
k^+_l(z)&=&K^+_{\bep_l} :\exp\left\{\sum_{m\not=0}\frac{(q^m-q^{-m})^2p^m}{1-p^m}\cE^{+l}_m (q^lz)^{-m}\right\}:.
\ena
Furthermore if we further require that $\al_{i,m}$ and $K_{\bep_i}^+$ satisfy
\bea
&&[g(P), \al_{i,m}]=[g(P+h),\al_{i,n}]=0,\lb{gboson}\\ 
&& [\hd, \al_{j,n}]=n\al_{j,n},\\
&&
[\al_{i,m},e_j(z)]=\frac{[a_{ij}m]_q}{m}\frac{1-p^m}{1-p^{*m}}
q^{-cm}z^m e_j(z),
\lb{bosonve}\\
&&
[\al_{i,m},f_j(z)]=-\frac{[a_{ij}m]_q}{m}z^m f_j(z)
,\lb{bosonvf}\\
&&K_{\bep_i}^{+}e_j(z)=q^{- <\al_j,h_{\bep_i}>}e_j(z)K_{\bep_i}^{+},\quad 
K_{\bep_i}^{+}f_j(z)=q^{ <\al_j,h_{\bep_i}>}f_j(z)K_{\bep_i}^{+}, 
\ena
$k^+_l(z)$ satisfy the remaining relations \eqref{dedfdkgl}, \eqref{kjej}-\eqref{fjkl}.

Now let us define $\psi^\pm_j(z)\ (1\leq j\leq N-1)$ by\footnote{Our $\psi^\pm_j(z)$ are $\psi^\mp_j(z)$ in \cite{JKOS,KK03}.}
\be
&&\psi^+_j(q^{-c/2}q^jz)=\kappa k^+_j(z)k^+_{j+1}(z)^{-1},\\
&&\psi^-_j(q^{-c/2}q^jz)=\kappa k^-_j(z)k^-_{j+1}(z)^{-1}.
\en
We have \bea
&&{\psi}_j^+(
q^{-\frac{c}{2}}
z)=K^+_{j}\exp\left(-(q-q^{-1})\sum_{n>0}\frac{\al_{j,-n}}{1-p^n}z^n\right)
\exp\left((q-q^{-1})\sum_{n>0}\frac{p^n\al_{j,n}}{1-p^n}z^{-n}\right),\lb{psip}
\ena
 and $\psi^-_j(z)=q^{2h_j}\psi^+_j(zpq^{-c})$
where we set $K^+_j=K^+_{\bep_j}K^{+-1}_{\bep_{j+1}}$.
\begin{prop}\lb{defUqp}
The elliptic algebra $U_{q,p}(\slnh)$ is characterized by \eqref{gboson}-\eqref{bosonvf} and the following relations.  For $g(P), g(P+h)\in 
\FF$, 
\bea
&&g({P+h})e_j(z)=e_j(z)g({P+h}),\quad g({P})e_j(z)=e_j(z)g(P-<Q_{\al_j},P>),\lb{gesln}\\
&&g({P+h})f_j(z)=f_j(z)g(P+h-<{\al_j},P+h>),\quad g({P})f_j(z)=f_j(z)g(P),\lb{gfsln}\\
&&[\hd, g(P+h,P)]=0,
\quad\lb{dgsln}\\
&&[\hd, e_j(z)]=-z\frac{\partial}{\partial z}e_j(z), \quad  [\hd, f_j(z)]=-z\frac{\partial}{\partial z}f_j(z),\quad \lb{dedfsln}\\
&&
z_1 \frac{(q^{a_{ij}}z_2/z_1;p^*)_\infty}{(p^*q^{-a_{ij}}z_2/z_1;p^*)_\infty}e_i(z_1)e_j(z_2)=
-z_2 \frac{(q^{a_{ij}}z_1/z_2;p^*)_\infty}{(p^*q^{-a_{ij}}z_1/z_2;p^*)_\infty}e_j(z_2)e_i(z_1),\lb{eesln}\\
&&
z_1 \frac{(q^{-b_{ij}}z_2/z_1;p)_\infty}{(pq^{a_{ij}}z_2/z_1;p)_\infty}f_i(z_1)f_j(z_2)=
-z_2 \frac{(q^{-a_{ij}}z_1/z_2;p)_\infty}{(pq^{a_{ij}}z_1/z_2;p)_\infty}f_j(z_2)f_i(z_1),\lb{ffsln}\\
&&[e_i(z_1),f_j(z_2)]=\frac{\delta_{i,j}}{q-q^{-1}}
\left(\delta(
q^{-c}
z_1/z_2)
\psi_j^-(
q^{\frac{c}{2}}
z_2)-
\delta(
q^c
z_1/z_2)
\psi_j^+(
q^{-\frac{c}{2}}
z_2)
\right),\lb{eifjsln}
\ena
\bea
&&\sum_{\sigma\in S_{a}}\prod_{1\leq m< k \leq a}
\frac{(p^*q^{2}{z_{\sigma(k)}}/{z_{\sigma(m)}}; p^*)_{\infty}}
{(p^*q^{-2}{z_{\sigma(k)}}/{z_{\sigma(m)}}; p^*)_{\infty}}\nn\\
&&\quad\times\sum_{s=0}^{a}(-1)^{s}
\left[\begin{array}{c}
a\cr
s\cr
\end{array}\right]_q\prod_{1\leq m \leq s}\frac{(p^*q^{a_{ij}}{w}/{z_{\sigma(m)}}; p^*)_{\infty}}
{(p^*q^{-a_{ij}}{w}/{z_{\sigma(m)}}; p^*)_{\infty}}
\prod_{s+1\leq m \leq a}\frac{(p^*q^{a_{ij}}{z_{\sigma(m)}}/{w}; p^*)_{\infty}}
{(p^*q^{-a_{ij}}{z_{\sigma(m)}}/{w}; p^*)_{\infty}}
\nn\\
&&\quad\times
e_{i}(z_{\sigma(1)})\cdots e_{i}(z_{\sigma(s)})e_{j}(w)e_{i}(z_{\sigma(s+1)}) \cdots 
e_{i}(z_{\sigma(a)})=0,\label{serree}
\ena
\bea
&&\sum_{\sigma\in S_{a}}\prod_{1\leq m< k \leq a
}
\frac{(pq^{-2}{z_{\sigma(k)}}/{z_{\sigma(m)}}; p)_{\infty}}
{(pq^{2}{z_{\sigma(k)}}/{z_{\sigma(m)}}; p)_{\infty}}\nn\\
&&\quad\times\sum_{s=0}^{a}(-1)^{s}\left[\begin{array}{c}
a\cr
s\cr
\end{array}\right]_q\prod_{1\leq  m \leq s}
\frac{(pq^{-a_{ij}}{w}/{z_\sigma(m)}; p)_{\infty}}{(pq^{a_{ij}}{w}/{z_{\sigma(m)}}; p)_{\infty}}
\prod_{s+1\leq m \leq a}\frac{(pq^{-a_{ij}}{z_{\sigma(m)}}/{w}; p)_{\infty}}{(pq^{a_{ij}}{z_{\sigma(m)}}/{w}; p)_{\infty}}\nn\\
&&\times f_{i}(z_{\sigma(1)})\cdots f_{i}(z_{\sigma(s)})f_{j}(w)f_{i}(z_{\sigma(s+1)}) \cdots f_{i}(z_{\sigma(a)})=0
\quad(i\neq j, a=1-a_{ij}),\label{serref}
\end{eqnarray}
\end{prop}

\begin{prop}
In the sense of analytic continuation, we have 
\bea
&&\psi^+_i(z_1)\psi^+_j(z_2)=\frac{\Theta_{p^*}(q^{a_{ij}}z_1/z_2)
\Theta_{p}(q^{-a_{ij}}z_1/z_2)} {\Theta_{p^*}(q^{-a_{ij}}z_1/z_2)
\Theta_{p}(q^{a_{ij}}z_1/z_2)}\psi^+_j(z_2) \psi^+_i(z_1),\\  
&&\psi^+_i(z_1)e_j(z_2)=q^{-a_{ij}}\frac{\Theta_{p^*}(q^{a_{ij}-c/2}z_1/z_2)} {\Theta_{p^*}(q^{-a_{ij}-c/2}z_1/z_2)}e_j(z_2) \psi^+_i(z_1),\\
&&\psi^+_i(z_1)f_j(z_2)=q^{a_{ij}}\frac{
\Theta_{p}(q^{-a_{ij}+c/2}z_1/z_2)} {\Theta_{p}(q^{a_{ij}+c/2}z_1/z_2)}f_j(z_2) \psi^+_i(z_1).
\ena
\end{prop}

Let $U_q(\gh)$ be the quantum affine algebra over $\C$ associated with the untwisted affine Lie algebra $\gh$  in the Drinfeld realization\cite{Dr} 
and $x^\pm_j(z), k^\pm_{0,l}(z)$ be the Drinfeld currents. 
See Appendix \ref{Uqgh} for the $\glnh$ case. The other cases can be found, for example in \cite{FKO}.    
Then $U_{q,p}(\gh)$ is a natural face type ( i.e. dynamical) 
elliptic deformation of $U_q(\gh)$ in the following sense. 
\begin{thm}\cite{FKO}
\be
&&{U_{q,p}(\gh)/pU_{q,p}(\gh)\cong (\FF\otimes_{\C}U_q(\gh))\sharp \C[\cR_Q]}
\en
by the following identification at $p=0$.
\be
&&e_j(z)=x^+_j(z)e^{-Q_{\al_j}},\quad f_j(z)=x^-_j(z), \quad k^\pm_l(z)=k^\pm_{0,l}(z)e^{-Q_{\bep_l}}. 
\en
 Here the smash product $\sharp$ is defined as follows.
\be
&&g(P,P+h)a\otimes e^{Q_\al} \cdot f(P,P+h)b\otimes e^{Q_\beta}\\
&&\quad= g(P,P+h)f(P-<Q_\al,P>,P+h-<Q_\al+\wt( a),P+h>)ab\otimes e^{Q_\al+Q_\beta}
\en
where $\wt({a})\in \bar{\hh}^*$ s.t. $q^h a q^{-h}=q^{<\wt({a}),h>}a$ for $a, b\in U_q(\gh), f(P), g(P)\in \FF, e^{Q_\al}, e^{Q_\beta}\in \C[\cR_Q]$. \end{thm}
\begin{dfn}
Let us introduce the multiplicative dynamical parameters $x=(x_1,\cdots,x_N), x_i=q^{2P_{\bep_i}}$.
We set ${U}_{q,x}(\gh)=U_{q,p}(\gh)/pU_{q,p}(\gh)$ and 
call  it the dynamical quantum affine algebra in the Drinfeld realization.  
\end{dfn}

\subsection{$E_{q,p}(\glnh)$}\lb{secEqp}
Let $\bar{L}_{ij, n}\ (n\in \Z, 1\leq i,j\leq N) $ be abstract symbols.   
We define $L^+(z)=\sum_{1\leq i,j\leq N}E_{ij}L^+_{ij}(z)$ by  
\bea
&&L^+_{ij}(z)= \sum_{n\in \Z} L_{ij, n} z^{-n},\qquad L_{ij, n}=p^{{\rm max}( n,0)} 
\bar{L}_{ij, n}.\lb{genL}
\ena
\begin{dfn}\lb{defEqp}
Let $R^{+}(z,s)$ be the same $R$ matrix as in Sec.\ref{edR}. 
The elliptic algebra  $E_{q,p}(\glnh)$  is a topological algebra over $\FF[[p]]$ generated by $\bar{L}_{ij, n}$, 
$
\hd$  
 and the central element $q^{\pm c/2}$ satisfying the following relations. 
\bea
&&R^{+(12)}(z_1/z_2,P+h)L^{+(1)}(z_1)L^{+(2)}(z_2)=L^{+(2)}(z_2)L^{+(1)}(z_1)R^{+*(12)}(z_1/z_2,P),\lb{RLL}\\
&&g({P+h})\bar{L}_{ij,n}=\bar{L}_{ij,n}\; g(P+h-<Q_{\bep_i},P+h>),\lb{lgr}\\
&&g({P})\bar{L}_{ij,n}=\bar{L}_{ij,n}\; g(P-<Q_{\bep_j},P>),\lb{rgr}\\
&&[\hd, L^+(z)]=-z\frac{ \partial}{\partial z}L^+(z),
\ena
where $g({P+h}), g(P)\in \FF$ and 
\be
L^{+(1)}(z)=L^{+}(z)\otimes \id, \qquad L^{+(2)}(z)=\id \otimes L^{+}(z). 
\en
\end{dfn}
We regard $L^+(z) \in \End V \otimes E_{q,p}(\glnh)$.
We treat \eqref{RLL} as a formal Laurent series in $z_1$ and $z_2$. Then the coefficients of $z_1, z_2$ are well defined in the $p$-adic topology.  See \cite{FIJKMY} for a similar formulation for the vertex type elliptic quantum algebra $\Aqp(\slth)$. 
Note also that due to the $RLL$-relation \eqref{RLL} the $L$-operator $L^+(z)$ is invertible. See Appendix \ref{ellqdet}. 

For later convenience we define $L^-(z)=\sum_{1\leq i,j\leq N}E_{ij}L^-_{ij}(z)$ by\cite{JKOStg}   
\bea
&&L^-(z)=\left(\Ad(q^{-2{\theta}_V(P)})\otimes \id\right)\left(q^{2{T}_{V}}L^+(z p^* q^{c})\right)
,\lb{LmfromLp}\\
&&{\theta}_V(P)=-\sum_{j=1}^{N-1}\left(\frac{1}{2}\pi_V(h_j)\pi_V(h^j)+P_{j}\pi_V(h^j)\right),\lb{thetaV}\\
&&{T}_{V}
=\sum_{j=1}^{N-1}\pi_V(h_j)\otimes h^j.\lb{TV}
\ena
Here $(\Ad X)Y=XYX^{-1}$, $h^j=h_{\bar{\Lambda}_j}\ (j\in I)$
, $\pi_V(h_j)=E_{jj}-E_{j+1j+1}$ 
and 
$\pi_V(h^j)=\sum_{i=1}^j \pi_V(h_{\bep_i})\ (j\in I)$. 
Then one can verify the following.  
\begin{prop}\lb{derRLL}
The $L$ operators $L^+(z)$ and $L^-(z)$ satisfy the following relations.  
\bea
&&\hspace{-1cm}R^{-(12)}(z_1/z_2,P+h)L^{-(1)}(z_1)L^{-(2)}(z_2)=L^{-(2)}(z_2)L^{-(1)}(z_1)R^{-*(12)}(z_1/z_2,P),\lb{mRLL}\\
&&\hspace{-1cm}R^{\pm(12)}(q^{\pm{c}}z_1/z_2,P+h)L^{\pm(1)}(z_1)L^{\mp(2)}(z_2)=L^{\mp(2)}(z_2)L^{\pm(1)}(z_1)R^{\pm*(12)}(q^{\mp{c}}z_1/z_2,P).
\lb{RLLpm}
\ena
\end{prop}  

\noindent
{\it proof)} Replace $z_i$ with $z_ip^*q^{c}\ (i=1,2)$ in \eqref{RLL}. 
Note  that \eqref{lgr}, \eqref{rgr} and \eqref{LmfromLp} yields 
\be
L^+(p^*q^cz)=
q^{2\frac{N-1}{N}}
\sum_{i,j}q^{-2(P+h)_{\bep_i}}q^{2P_{\bep_j}}E_{ij}L^-_{ij}(z). 
\en 
By a componentwise comparison we obtain
\be
&&R^+(z_1/z_2,P+h)L^-(z_1)L^-(z_2)=L^-(z_2)L^-(z_1)R^{+*}(z_1/z_2,P).
\en
Then noting \eqref{trhopm} and \eqref{trhorhos}, we obtain \eqref{mRLL}.

Similarly let us replace $z_1$ by $z_1p^*q^{c}$ in \eqref{RLL}. 
Noting  $p^*q^c=pq^{-c}$, the components of $R^+$ are changed as 
\bea
&&\trho^+(z p q^{-c})=q^{-2\frac{N-1}{N}}\trho^-(z q^{-c}), \nn\\
&&\tb(z p q^{-c},s)=q^2\tb(z q^{-c},s), \qquad \tbb(z p q^{-c})=q^2\tbb(z q^{-c}),\lb{shiftR}\\
&&\tc(z p q^{-c},\pm s)=q^{\mp 2s +2}\tc(z q^{-c},\pm s) \nn
\ena
and similarly for $R^{+*}$. Then from \eqref{tRmat} and \eqref{trhopm}, 
we obtain
the second (lower sign) relation in \eqref{RLLpm}.  
Note that a factor arising from the action of $\Ad(q^{-2{\theta}_V(P)})\otimes \id$ on the $L$-operators cancels the extra factors in \eqref{shiftR}.  

To obtain the first relation in \eqref{RLLpm}, exchange $z_1$ and $z_2$ in the second relation of \eqref{RLLpm}. Then  we have
\be
&&R^-(q^{-c}z_2/z_1,P+h)^{-1}L^+(z_1)L^-(z_2)=L^-(z_2)L^+(z_1)R^{-*}(q^cz_2/z_1,P)^{-1}.
\en  
Using \eqref{tRpRm}, we obtain the desired result. 
{\qed}

\noindent
{\it Remark. } \   
We can expand  \eqref{RLL} and \eqref{mRLL} in  both $z=z_1/z_2$ and $z^{-1}=z_2/z_1$. 
However  \eqref{RLLpm} admits an expansion only in  
 $z$ (resp. $z^{-1}$) for the upper (resp. lower) sign case 
 for the sake of the well-definedness in the $p$-adic topology. 
 It is instructive to compare this with the trigonometric case \cite{DF}.

In the component form \eqref{RLL}, \eqref{mRLL} and \eqref{RLLpm} are 
\bea
 &&\hspace{-1cm}\sum_{i',j'}R^\pm(z_1/z_2,P+h)_{ij}^{i'j'}L^\pm_{i'i''}(z_1)L^\pm_{j'j''}(z_2)
= \sum_{i',j'}L^\pm_{jj'}(z_2)L^\pm_{ii'}(z_1)R^{\pm*}(z_1/z_2,P)_{i'j'}^{i''j''},
\lb{hRLLcom}\\
&&\hspace{-1.5cm}\sum_{i',j'}R^\pm(q^{\pm c}z_1/z_2,P+h)_{ij}^{i'j'}L^\pm_{i'i''}(z_1)L^\mp_{j'j''}(z_2)
= \sum_{i',j'}L^\mp_{jj'}(z_2)L^\pm_{ii'}(z_1)R^{\pm*}(q^{\mp c}z_1/z_2,P)_{i'j'}^{i''j''},
\lb{hRLLpmcom}
\ena
We call \eqref{hRLLcom} the $(i,j), (i'',j'')$ component of \eqref{RLL}, etc. 

\noindent
{\it Remark.}\  In order to obtain a `fully' dynamical $RLL$-relations used  in \cite{Felder, FV} with a central extension 
one may introduce the $L$-operators $L^\pm(z,P)$ related to our $L^\pm(z)$ by\cite{JKOS,KK03}
\bea
&&L^\pm(z,P)=L^\pm(z)e^{\sum_{i=1}^N\pi_V(h_{\vep_i})\otimes Q_{\bep_i}},\lb{DL}
\ena
where $\pi_V(h_{\epsilon_i})=E_{i,i}$.
In fact from \eqref{lgr} and \eqref{rgr} we have  
\bea
&&[L^\pm_{ij}(z,P), f(P)]=0,\qquad \lb{LPcomm}\\
&&g(h)L^\pm_{ij}(z,P)=L^\pm_{ij}(z,P)g(h-<\alpha_{ij},h>),\lb{Lhcomm}\\
&&[\hd, L^\pm(z)]=-z\frac{ \partial}{\partial z}L^\pm(z).\lb{Ldcomm}
\ena
 \eqref{LPcomm} indicates that  $L^+(z,P)$ is independent of  
$\C[\cR_Q]$. 
Furthermore from \eqref{HA2}, \eqref{RLL} \eqref{mRLL} and \eqref{RLLpm}, $L^\pm(z,P)$ satisfy the following full dynamical $RLL$-relations 
\bea
&&R^{\pm(12)}(z_1/z_2,{P+h})L^{\pm(1)}(z_1,P)L^{\pm(2)}(z_2,{P+\pi_V(h)^{(1)}})\nn\\
&&\qquad\qquad=L^{\pm(2)}(z_2,P)L^{\pm(1)}(z_1,{P+\pi_V(h)^{(2)}})R^{\pm*(12)}(z_1/z_2,P),
\lb{DRll}\\
&&R^{\pm(12)}(q^{\pm{c}}z_1/z_2,{P+h})L^{\pm(1)}(z_1,P)L^{\mp(2)}(z_2,{P+\pi_V(h)^{(1)}})\nn\\
&&\qquad\qquad=L^{\mp(2)}(z_2,P)L^{\pm(1)}(z_1,{P+\pi_V(h)^{(2)}})R^{\pm*(12)}(q^{\mp{c}}z_1/z_2,P).
\lb{DRllpm}
\ena
Here the generators are clear. If we set  $L^\pm(z,P)=\sum_{i,j}E_{i,j}L^\pm_{ij}(z,P)$ 
with $L^\pm_{ij}(z,P)=\sum_{m\in \Z}L^\pm_{ij, n}(P)z^{- n}$, then from \eqref{genL} we have 
$L^\pm_{ij, n}(P)=L^\pm_{ij,n}e^{-Q_{\bep_j}}$.

\noindent
{\it Remark.}\ The dynamical $RLL$ relations \eqref{DRll}-\eqref{DRllpm} coincides with those derived from the 
universal DYBE for $\Bqla(\gh)$ in \cite{JKOStg,JKOS}.

\subsection{Reflection equations}
Following \cite{RS}, let us set
\be
&&\cL(z)=L^+(zq^{c})L^-(z)^{-1}.
\en
Then using \eqref{RLL}, \eqref{mRLL}-\eqref{RLLpm}, one can show the following relations.
\begin{prop}
\bea
&&R^{+(12)}(z_1/z_2,P+h)\cL^{(1)}(z_1)R^{+(21)}(q^{2c}z_2/z_1,P+h)\cL^{(2)}(z_2)\nn\\
&&\qquad\qquad=\cL^{(2)}(z_2)R^{+(12)}(q^{2c}z_1/z_2,P+h)\cL^{(1)}(z_1)R^{+(21)}(z_2/z_1,P+h),\nn\\
&&R^{+(12)}(z_1/z_2,P+h)L^{+(1)}(z_1q^c)\cL^{(2)}(z_2)
=\cL^{(2)}(z_2)R^{+(12)}(q^{2c}z_1/z_2,P+h)L^{+(1)}(z_1q^c).\nn
\ena
\end{prop}


\subsection{The trigonometric limit}
Let us consider the trigonometric counterpart of $E_{q,p}(\gh)$ according to an idea  described in \cite{Konno09}.
 Set $L^\pm_{0;ij}(z)=L_{ij}^\pm(z)|_{p=0}$. From \eqref{genL} and \eqref{LmfromLp},
we have  
\be
&&L^\pm_{0;ij}(z)=\sum_{m\in \Z_{\geq0}}L^\pm_{0;ij,\mp m}z^{\pm m}
\qquad (1\leq i,j\leq N)
\en
where
\be
&&L^+_{0;ij,-m}=\bar{L}_{ij,-m}|_{p=0},\qquad 
L^-_{0;ij,m}=q^{2(P+h)_{\bep_i}}\bar{L}_{ij,m}|_{p=0}q^{-2P_{\bep_j}}q^{cm}\qquad (m\in \Z_{\geq 0}).
\en
Let $R^\pm_0(z,s)$ be the trigonometric dynamical $R$ matrix in \eqref{trigDR}. 
From \eqref{RLL}, \eqref{mRLL} and \eqref{RLLpm}, ${L_0^\pm(z)=\sum_{1\leq i,j\leq N}E_{ij}L^\pm_{0;ij}(z)}$ satisfy for $g({P+h}), g(P)\in \FF$
\bea
&&\hspace{-1cm}R_0^{\pm(12)}(z_1/z_2,{P+h})L_0^{\pm(1)}(z_1)L_0^{\pm(2)}(z_2)=L_0^{\pm(2)}(z_2)L_0^{\pm(1)}(z_1)R_0^{\pm(12)}(z_1/z_2,P),\lb{Rlltri}\\
&&\hspace{-1cm}R_0^{\pm(12)}(q^{\pm{c}}z_1/z_2,{P+h})L_0^{\pm(1)}(z_1)L_0^{\mp(2)}(z_2)=L_0^{\mp(2)}(z_2)L_0^{\pm(1)}(z_1)R_0^{\pm(12)}(q^{\mp{c}}z_1/z_2,P).
\lb{Rllpm}\\
&&g({P+h}){L}^\pm_{0;ij}(z)={L}^\pm_{0;ij}(z)g(P+h-<Q_{\bep_i},P+h>),\lb{lgrtri}\\
&&g({P}){L}^\pm_{0;ij}(z)={L}^\pm_{0;ij}(z)g(P-<Q_{\bep_j},P>),\lb{rgrtri}\\
&&[\hd, L_0^\pm(z)]=-z\frac{ \partial}{\partial z}L_0^\pm(z). \lb{dLtri}
\ena
\begin{dfn}
Let $x=(x_1,\cdots,x_N), x_i=q^{2P_{\bep_i}}$ as before. We denote by ${U}^R_{q,x}(\glnh)$  the unital associative algebra over $\FF$ generated by $L^\pm_{0;ij,\mp m}\ (m\in \Z_{\geq 0}), 
\widehat{d}$ and the central element $q^{\pm c/2}$ subject to \eqref{Rlltri}-\eqref{dLtri}. We call ${U}^R_{q,x}(\glnh)$ the dynamical quantum affine algebra in the FRST formulation.  
\end{dfn}
Hence we have 
\begin{prop}
\be
&&E_{q,p}(\glnh)/pE_{q,p}(\glnh)\cong {U}^R_{q,x}(\glnh).
\en
\end{prop}

In order to clarify a relation between the dynamical  ${U}^R_{q,x}(\glnh)$ and the usual quantum affine algebra $U^R_q(\glnh)$ in the FRST 
formulation\cite{RS}, one needs to further remove the $\C[\cR_Q]$ dependence from  ${U}^R_{q,x}(\glnh)$. 
This can be done by considering the algebra generated by the trigonometric limit $L^\pm_0(z,P)$ of  $L^\pm(z,P)$ in \eqref{DL}. Then $L^\pm_0(z,P)$ satisfy the same relations as \eqref{LPcomm}-\eqref{Ldcomm} as well as  the trigonometric limit of the dynamical $RLL$-relations \eqref{DRll}-\eqref{DRllpm}, where   $R^\pm(z,s)$ and $R^{*\pm}(z,s)$ are replaced by $R^\pm_0(z,s)$ and $R_0^{*\pm}(z,s)$, respectively. 
We set 
$L^\pm_{0;ij}(z,P)=\sum_{m\in \Z_{\geq 0}}L^\pm_{0;ij,\mp m}(P)z^{\pm m}$ and denote  by $\widetilde{U}^R_{q,x}(\glnh)$ the unital associative algebra over $\FF$ generated by $L^\pm_{0;ij,m}(P)$.
Then we have 
\be
&&{U}^R_{q,x}(\glnh)\cong \widetilde{U}^R_{q,x}(\glnh)\sharp \C[\cR_Q]. 
\en

Recall that $R^\pm_0(z,P)_{ij}^{kl}$ can be expanded to a formal power series in $x_{i,j}=q^{2P_{i,j}}$ and the 0-th order term gives the trigonometric $R$ matrix in \eqref{trigR}. We assume the same property for  $L_0^\pm(z,P)$. Let 
$L_0^\pm(z,P)=\sum_{|k|=0}^\infty \sum_{k\in \N^N}L^\pm_{0}(z;k)x^k, \ L^\pm_{0}(z;k)=\sum_{1\leq i,j\leq N}E_{ij}L^\pm_{0;ij}(z,k)$, 
where  $k=(k_1,\cdots,k_N), |k|=k_1+\cdots+k_N$ and $x^k=x_1^{k_1}\cdots x_N^{k_N}$. Then $L^\pm_{0}(z;0)=\sum_{1\leq i,j\leq N}E_{ij}L^\pm_{0;ij}(z,0)$,  $L^\pm_{0;ij}(z;0)=\sum_{m\in \Z_{\geq 0}}L^\pm_{0;ij,\mp m}(0)z^{\pm m}$ satisfy the same $RLL$-relations as the quantum affine algebra $U^R_q(\glnh)$ over $\C$ in the FRST formulation\cite{RS}. 
Hence
\be
\widetilde{U}^R_{q,x}(\glnh)/(\sum_i x_i\widetilde{U}^R_{q,x}(\glnh))\cong  U^R_q(\glnh).
\en
　


\section{Hopf Algebroid Structure}
In this section, we introduce an $\h$-Hopf algebroid structure\cite{EV,KR,Konno09} into  
the elliptic algebras $E_{q,p}(\glnh)$ and $U_{q,p}(\glnh)$, and formulate them as elliptic 
quantum groups.

\subsection{$U_{q,p}(\glnh)$ and $E_{q,p}(\glnh)$ as $H$-Algebras}\lb{Halgebra}

Let $\cA$ be an associative algebra, $\cH$ be a 
commutative subalgebra of 
$\cA$, and $\cM_{\cH^*}$ be the 
field of meromorphic functions on $\cH^*$ the dual space of $\cH$. 

\begin{dfn}
An associative algebra $\cA$ with 1 is said to be an $\cH$-algebra, if it is bigraded over 
$\cH^*$, $\ds{\cA=\bigoplus_{\alpha,\beta\in \H^*} \cA_{\al\beta}}$, and equipped with two 
algebra embeddings $\mu_l, \mu_r : \cM_{\cH^*}\to \cA_{00}$ (the left and right moment maps), such that 
\be
\mu_l(\hf)a=a \mu_l(T_\al \hf), \quad \mu_r(\hf)a=a \mu_r(T_\beta \hf), \qquad 
a\in \cA_{\al\beta},\ \hf\in \cM_{\cH^*},
\en
where $T_\al$ denotes the automorphism $(T_\al \hf)(\la)=\hf(\la+\al)$ of $\cM_{\cH^*}$.
\end{dfn}
\begin{dfn}
An $\cH$-algebra homomorphism is an algebra homomorphism $\pi:\cA\to \cB$ between two $\cH$-algebras $\cA$ and $\cB$ preserving the bigrading and the moment maps, i.e. $\pi(\cA_{\al\beta})\subseteq \cB_{\al\beta}$ for all $\al,\beta \in \cH^*$ and $\pi(\mu^\cA_l(\hf))=\mu^\cB_l(\hf), \pi(\mu^\cA_r(\hf))=\mu^\cB_r(\hf)$. 
\end{dfn}

Let $\cA$ and $\cB$ be two $\cH$-algebras. The tensor product $\cA {\widetilde{\otimes}}\cB$ is the $\cH^*$-bigraded vector space with 
\be
 (\cA {\widetilde{\otimes}}\cB)_{\al\beta}=\bigoplus_{\gamma\in\H^*} (\cA_{\al\gamma}\otimes_{\cM_{\cH^*}}\cB_{\gamma\beta}),
\en
where $\otimes_{\cM_{\cH^*}}$ denotes the usual tensor product 
modulo the following 
relation.
\bea
\mu_r^\cA(\hf) a\otimes b=a\otimes\mu_l^\cB(\hf) b, \qquad a\in \cA, 
b\in \cB, \hf\in \cM_{\cH^*}.\lb{AtotB}
\ena
The tensor product $\cA {\widetilde{\otimes}}\cB$ is again an $\cH$-algebra with the multiplication $(a\otimes b)(c\otimes d)=ac\otimes bd$ and the moment maps 
\be
\mu_l^{\cA {\widetilde{\otimes}}\cB} =\mu_l^\cA\otimes 1,\qquad \mu_r^{\cA {\widetilde{\otimes}}\cB} =1\otimes \mu_r^\cB.
\en

Let $\cD$ be the algebra of automorphisms $\cM_{\cH^*}\to \cM_{\cH^*}$ 
\be
\cD&=&\{\ \sum_i \hf_i T_{\beta_i}\ |\ \hf_i\in \cM_{\cH^*},\ \beta_i\in \cH^*\ \}.
\en
Equipped  with the bigrading 
 $\cD_{\al\al}=\{\  \hf T_{-\al}\ |\ \hf\in \cM_{\cH^*},\ \al\in \cH^*\ \}$, 
 $\cD_{\al\beta}=0\ (\al\not=\beta)$ 
 and the moment maps $\mu^{\cD}_l, \mu^{\cD}_r : \cM_{\cH^*}\to \cD_{00}$ 
 defined by 
$\mu^{\cD}_l(\hf)=\mu^{\cD}_r(\hf)=\hf T_0$, $\cD$ is an $\cH$-algebra.
 For any $\cH$-algebra $\cA$, we have the canonical 
isomorphism as an $\cH$-algebra 
\bea
&&\cA\cong \cA\tot \cD\cong  \cD\tot \cA \lb{Diso} 
\ena
by $a\cong a\tot T_{-\beta}\cong T_{-\al}\tot a$ for all $a\in \cA_{\al\beta}$.

Now let $H$ be the same as defined in Sec.\ref{defellRmat} and  take $\cH=H$. 
\begin{prop}\lb{UqpHalg}
The $\cU=U_{q,p}(\glnh)$ is an $H$-algebra by
\bea
&&\cU=\bigoplus_{\al,\beta\in H^*} \cU_{\al,\beta}\lb{bigradingUqp}\\
&&\cU_{\al\beta}=\left\{a\in \cU \left|\ q^{P+h}a q^{-(P+h)}=q^{<\al,P+h>}a,\quad q^{P}a q^{-P}=q^{<\beta,P>}a,\quad \forall P+h, P\in H\right.\right\}\nn
\ena
and $\mu_l, \mu_r : \FF \to \cU_{0,0}$ defined by \cite{Konno09}
\be
&&\mu_l(\hf)=f(P+h,p)\in \FF[[p]],\qquad \mu_r(\hf)=f(P,p^*)\in \FF[[p]].
\en
\end{prop}

\begin{prop}\lb{EqpHalg}
The $\cE=E_{q,p}(\glnh)$ is an $H$-algebra by 
\bea
&&\cE=\bigoplus_{\alpha,\beta\in H^*} \cE_{\alpha,\beta},\lb{bigradingE}\\
&&\cE_{\alpha,\beta}=
\left\{\ a\in \cE\ \left|\ 
q^{P+h}aq^{-(P+h)}=q^{<\al,P+h>}a, \quad 
q^{P}aq^{-P}=q^{<\beta,P>}a, \quad  \forall P+h,  P\in H
\right.\right\}\nn
\ena
and $\mu_l, \mu_r : \FF \to \cE_{0,0}$ defined by the same $\mu_l, \mu_r$ as in $\cU$. 
Note that  $\bar{L}_{ij,n}\in (E_{q,p})_{-Q_{\bep_i}, -Q_{\bep_j}}$.

\end{prop}

We regard $T_\al=e^{-Q_\al}\in \C[\cR_Q]$ as the shift operator $\FF[[p]]\to \FF[[p]]$ 
\be
&&(T_\al \mu_r(\widehat{f}))=e^{-Q_\al}f(P,p^*)e^{Q_\al}={f}(P+<Q_\al,P>,p^*),\\
&&(T_\al \mu_l(\widehat{f}))=e^{-Q_\al}f(P+h,p)e^{Q_\al}={f}(P+h+<Q_\al,P+h>,p).
\en
Then $\cD=\FF\otimes_\C \C[\cR_Q]$ becomes the $H$-algebra  having the property \eqref{Diso} 
for $\cA=\cU, \cE$. 

Hereafter we abbreviate 
$f(P+h,{p})$ and $f(P,{p^*})$ as $f(P+h)$ and
 $f^*(P)$, respectively.

\subsection{$H$-Hopf algebroids $E_{q,p}(\glnh)$ and $U_{q,p}(\glnh)$}

Let us first recall the $\cH$-Hopf algebroid following \cite{EV,KR}. 

\begin{dfn}
An $\cH$-bialgebroid is an $\cH$-algebra $\cA$ equipped with two $\cH$-algebra homomorphisms 
$\Delta:\cA\to \cA{{\tot}}\cA$ (the comultiplication) and $\vep : \cA\to \cD$ (the counit) such that 
\be
&&(\Delta \tot \id)\circ \Delta=(\id \tot \Delta)\circ \Delta,\\
&&(\vep \tot \id)\circ\Delta =\id =(\id \tot \vep)\circ \Delta,
\en
under the identification \eqref{Diso}.
\end{dfn}
 
\begin{dfn}
\lb{defS}
An $\cH$-Hopf algebroid is an $\cH$-bialgebroid $\cA$ equipped with a $\C$-linear map $S : \cA\to \cA$ (the antipode), such that 
\be
&&S(\mu_r(\hf)a)=S(a)\mu_l(\hf),\quad S(a\mu_l(\hf))=\mu_r(\hf)S(a),\quad \forall a\in \cA, \hf\in \cM_{\cH^*},\\
&&m\circ (\id \tot S)\circ\Delta(a)=\mu_l(\vep(a)1),\quad \forall a\in \cA,\\
&&m\circ (S\tot\id  )\circ\Delta(a)=\mu_r(T_{\al}(\vep(a)1)),\quad \forall a\in \cA_{\al\beta},
\en
where $m : \cA{{\tot}} \cA \to \cA$ denotes the multiplication and $\vep(a)1$ is the result of applying the difference operator $\vep(a)$ to the constant function $1\in \cM_{\cH^*}$.
\end{dfn}

\noi
{\it Remark.}\cite{KR} Definition \ref{defS} yields that the antipode of an $\cH$-Hopf algebroid 
 uniquely exists and gives the algebra antihomomorphism.

The $\cH$-algebra $\cD$ is an $\cH$-Hopf algebroid with 
$\Delta_\cD : \cD\to \cD\tot \cD,\ \vep_\cD: \cD \to \cD,\ 
S_\cD : \cD \to \cD$ defined by 
\be
&&\Delta_\cD(\hf T_{-\al})=\hf T_{-\al} \tot T_{-\al},\\
&&\vep_\cD=\id,
\qquad  S_\cD(\hf T_{-\al})=T_{\al}\hf=(T_{\al}\hf)T_{\al}.
\en

Now let us consider the $H$-algebras $\cE$ and $\cU$. 
Let us first consider the $\h$-Hopf algebroid structure on $\cE$. 
We define two $\h$-algebra homomorphisms, the co-unit $\vep : \cE\to \cD$ and the co-multiplication $\Delta : \cE\to \cE \widetilde{\otimes}\cE$ by
\bea
&&\vep(L_{ij,n})=\delta_{i,j}\delta_{n,0}{T}_{{\bep_i} }\quad (n\in \Z),
\qquad \vep(e^Q)=e^Q,\lb{counitUqp}\\
&&\vep(\mu_l({\hf}))= \vep(\mu_r(\hf))=\widehat{f}T_0, \lb{counitf}\\
&&\Delta(L^+_{ij}(z))=\sum_{k}L^+_{ik}(z)\widetilde{\otimes}
L^+_{kj}(z),\lb{coproUqp}\\
&&\Delta(e^{Q})=e^{Q}\tot e^{Q},\qquad \Delta(\hd)=\hd\ \tot1+1\tot \hd,\\
&&\Delta(\mu_l(\hf))=\mu_l(\hf)\widetilde{\otimes} 1,\quad \Delta(\mu_r(\hf))=1\widetilde{\otimes} \mu_r(\hf).\lb{coprof}
\ena
One can check that $\Delta$ preserves the relation in Definition \ref{defEqp}. 

\begin{lem}\lb{counitcopro}
The maps $\vep$ and $\Delta$ satisfy
\bea
&&(\Delta\tot \id)\circ \Delta=(\id \tot \Delta)\circ \Delta,\lb{coaso}\\
&&(\vep \tot \id)\circ\Delta =\id =(\id \tot \vep)\circ \Delta.\lb{vepDelta}
\ena
\end{lem}
\noi
{\it Proof.} Straightforward. 
\qed

We also have the following formulae. 
\begin{prop}\lb{Coprofsoverf}
\bea
&&\Delta\left(\frac{f(P,p^*)}{f(P+h,p)}\right)=\frac{f(P,p^*)}{f(P+h,p)}\tot
\frac{f(P,p^*)}{f(P+h,p)}.\lb{fsf}
\ena
\end{prop}
Hence $(\cE, \Delta, \cM_{H^*}, \mu_l, \mu_r, \vep)$ is a $H$-bialgebroid. 

We define  an algebra antihomomorphism (the antipode) $S : \cE\to \cE$ by
\bea
&&S(L^+_{ij}(z))=(L^+(z)^{-1})_{ij},\lb{SonL}\\
&&S(e^Q)=e^{-Q},\quad S(\mu_r(\hf))=\mu_l(\hf),\quad S(\mu_l(\hf))=\mu_r(\hf).
\ena
The explicit formula for \eqref{SonL}  in terms of the components of the $L$-operator is given in Appendix \ref{ellqdet}.  
Then $S$ preserves the $RLL$ relation \eqref{RLL} and satisfies the antipode axioms.
We hence obtain 
\begin{thm}
The $\h$-algebra $\cE$ equipped with $(\Delta,\vep,S)$ is an $\h$-Hopf algebroid. 
\end{thm}

\begin{dfn}
We call the $\h$-Hopf algebroid $(\cE,\h,{\cM}_{\h^*},\mu_l,\mu_r,\Delta,\vep,S)$ the  elliptic quantum group $E_{q,p}(\glnh)$. 
\end{dfn}

\noindent
{\it Remark.}\  The coproduct for $L^+(z,P)$ used in \cite{Felder,FV} is essentially obtained from \eqref{coproUqp} via \eqref{DL}:
\be
&&\Delta(L^+(z,P))=L^+(z,P)\otimes   L^+(z,P+h^{(1)}).
\en

By making use of the isomorphism between $\cU$ and $\cE$ given in Sec.\ref{Iso}, we can define the $L$-operators of $\cU$ by identifying them with  those of $\cE$ in \eqref{def:lhat}. Then the $H$-Hopf algebroid structure of $\cU$ is  defined by using  the same 
$\Delta, \vep, S$ as $\cE$.  See \cite{Konno09} for the $\slth$ case.

\begin{dfn}
We call the $\h$-Hopf algebroid $(\cU,\h,{\cM}_{\h^*},\mu_l,\mu_r,\Delta,\vep,S)$ the  elliptic quantum group $U_{q,p}(\glnh)$. 
\end{dfn}
Hence the isomorphism obtained in Sec.\ref{Iso} can be  extended to as an  $H$-Hopf algebroid. 

\noindent
{\it Remark.}\  $U_{q,p}(\gh)$ admits another co-algebra structure through 
another coproduct called the Drinfeld coproduct\cite{JKOS,Konno14}. 

\section{Dynamical Representations}\lb{DRep}

\subsection{Definition}
We summarize some basic facts on the dynamical representation of $U_{q,p}(\glnh)$. 
Most of them can be extended to the arbitrary untwisted affine Lie algebra $\gh$ case\cite{FKO}. 

Let us consider a vector space $\hV$ over $\FF[[p]]$, which is  
${H}$-diagonalizable, i.e.  
\be
&&\hV=\bigoplus_{\la,\mu\in {H}^*}\hV_{\la,\mu},\ \hV_{\la,\mu}=\{ v\in \hV\ |\ q^{P+h}\cdot v=q^{<\la,P+h>} v,\ q^{P}\cdot v=q^{<\mu,P>} v\ \forall 
P+h, P\in 
{H}\}.
\en
Let us define the $H$-algebra $\cD_{H,\hV}$ of the $\C$-linear operators on $\hV$ by
\be
&&\cD_{H,\hV}=\bigoplus_{\al,\beta\in {H}^*}(\cD_{H,\hV})_{\al\beta},\\
&&\hspace*{-10mm}(\cD_{H,\hV})_{\al\beta}=
\left\{\ X\in \End_{\C}\hV\ \left|\ 
\mmatrix{ f(P+h)X=X f(P+h+<\alpha,P+h>),\qquad\qquad\cr 
f(P)X=X f(P+<\beta,P>),\ \forall f(P), f(P+h)\in \FF[[p]],\cr
 X\cdot\hV_{\la,\mu}\subseteq 
 \hV_{\la+\al,\mu+\beta} \qquad\qquad\qquad\qquad\qquad\qquad \cr}  
 \right.\right\},\\
&&\mu_l^{\cD_{H,\hV}}(\widehat{f})v=f(<\la,P+h>,p)v,\quad 
\mu_r^{\cD_{H,\hV}}(\widehat{f})v=f(<\mu,P>,p^*)v,\quad \widehat{f}\in \FF[[p]],
\ v\in \hV_{\la,\mu}.
\en
\begin{dfn}
We define a dynamical representation of $U_{q,p}(\glnh)$ on $\hV$ to be  
 an $H$-algebra homomorphism ${\pi}: U_{q,p}(\glnh) 
 \to \cD_{H,\hV}$. By the action $\pi$ of $U_{q,p}(\glnh)$ we regard $\hV$ as a 
$U_{q,p}(\glnh)$-module. 
\end{dfn}

\begin{dfn}
For $k\in \C$, we say that a $U_{q,p}(\glnh)$-module has  level $k$ if $q^{\pm c/2}$ acts  
as the scalar $q^{\pm k/2}$ on it.  
\end{dfn}

\begin{dfn}
Let ${\goth H}$, ${\goth N}_+, {\goth N}_-$ be the subalgebras of 
$U_{q,p}(\glnh)$ 
generated by 
$q^{\pm c/2}, d, 
k_{i,0}\ (1\leq i\leq N)$, by $k_{i,n}\ ( 1\leq i\leq N, n\in \Z_{>0})$,  $e_{j, n}\ (j\in I, n\in \Z_{\geq 0})$  
$f_{j, n}\ (j\in I, n\in \Z_{>0})$ and by $k_{i,-n}\ (1\leq i\leq N, n\in \Z_{>0}),\ e_{j, -n}\ (j\in I, n\in \Z_{> 0}),\ f_{j, -n}\ (i\in I, n\in \Z_{\geq 0})$, respectively.   
\end{dfn}

\begin{dfn}
For $k\in\C$, $\la, \mu\in H^*$, 
a dynamical $U_{q,p}(\gh)$-module $\hV(\la,\mu)$ is called the 
level-$k$ highest weight module with the highest weight $(\la,\mu)$, if there exists a vector 
$v\in \hV(\la,\mu)$ such that
\be
&&\hV(\la,\mu)=U_{q,p}(\gh)\cdot v,\qquad {\goth N}_+\cdot v=0,\qquad q^{\pm c/2}\cdot v=q^{\pm k/2}v, \\
&&k_{i,0}\cdot v=q^{-<\la-\mu,h_{\bep_i}>}v, 
\quad  f({P})\cdot v =f({<\mu,P>})v,\quad f({P+h})\cdot v =f({<\la,P+h>})v.
\en
\end{dfn}

\subsection{The evaluation $H$-algebra homomorphism}\lb{IndDRep}

Let $k^\pm_{0,i}(z), x^\pm_j(z)$ be the Drinfeld currents of the quantum affine algebra $U_q(\glnh)$. See Appendix \ref{Uqgh}. 
Let us introduce the currents $u^+_{\vep_i}(z,p)\in (U_q(\glnh)[[p]])[[z]], u^-_{\vep_i}(z,p)\in (U_q(\glnh)[[p]])[[z^{-1}]]\ (1\leq i\leq N)$ by 
\bea
&&u^+_{\vep_i}(z,p)=\prod_{n=1}^\infty \left(k^-_{i,0} \cdot k^+_{0,i}(p^{*n}q^{c-i}z)\right),\\
&&u^-_{\vep_i}(z,p)=\prod_{n=1}^\infty \left(k^+_{i,0} \cdot k^-_{0,i}(p^{-n}q^{c-i}z)\right).
\ena
We also set 
\bea
&&u_j^\pm(z,p)=u^\pm_{\vep_j}(z,p)u^\pm_{\vep_{j+1}}(qz,p)^{-1}\qquad (1\leq j\leq N-1).
\ena
These are well defined elements in $(U_{q}(\glnh)[[p]])[[z,z^{-1}]]$ in the $p$-adic 
topology. 

Now let us define the `dressed' currents $x_j^\pm(z,p)\ (1\leq j\leq N-1)$, $k^{\pm}_{i}(z,p)\ (1\leq i\leq N)$ by 
\bea
&& x_j^+(z,p)=u_j^+(z,p)x_j^+(z) e^{-Q_{\al_j}}, 
\lb{gdress1}
\\
&&x_j^-(z,p)=x_j^-(z)u_j^-(z,p),
\lb{gdress2}
\\
&&k^{+}_{i}(z,p)
=u_{\vep_i}^+(q^{-c+j}z,p)k^+_{0,i}(z)u_{\vep_i}^-(q^{j}z,p)e^{-Q_{\bep_i}},\lb{gdress3} 
\\
&&k^{-}_{i}(z,p)
=u_{\vep_i}^+(q^{j} z,p)k^-_{0,i}(z)u_{\vep_i}^-(q^{-c+j}z,p)e^{-Q_{\bep_i}}.
\lb{gdress4}
\ena

\begin{thm}\lb{evHhom}
The map  $\phi_p: U_{q,p}(\glnh)[[z,z^{-1}]] \to (\FF[[p]]\otimes_\C U_{q}(\glnh))[[z,z^{-1}]]\sharp\C[\cR_Q]$ defined by 
\be
&& e_i(z) \mapsto x^+_i(z,p) , \qquad f_i(z) \mapsto  x^-_i(z,p), \qquad k^\pm_i(z) \mapsto k^\pm_i(z,p)
\en
is an $H$-algebra homomorphism. We call $\phi_p$ the evaluation $H$-algebra homomorphism. 
\end{thm}
\noindent
{\it Proof.}\ Direct calculations using Lemma \ref{commuu}. \qed

Let $(\varphi_V, V)$ be a representation of $U_q(\glnh)$. We assume $V$ is an ${\hh}$-diagonalizable vector space over 
$\C$.  We set $V_{\FF[[p]]}=\FF[[p]]\otimes_\C V$. 
Let $V_Q$ be a vector space over $\C$, 
on which an action of $e^{ Q}$ is defined appropriately.  
Two important examples of $V_Q$ are  
$ V_Q=\C 1$ and $V_Q=\C[\cR_Q]$, where 
$1$ denotes the vacuum state satisfying $e^{Q}.1=1$. 
Let us consider the vector space $\hV_{\FF[[p]]}= V_{\FF[[p]]}\otimes_{\C} V_Q$, on which 
the actions of $f(P,h,p)\in \FF[[p]]$ and $e^{ Q}$ are defined as follows. 
For $v\otimes \xi\in V\otimes V_Q$,
\be
&&f(P,h,p).(v\otimes \xi)=f(P,{\rm wt}(v),p)v\otimes \xi,\\
&&e^{Q}.(f(P,h,p)v\otimes \xi)=f(P-<Q,P>,h,p)v\otimes e^{Q}\xi,
\en
where $h.v={\rm wt}(v)v$. 
We extend $\varphi_V:U_q(\glnh)\to \End_{\C} V$ 
to a dynamical representation $\varphi_V:(\FF[[p]]\otimes_C U_q(\glnh))\sharp \C[\cR_Q]\to \cD_{H,\hV_{\FF[[p]]}}$ by
\be
&&\varphi_V(f(P))=f(P),\qquad \varphi_V(e^{Q_\al})=e^{Q_\al}\qquad \forall e^{Q_\alpha}\in \C[\cR_Q]. 
\en
 Note that  if we specialize $p$ as 0, $\widehat{V}_{\FF[[p]]} $ becomes 
$V_{\FF}\otimes V_Q$, where $V_\FF=\FF\otimes_\C V$. 

Then from Theorem \ref{evHhom} we obtain the following. 
\begin{cor}\lb{dynamicalrep}
A map ${\varphi}^p_V=\varphi_V\circ \phi_p
: U_{q,p}(\glnh)\to \cD_{H,\hV_{\FF[[p]]}}$  
 gives a dynamical representation of 
$U_{q,p}(\glnh)$ on $\hV_{\FF[[p]]}$. We call $(\varphi^p_V,\widehat{V}_{\FF[[p]]})$ the evaluation dynamical representation. 

\end{cor}
Due to this corollary, any representation of $U_q(\glnh)$ admits an `elliptic
 and dynamical deformation' for generic $p$.  One can easily extend this to any untwisted affine Lie algebra case by 
 using the evaluation homomorphism given in Appendix A of \cite{JKOS} . See \cite{Konno09} for $\slth$ case. 



\section{Isomorphism Between $U_{q,p}(\glnhbig)$ and $E_{q,p}(\glnhbig)$}\lb{Iso}

We introduce the Gauss components of the $L$ operator of $\cE=E_{q,p}(\glnh)$ and 
the half currents of $\cU=U_{q,p}(\glnh)$. Then we show the isomorphism between $\cU$ and $\cE$. 
 
\subsection{The $L$-operators of $\cE$}

Let us set
\be
&&\cE^\pm=\left\{A(z)\in \cE[[p]][[z,z^{-1}]]\ \left|\ A(z)\in \cE[[z^{\pm 1}]] \mod p\cE[[p]][[z,z^{-1}]]
\ \right\}\right. .
\en
Then it is easy to show 
\begin{lem}\lb{AB}
For $A(z), B(z)\in \cE^\pm$, the product $A(z)B(z)$ is a well-defined element in  $\cE^\pm$ in the $p$-adic topology, respectively. Conversely, if $A(z), B(z)\in \cE[[p]][[z,z^{-1}]]$ satisfy $A(z)B(z)\in \cE^\pm$, then   
$A(z), B(z)\in \cE^\pm$, respectively.
\end{lem}

\begin{df}~~We define the Gauss components $E^\pm_{l,j}(z), F^\pm_{j,l}(z), K^\pm_m(z) \ (1\leq j<l\leq N, 1\leq m\leq N)$ of the $L$-operator
${L}^\pm(z)$ of  $\cE$ as follows.
\begin{eqnarray}
&&{L}^\pm(z)=
\left(\begin{array}{ccccc}
1&F_{1,2}^\pm(z)&F_{1,3}^\pm(z)&\cdots&F_{1,N}^\pm(z)\\
0&1&F_{2,3}^\pm(z)&\cdots&F_{2,N}^\pm(z)\\
\vdots&\ddots&\ddots&\ddots&\vdots\\
\vdots&&\ddots&1&F_{N-1,N}^\pm(z)\\
0&\cdots&\cdots&0&1
\end{array}\right)\left(
\begin{array}{cccc}
K^\pm_1(z)&0&\cdots&0\\
0&K^\pm_2(z)&&\vdots\\
\vdots&&\ddots&0\\
0&\cdots&0&K^\pm_{N}(z)
\end{array}
\right)\nn\\
&&\qquad\qquad\qquad\qquad\qquad\qquad\qquad\times
\left(
\begin{array}{ccccc}
1&0&\cdots&\cdots&0\\
E^\pm_{2,1}(z)&1&\ddots&&\vdots\\
E^\pm_{3,1}(z)&
E^\pm_{3,2}(z)&\ddots&\ddots&\vdots\\
\vdots&\vdots&\ddots&1&0\\
E^\pm_{N,1}(z)&E^\pm_{N,2}(z)
&\cdots&E^\pm_{N,N-1}(z)&1
\end{array}
\right).\lb{def:lhat}
\ena
In particular we call $E^\pm_{j+1,j}(z), F^\pm_{j,j+1}(z), K^\pm_m(z)$ the basic Gauss components.
\end{df}
\noindent
{\it Remark.}\lb{Gausscom} By definition the matrix elements $L^\pm_{i,j}(z)$ are the elements 
in $\cE^\pm$, respectively. Then from Lemma \ref{AB}, the matrix elements $E_{l,j}^\pm(z), F_{j,l}^\pm(z), K_m^\pm(z)$ of the right hand side of \eqref{def:lhat} are elements in $\cE^\pm$, respectively and their products are well defined formal Laurent series in $z$ in the $p$-adic topology.  
In addition, since  $L^\pm(z)$ are invertible, $K_m^\pm(z)\ (1\leq m\leq N)$ are invertible.  Therefore all the components $E_{l,j}^\pm(z), F_{j,l}^\pm(z)$ and
$K_m^\pm(z)$ are determined uniquely by $L^\pm_{ij}(z)$, respectively.

Hence we define the coefficients of the Gauss components $E_{l,j}^+(z), F_{j,l}^+(z), K_m^+(z)$ as follows. 
\begin{dfn}\lb{modeGC}
\bea
&&E^+_{l,j}(z)=\sum_{n\in \Z_{\geq 0}}E^+_{l,j,-n}z^n+\sum_{n\in \Z_{> 0}}E^+_{l,j,n}p^nz^{-n},\lb{expEp}\\
&&F^+_{l,j}(z)=\sum_{n\in \Z_{\geq 0}}F^+_{l,j,-n}z^n+\sum_{n\in \Z_{> 0}}F^+_{l,j,n}p^nz^{-n},\lb{expFp}\\
&&K^+_{j}(z)=\sum_{n\in \Z_{\geq 0}}K^+_{j,-n}z^n+\sum_{n\in \Z_{> 0}}K^+_{j,n}p^nz^{-n}.\lb{expKp}
\end{eqnarray}
\end{dfn}
In addition, from the definition of $L^-(z)$ \eqref{LmfromLp}, we have 
\bea
&&E^-_{j,i}(z)=q^{2P_{\bep_j}}E^+_{j,i}(zpq^{-c})q^{-2P_{\bep_i}},\lb{defEm}\\
&&F^-_{i,j}(z)=q^{2(P+h)_{\bep_i}}F^+_{i,j}(zpq^{-c})q^{-2(P+h)_{\bep_j}},\lb{defFm}\\
&&K^-_{i}(z)=q^{\frac{2(N-1)}{N}}q^{2(P+h)_{\bep_i}}K^+_{i}(zpq^{-c})q^{-2P_{\bep_i}}. \lb{defKm}
\ena
Hence we define  
\begin{dfn}
\bea
&&E^-_{j,i,n}=q^{2P_{\bep_j}}E^+_{j,i,n}q^{-2P_{\bep_i}},\quad 
F^-_{i,j,n}=q^{2(P+h)_{\bep_i}}F^+_{i,j,n}q^{-2(P+h)_{\bep_j}},\quad\nn\\ 
&&K^-_{i,n}=q^{\frac{2(N-1)}{N}}q^{2(P+h)_{\bep_i}}K^+_{i,n}q^{-2P_{\bep_i}}, 
\lb{mhccom}
\ena
for $n\in \Z$. 
\end{dfn}
Then  we have 
\bea
&&E^-_{j,i}(z)=\sum_{n\in \Z_{>0}}E^-_{l,j,-n}p^n(q^{-c}z)^n+\sum_{n\in \Z_{\geq 0}}E^-_{l,j,n}(q^{-c}z)^{-n},\lb{expEm}\\
&&F^-_{i,j}(z)=\sum_{n\in \Z_{>0}}F^-_{l,j,-n}p^n(q^{-c}z)^n+\sum_{n\in \Z_{\geq 0}}F^-_{l,j,n}(q^{-c}z)^{-n},\lb{expFm}\\
&&K^-_{i}(z)=\sum_{n\in \Z_{>0}}K^-_{i,-n}p^n(q^{-c}z)^n+\sum_{n\in \Z_{\geq 0}}K^-_{i,n}(q^{-c}z)^{-n}. \lb{expKm}
\ena

\subsection{Subalgebras}
For  $1\leq l<N$, let us define the reduced $R$-matrix and $L$-operators by
\bea
&&R^\pm_l(z,s)=(R^\pm(z,s)_{ij}^{i'j'})_{l\leq i,j, i',j'\leq N},\lb{Rl}\\
&&L^\pm_l(z)=(L^\pm_{ij}(z))_{l\leq i,j\leq N}.\lb{Ll}
\ena
Note that up to overall factors $R^\pm_l(z,s)$ are the elliptic dynamical $R$ matrix of type $A_{N-l}^{(1)}$. 
Note also that   
if  $R^\pm(z,s)_{ij}^{i'j'}\not=0$ for $1\leq i,j\leq l$ (resp. $l\leq i,j\leq N$), then $1\leq i',j'\leq l$ (resp. $l\leq i',j'\leq N$). 
Hence we obtain
\begin{prop}\lb{L_l}
The reduced $L$-operators $L_l^\pm(z)$ satisfy 
\bea
&&\hspace{-1.5cm}R^{\pm(12)}_l(z_1/z_2,P+h)L^{\pm(1)}_l(z_1)L^{\pm(2)}_l(z_2)=L^{\pm(2)}_l(z_2)L^{\pm(1)}_l(z_1)R^{-*(12)}_l(z_1/z_2,P),\lb{mRLLl}\\
&&\hspace{-1.5cm}R^{\pm(12)}_l(q^{\pm{c}}z_1/z_2,P+h)L^{\pm(1)}_l(z_1)L^{\mp(2)}_l(z_2)=L^{\mp(2)}_l(z_2)L^{\pm(1)}_l(z_1)R^{\pm*(12)}_l(q^{\mp{c}}z_1/z_2,P).\lb{RLLpml}
\ena
Therefore $L^+_l(z), \hd, q^{\pm c/2}$ generate a subalgebra of $E_{q,p}(\glnh)$, which is isomorphic to $E_{q,p}(\widehat{\gl}_{N-l+1})$.
\end{prop}
Note that from \eqref{def:lhat} we have 
\bea
&&\hspace{-1cm}L^\pm_l(z)\nonumber\\
&&\hspace{-1cm}=
\left(\begin{array}{ccccc}
1&F_{l,l+1}^\pm(z)&F_{l,l+2}^\pm(z)&\cdots&F_{l,N}^\pm(z)\\
0&1&F_{l+1,l+2}^\pm(z)&\cdots&F_{l+1,N}^\pm(z)\\
\vdots&\ddots&\ddots&\ddots&\vdots\\
\vdots&&\ddots&1&F_{N-1,N}^\pm(z)\\
0&\cdots&\cdots&0&1
\end{array}\right)\left(
\begin{array}{cccc}
K^\pm_l(z)&0&\cdots&0\\
0&K^\pm_{l+1}(z)&&\vdots\\
\vdots&&\ddots&0\\
0&\cdots&0&K^\pm_{N}(z)
\end{array}
\right)\nn\\
&&\qquad\qquad\qquad\qquad\qquad\times
\left(
\begin{array}{ccccc}
1&0&\cdots&\cdots&0\\
E^\pm_{l+1,l}(z)&1&\ddots&&\vdots\\
E^\pm_{l+2,l}(z)&
E^\pm_{l+2,l+1}(z)&\ddots&\ddots&\vdots\\
\vdots&\vdots&\ddots&1&0\\
E^\pm_{N,l}(z)&E^\pm_{N,l+1}(z)
&\cdots&E^\pm_{N,N-1}(z)&1
\end{array}
\right).\lb{Laa}
\ena
Hence
\be
&&\hspace{-1cm}L^\pm_l(z)^{-1}\nn\\
&&\hspace{-1cm}=\left(
\begin{array}{cccc}
K_l^{\pm-1}&-K_l^{\pm-1}F_{l,l+1}^\pm&K_l^{\pm -1}x_l&\cdots\\
-E_{l+1,l}^\pm K_l^{\pm-1}&
E_{l+1,l}^\pm K_l^{\pm -1}F_{l,l+1}^\pm+K_{l+1}^{\pm -1}
&
-E_{l+1,l}^\pm K_l^{\pm -1}x_l-K_{l+1}^{\pm -1}F_{l+1,l+2}^\pm
&\cdots\\
y_lK_l^{\pm -1}&-y_lK_l^{\pm -1}F_{l,l+1}^\pm-E_{l+2,l+1}^\pm K_{l+1}^{\pm -1}&\ddots&\vdots \\
\vdots &\vdots &\cdots &\ddots \\
\end{array}
\right),
\en
where we omitted the argument $z$ and set
\begin{eqnarray}
x_l(z)=F_{l,l+1}^\pm(z)F_{l+1,l+2}^\pm(z)-F_{l,l+2}^\pm(z),\\
y_l(z)=E_{l+2,l+1}^\pm(z)E_{l+1,l}^\pm(z)-E_{l+2,l}^\pm(z).
\end{eqnarray}

Furthermore, for $l<m\leq N$ let us define 
\bea
&&R^{m\pm}_l(z,s):=(R^\pm(z,s)_{ij}^{i'j'})_{l\leq i,j, i',j'\leq m},\lb{Rlm}\\
&&(L^\pm_l(z)^{-1})^m:=((L^\pm_l(z)^{-1})_{ij})_{l\leq i,j\leq m}.\lb{Linvlm}
\ena
\begin{prop}
If we set
\bea
&&\hspace{-1cm}L^{m\pm}_l(z)\nonumber\\
&=&
\left(\begin{array}{ccccc}
1&F_{l,l+1}^\pm(z)&F_{l,l+2}^\pm(z)&\cdots&F_{l,m}^\pm(z)\\
0&1&F_{l+1,l+2}^\pm(z)&\cdots&F_{l+1,m}^\pm(z)\\
\vdots&\ddots&\ddots&\ddots&\vdots\\
\vdots&&\ddots&1&F_{m-1,m}^\pm(z)\\
0&\cdots&\cdots&0&1
\end{array}\right)\left(
\begin{array}{cccc}
K^\pm_l(z)&0&\cdots&0\\
0&K^\pm_{l+1}(z)&&\vdots\\
\vdots&&\ddots&0\\
0&\cdots&0&K^\pm_{m}(z)
\end{array}
\right)\nn\\
&&\qquad\qquad\qquad\qquad\qquad\times
\left(
\begin{array}{ccccc}
1&0&\cdots&\cdots&0\\
E^\pm_{l+1,l}(z)&1&\ddots&&\vdots\\
E^\pm_{l+2,l}(z)&
E^\pm_{l+2,l+1}(z)&\ddots&\ddots&\vdots\\
\vdots&\vdots&\ddots&1&0\\
E^\pm_{N,l}(z)&E^\pm_{m,l+1}(z)
&\cdots&E^\pm_{m,m-1}(z)&1
\end{array}
\right).\lb{Llm}
\ena
Then  we have
\bea
&&(L^+_l(z)^{-1})^m=L^{m\pm}_l(z)^{-1}.\lb{Llmeq}
\ena
\end{prop}
Note that 
\be
&&L^{m\pm}_l(z)\not=(L^\pm_{ij}(z))_{l\leq i,j\leq m}.
\en
Hence we have
\begin{lem}
The restriction of the relations 
\bea
&&\hspace{-1.5cm}L^{\pm(1)}_l(z_1)^{-1}L^{\pm(2)}_l(z_2)^{-1}R^{\pm(12)}_l(z,P+h)=R^{-*(12)}_l(z,P)
L^{\pm(2)}_l(z_2)^{-1}L^{\pm(1)}_l(z_1)^{-1},\lb{RLLlinv}\\
&&\hspace{-1.5cm}L^{\pm(1)}_l(z_1)^{-1}L^{\mp(2)}_l(z_2)^{-1}R^{\pm(12)}_l(zq^{\pm{c}},P+h)=
R^{\pm*(12)}_l(zq^{\mp{c}},P)L^{\mp(2)}_l(z_2)^{-1}L^{\pm(1)}_l(z_1)^{-1}
\lb{RLLpmlinv}
\ena
to the $(i,j), (i'',j'')$ components  with $l\leq i, j, i'', j''\leq m$ are 
equivalent to
\bea
&&\hspace{-1.5cm}R^{m\pm(12)}_l(z,P+h)L^{m\pm(1)}_l(z_1)L^{m\pm(2)}_l(z_2)
=L^{m\pm(2)}_l(z_2)L^{m\pm(1)}_l(z_1)R^{m\pm*(12)}_l(z,P+h),\lb{mRLLlm}\\
&&\hspace{-1.5cm}R^{m\pm(12)}_l(zq^{\pm{c}},P+h)L^{m\pm(1)}_l(z_1)L^{m\mp(2)}_l(z_2)
=L^{m\mp(2)}_l(z_2)L^{m\pm(1)}_l(z_1)R^{m\pm*(12)}_l(zq^{\mp{c}},P),
\lb{RLLpmlm}
\ena
where $z=z_1/z_2$. 
\end{lem}
Therefore we obtain the following statement.
\begin{thm}\lb{L_l^m}
$L^{m+}_l(z),  \hd, q^{\pm c/2}$ generate the subalgebra of $E_{q,p}(\glnh)$, which is isomorphic to $E_{q,p}(\widehat{\gl}_{m-l+1})$. 
\end{thm}
From \eqref{mRLLlm} and \eqref{RLLpmlm}, we have
\begin{prop}
\bea
&&\hspace{-1cm}L^{m\pm(2)}_l(z_2)^{-1}R^{m\pm(12)}_l(z,P+h)L^{m\pm(1)}_l(z_1)
=L^{m\pm(1)}_l(z_1)R^{m\pm*(12)}_l(z,P)L^{m\pm(2)}_l(z_2)^{-1},
\lb{mRLLlminv}\\
&&\hspace{-1cm}L^{m\pm(1)}_l(z_1)^{-1}L^{m\pm(2)}_l(z_2)^{-1}R^{m\pm(12)}_l(z,P+h)
=R^{m\pm*(12)}_l(z,P)L^{m\pm(2)}_l(z_2)^{-1}L^{m\pm(1)}_l(z_1)^{-1},
\lb{mRLLlminvinv}\\
&&\hspace{-1cm}L^{m\mp(2)}_l(z_2)^{-1}R^{m\pm(12)}_l(zq^{\pm{c}},P+h)L^{m\pm(1)}_l(z_1)
=L^{m\pm(1)}_l(z_1)R^{m\pm*(12)}_l(zq^{\mp{c}},P)L^{m\mp(2)}_l(z_2)^{-1},
\lb{RLLpmlminv}\\
&&\hspace{-1cm}L^{m\pm(1)}_l(z_1)^{-1}L^{m\mp(2)}_l(z_2)^{-1}R^{m\pm(12)}_l(zq^{\pm{c}},P+h)=R^{m\pm*(12)}_l(zq^{\mp{c}},P)L^{m\mp(2)}_l(z_2)^{-1}L^{m\pm(1)}_l(z_1)^{-1}.\nn\\
&&\lb{RLLpmlminvinv}
\ena

\end{prop}
\begin{lem}\lb{KKKE}
For $2\leq l+1< m \leq N$, we have
\bea
&&\rho^\pm(z)\bb(z)K^\pm_l(z_2)^{-1}K^\pm_m(z_1)
=\rho^{*\pm}(z)\bb(z)^*K^\pm_m(z_1)K^\pm_l(z_2)^{-1},\lb{KlKm}\\
&&\rho^\pm(q^{\pm c}z)\bb(q^{\pm c}z)K^\mp_l(z_2)^{-1}K^\pm_m(z_1)
=\rho^{*\pm}(q^{\mp c}z)\bb(q^{\mp c}z)^*K^\pm_m(z_1)K^\mp_l(z_2)^{-1},\lb{KlKmpm}\\
&&[E^\pm_{l+1,l}(z_1), K^\pm_m(z_2)]=0=[E^\pm_{m,l+1}(z_1), K^\pm_l(z_2)],
\lb{EKcomm} \\
&&[E^\mp_{l+1,l}(z_2), K^\pm_m(z_1)]=0=[E^\pm_{m,l+1}(z_1), K^\mp_l(z_2)]
\lb{EKpmcomm} 
\ena
where $z=z_1/z_2$. 
\end{lem}
\noindent
{\it Proof.}  \eqref{KlKm} and \eqref{KlKmpm} follows from the $(m,l), (m,l)$ component of \eqref{mRLLlminv} and \eqref{RLLpmlminv}, respectvely. Similarly, \eqref{EKcomm} ( resp. \eqref{EKpmcomm} ) follows from \eqref{KlKm}  and the $(m,l+1), (m,l)$ component of \eqref{mRLLlminv}(resp. \eqref{KlKmpm} and the same component of \eqref{RLLpmlminv}  ).
\qed

Other relations among the basic Gauss components are given in Appendix \ref{hc2ellc}. 

The following lemma indicates that the whole Gauss components of $L^\pm(z)$ can be determined recursively by the basic ones. 
\begin{lem}\lb{recursion}
Let $I_{a,b}=\{\ (j,k)\ |\  a\leq j\leq b-1,\  j+1\leq k\leq b \}\setminus \{(a,b)\} $. 
For $2\leq  l+1<m\leq N$,  $E^+_{m,l}(z)$ (resp. $F^+_{l,m}(z)$) is determined by $\{E^+_{m,l}(z) \ (l,m)\in I_{l,m}, 
K^+_j(z)\ l\leq j\leq m \}$ (resp. $\{F^+_{l,m}(z) \ (l,m)\in I_{l,m}, 
K^+_j(z)\ l\leq j\leq m \}$). 
\end{lem}
\noindent
{\it Proof.} 
Let us consider $E^+_{m,l}(z)$. The $F^+_{l,m}(z)$ case is similar.  
From \eqref{KlKm} and the $(m, l+1), (l+1,l)$ component of \eqref{mRLLlminv} 
we have 
\bea
&&-\bar{b}^*(z)E_{l+1,l}^\pm(z_2^\pm)E_{m,l+1}^\pm(z_1^\pm)\nn\\
&&=E_{m,l}^\pm(z_1^\pm)c^*(z,P_{l,l+1})
-E_{m,l+1}^\pm(z_1^\pm)E_{l+1,l}^\pm(z_2^\pm)
+\bar{c}^*(z,P_{l+1,m})(L_l^{m\pm}(z_2^\pm)^{-1})_{ml}K^\pm_l(z_2^\pm)\nn\\
&&+\sum_{l+2\leq k\leq m-1}E_{m,k}^\pm(z_1^\pm)\bar{c}^*(z,P_{l+1,k})
(L_l^{m\pm}(z_2^\pm)^{-1})_{kl}K^\pm_l(z_2^\pm)
.\lb{E21EN2}
\ena
Similarly, from \eqref{KlKmpm} and the $(m,l+1), (l+1,l)$ component of \eqref{RLLpmlminv}, we have 
\bea
&&-\bar{b}^*(z)E_{l+1,l}^\mp(z_2^\mp)E_{m,l+1}^\pm(z_1^\pm)\nn\\
&&=E_{m,l}^\pm(z_1^\pm)c^*(z,P_{l,l+1})
-E_{m,l+1}^\pm(z_1^\pm)E_{l+1,l}^\mp(z_2^\mp)
+\bar{c}^*(z,P_{l+1,m})(L_l^{m\mp}(z_2^\mp)^{-1})_{ml}K^\mp_l(z_2^\mp)\nn\\
&&+\sum_{l+2\leq k\leq m-1}E_{m,k}^\pm(z_1^\pm)\bar{c}^*(z,P_{l+1,k})
(L_l^{m\mp}(z_2^\mp)^{-1})_{kl}K^\mp_l(z_2^\pm)
.\lb{E21mEN2p}
\ena
Subtracting \eqref{E21mEN2p}  with the upper and lower signs reversed from \eqref{E21EN2}, 
we have 
\bea
&&E_{m,l}^\pm(z_1^\pm)-E_{m,l}^\mp(z_1^\mp)\nn\\
&&=\frac{1}{{c}^*(z,P_{l,l+1})}\Biggl(\left(E_{m,l+1}^\pm(z_1^\pm)-E_{m,l+1}^\mp(z_1^\mp)\right)E_{l+1,l}^\pm(z_2^\pm)\nn\\
&&\qquad\qquad\qquad\qquad -\bar{b}^*(z)E_{l+1,l}^\pm(z_2^\pm)\left(
E_{m,l+1}^\pm(z_1^\pm)-E_{m,l+1}^\mp(z_1^\mp)\right)\nn\\
&&\left.\quad-\sum_{l+2\leq k\leq m-1}(E^\pm_{m,k}(z_1^\pm)-E^\mp_{m,k}(z_1^\mp))
\bar{c}^*(z,P_{l+1,k})(L_l^{m\pm}(z_2^\pm)^{-1})_{kl}K^\pm_l(z_2^\pm)
\right).
\lb{EN1byE21EN2p}
\ena
Then due to \eqref{expEp} and \eqref{expEm}, in each sign case the right hand side of \eqref{EN1byE21EN2p}  determines both $E^+_{m,l}(z)$ and $E^-_{m,l}(z)$ 
 uniquely as formal Laurent series in $z$.   
 \qed

\noindent
{\it Remark.}\ 
The upper and the lower sign cases of \eqref{EN1byE21EN2p} give 
two expressions for $E_{m,l}^+(z_1^+)-E_{m,l}^-(z_1^-)$. 
It is instructive to derive their consistency condition. 
Equating them, we obtain 
 \be
&& \left(E_{m,l+1}^+(z_1^+)-E_{m,l+1}^-(z_1^-)\right)E_{l+1,l}^+(z_2^+)\nn\\
&&\qquad\qquad\qquad\qquad -\bar{b}^*(z)E_{l+1,l}^+(z_2^+)\left(
E_{m,l+1}^+(z_1^+)-E_{m,l+1}^+(z_1^-)\right)\nn\\
&&\qquad\qquad-\sum_{l+2\leq k\leq m-1}(E^+_{m,k}(z_1^+)-E^-_{m,k}(z_1^-))
\bar{c}^*(z,P_{l+1,k})(L_l^{m+}(z_2^+)^{-1})_{kl}K^+_l(z_2^+)\\
&&=\left(E_{m,l+1}^+(z_1^+)-E_{m,l+1}^-(z_1^-)\right)E_{l+1,l}^-(z_2^-)\nn\\
&&\qquad\qquad\qquad\qquad -\bar{b}^*(z)E_{l+1,l}^-(z_2^-)\left(
E_{m,l+1}^+(z_1^+)-E_{m,l+1}^-(z_1^-)\right)\nn\\
&&\qquad\qquad-\sum_{l+2\leq k\leq m-1}(E^+_{m,k}(z_1^+)-E^-_{m,k}(z_1^-))
\bar{c}^*(z,P_{l+1,k})(L_l^{m-}(z_2^-)^{-1})_{kl}K^-_l(z_2^-).
\en
Hence the consistency condition is 
 \be
&& \left(E_{m,l+1}^+(z_1^+)-E_{m,l+1}^-(z_1^-)\right)
\left(E_{l+1,l}^+(z_2^+)-E_{l+1,l}^-(z_2^-)\right)\nn\\
&&-\bar{b}^*(z)\left(E_{l+1,l}^+(z_2^+)-E_{l+1,l}^-(z_2^-)\right)
\left(E_{m,l+1}^+(z_1^+)-E_{m,l+1}^+(z_1^-)\right)\nn\\
&&\hspace{-1cm}=\sum_{l+2\leq k\leq m-1}(E^+_{m,k}(z_1^+)-E^-_{m,k}(z_1^-))
\bar{c}^*(z,P_{l+1,k})\left((L_l^{m+}(z_2^+)^{-1})_{kl}K^+_l(z_2^+)
-(L_l^{m-}(z_2^-)^{-1})_{kl}K^-_l(z_2^-)\right).
\en
In particular, for $m=l+2$ we have
\bea
&& \left(E_{l+2,l+1}^+(z_1^+)-E_{l+2,l+1}^-(z_1^-)\right)
\left(E_{l+1,l}^+(z_2^+)-E_{l+1,l}^-(z_2^-)\right)\nn\\
&&=\bar{b}^*(z)\left(E_{l+1,l}^+(z_2^+)-E_{l+1,l}^-(z_2^-)\right)
\left(E_{l+2,l+1}^+(z_1^+)-E_{l+2,l+1}^+(z_1^-)\right)\lb{elelp1}
\ena 
for $1\leq l\leq N-2$. 
In Appendix \ref{hc2ellc} we identify these relations with  the commutation relations of the total elliptic currents of $U_{q,p}(\glnh)$.
This provides an example suggesting the injectivity of the $H$-algebra homomorphism $\Phi: \cU\to \cE$ given in the next subsection.  

\subsection{The half currents of $\cU$}

Let us define the basic half currents of $\cU$ as follows. For $1\leq j\leq N-1$
\bea
e^+_{j+1,j}(z)&=&
\frac{{a_{j+1,j}^*\tes(q^2)}}{(p^*;p^*)_\infty^3}\left(\sum_{m\geq 0} e_{j,-m}\frac{1}{1-q^{-2(P_{\al_j}-1)}p^{*m}}(zq^{j-c})^{m}\right.\nn\\
&&\left. \qquad\qquad\qquad -\sum_{m>0} e_{j,m}\frac{q^{2(P_{\al_j}-1)}p^{*m}}{1-q^{2(P_{\al_j}-1)}p^{*m}}(zq^{j-c})^{-m}\right),\\
f^+_{j,j+1}(z)&=&
\frac{a_{j,j+1}\te(q^2)}{(p;p)_\infty^3}\left(\sum_{m\geq 0} f_{j,-m}\frac{1}{1-q^{2((P+h)_{\al_j}-1)}p^{m}}(zq^{j})^{m}\right.\nn\\
&&\left. \qquad\qquad\qquad -\sum_{m>0} f_{j,m}\frac{q^{-2((P+h)_{\al_j}-1)}p^{m}}{1-q^{-2((P+h)_{\al_j}-1)}p^{m}}(zq^{j})^{-m}\right),\\
e^-_{j+1,j}(z)&=&q^{2P_{\bep_{j+1}}}e_{j+1,j}^+(zp^*q^c)q^{-2P_{\bep_j}} \nn\\
&=&
\frac{a_{j+1,j}^*\tes(q^2)}{(p^*;p^*)_\infty^3}\left(\sum_{m\geq 0} e_{j,-m}\frac{q^{-2(P_{\al_j}-1)}p^{*m}}{1-q^{-2(P_{\al_j}-1)}p^{*m}}(zq^{j})^{m}\right.\nn\\
&&\left. \qquad\qquad\qquad\qquad\qquad -\sum_{m>0} e_{j,m}\frac{1}{1-q^{2(P_{\al_j}-1)}p^{*m}}(zq^{j})^{-m}\right),\\
f^-_{j,j+1}(z)&=&q^{2(P+{h})_{\bep_j}}f_{j,j+1}^+(zpq^{-c})q^{-2(P+{h})_{\bep_{j+1}}}\nn\\
&=&
\frac{a_{j,j+1}\te(q^2)}{(p;p)_\infty^3}\left(\sum_{m\geq 0} f_{j,-m}\frac{q^{2((P+h)_{\al_j}-1)}p^{m}}{1-q^{2((P+h)_{\al_j}-1)}p^{m}}(zq^{j-c})^{m}\right.\nn\\
&&\left. \qquad\qquad\qquad\qquad -\sum_{m>0} f_{j,m}\frac{1}{1-q^{-2((P+h)_{\al_j}-1)}p^{m}}(zq^{j-c})^{-m}\right),
\ena
where $a^*_{j+1,j}$ and $a_{j,j+1}$ are constants given by
\bea
&&a^*_{j+1,j}=q^{-1}\frac{(p^*;p^*)_\infty}{(p^*q^2;p^*)_\infty} ,\qquad a_{j,j+1}=q^{-1}\frac{(p;p)_\infty}{(pq^{-2};p)_\infty}.
\ena
We then obtain
\be
&&e^+_{j+1,j}(z^+)-e^-_{j+1,j}(z^-)=\frac{a_{j+1,j}^*\tes(q^{2})}{(p^*;p^*)_\infty^3}e_j(zq^{j-c/2}),\\
&&f^+_{j,j+1}(z^-)-f^-_{j,j+1}(z^+)=\frac{a_{j,j+1}\te(q^{2})}{(p;p)_\infty^3}f_j(zq^{j-c/2}).
\en
Here we set $z^\pm=z q^{\pm c/2}$. 

Note that at $p= 0$ 
\bea
e^+_{j+1,j}(z)&=&
a_{j+1,j}^*(1-q^{2})\left(e_{j,0}\frac{1}{1-q^{-2(P_{\al_j}-1)}}+\sum_{m> 0} e_{j,-m}(zq^{j-c})^m \right),\\
e^-_{j+1,j}(z)&=&
a_{j+1,j}^*(1-q^{2})\left(e_{j,0}\frac{q^{-2(P_{\al_j}-1)}}{1-q^{-2(P_{\al_j}-1)}}-\sum_{m\geq 0} e_{j,m}(zq^{j})^{-m}\right),\\
f^+_{j,j+1}(z)&=&
a_{j,j+1}(1-q^{2})\left(f_{j,0}\frac{1}{1-q^{2((P+h)_{\al_j}-1)}}+\sum_{m\geq 0} f_{j,-m}(zq^{j})^m\right),\\
f^-_{j,j+1}(z)&=&
a_{j,j+1}(1-q^{2})\left(f_{j,0}\frac{q^{2((P+h)_{\al_j}-1)}}{1-q^{2((P+h)_{\al_j}-1)}}-\sum_{m\geq 0} f_{j,m}(zq^{j-c})^{-m}\right).
\ena

Noting  \eqref{defenfngl} and the expansion formula 
\bea
\frac{\te(q^{2s}z)(p;p)^3_\infty}{\te(q^{2s})\te(z)}
&=&\sum_{n\in \Z}\frac{1}{1-q^{2s}p^n}z^{n}\nn\\
&=&\sum_{l\in \Z_{\geq 0}}\left(\frac{q^{2sl}}{1-{p^l}{z}}-
\frac{q^{-2s(l+1)}p^{l+1}/z}{1-p^{l+1}/z}\right),\lb{zpositiveexp}
\ena
for $|p|<|z|<1$,
 one can express the basic half currents as follows.
\begin{prop}\lb{inthc}
\bea
&&e_{j+1,j}^+(z)=
a_{j+1,j}^*\oint_{C^*} \frac{dz'}{2\pi i z'}
e_{j}(z')\frac{\tes(z q^{j-c}q^{2(1-
P_{\al_j})}/z')\tes(q^2)}
{\tes(zq^{j-c}/z')
\tes(q^{2(P_{\al_j}-1)})},
\\
&&f_{j,j+1}^+(z)
=a_{j,j+1}\oint_C \frac{dz'}{2\pi i z'}
f_{j}(z')\frac{\te(zq^{j}q^{2( (P+h)_{\al_j}-1)}/z')\te(q^2)}{\te(zq^{j}/z')\te(q^{2((P+h)_{\al_j}-1)})},
\end{eqnarray}
where $C^*: |q^{j-c}z|<|z'|<|p^{*-1}q^{j-c}z|$, $C: |q^{j}z|<|z'|<|p^{-1}q^{j}z|$. 
\end{prop}

\begin{prop}\cite{JKOS,KK03}\lb{KEF}~~~The basic half currents $e_{j+1,j}^+(z), 
f_{j,j+1}^+(z) \ (j\in I)$ and
$k^+_l(z)\ (1\leq l\leq N)$ satisfy the following relations. 
\begin{eqnarray}
%
%
%
%
&&k_{j+1}^+(z_1)^{-1}e_{j+1,j}^+(z_2)k_{j+1}^+(z_1)
=e_{j+1,j}^+(z_2)\frac{1}{\bar{b}_{}^*(z_1/z_2)}
-e_{j+1,j}^+(z_1)\frac{c_{}^*(z,P_{j,j+1})}{\bar{b}_{}^*(z)},\lb{klelj}\\
%
%
%
%
&&k_{j+1}^+(z_1)f_{j,j+1}^+(z_2)
k_{j+1}^+(z_1)^{-1}=
\frac{1}{\bar{b}_{}(z)}f_{j,j+1}^+(z_2)-
\frac{\bar{c}_{}(z,(P+h)_{j,j+1})}{
\bar{b}_{}(z)}f_{j,j+1}^+(z_1),
\lb{klfjl}\\
%
%
%
%
&&\frac{1}{\bar{b}^*(1/z)}e^+_{j+1,j}(z_1)e^+_{j+1,j}(z_2)-e^+_{j+1,j}(z_2)^2\frac{{c}^*(1/z,P_{j,j+1}-2)}{\bar{b}^*(1/z)}\nn\\
&&\qquad\qquad=\frac{1}{\bar{b}^*(z)}
e^+_{j+1,j}(z_2)e^+_{j+1,j}(z_1)
-e^+_{j+1,j}(z_1)^2\frac{{c}^*(z,P_{j,j+1}-2)}{\bar{b}^*(z)},\lb{eljelj}
\ena
%
%
%
%
\bea
&&\frac{1}{\bar{b}(z)}f^+_{j,j+1}(z_1)f^+_{j,j+1}(z_2)-f^+_{j,l}(z_1)^2\frac{\bar{c}(z,(P+h)_{j,j+1}-2)}{\bar{b}(z)}\nn\\
&&\qquad\qquad=\frac{1}{\bar{b}(1/z)}
f^+_{j,j+1}(z_2)f^+_{j,j+1}(z_1)-
f^+_{j,l}(z_2)^2\frac{\bar{c}(1/z,(P+h)_{j,j+1}-2)}{\bar{b}(1/z)},
\lb{fjlfjl}
\\
%
%
%
%
%
%
%
%
%
%
%
%
&&[e^+_{j+1,j}(z_1),f^+_{j,j+1}(z_2)]
=k^+_{j}(z_2)k^+_{j+1}(z_2)^{-1}
\frac{\bar{c}^*(z,P_{j,j+1}-1)}{\bar{b}^*(z)} \nn \\
&&\qquad\qquad\qquad\qquad\qquad-k^+_{j+1}(z_1)^{-1}k^+_{j}(z_1)
\frac{\bar{c}(z,(P+h)_{j,j+1}-1)}{\bar{b}(z)},\lb{elfjl}
%
%
%
\end{eqnarray}
where $z=z_1/z_2$. 

\end{prop}

\noindent
{\it Proof.} Direct calculation using Proposition \ref{inthc}, 
the  relations in Definition \ref{defUqpgl} and \eqref{kjkja}-\eqref{ffthetajjp1}. 
\qed

For $1\leq j\leq N-1$, let us consider the subalgebra $U^{(j)}_{q,p}(\glth)$ of  $\cU$ 
generated by $e_{j}(z), f_{j}(z), k_j^+(z)$, $k_{j+1}^+(z), q^{\pm c/2}, \widehat{d}$.  
Let us define the $L$-operator by the associated basic half currents by
\be
\cL_j^{\pm}(z)&=&\mat{1&f_{j,j+1}^\pm(z)\cr 0&1\cr}\mat{k_j^\pm(z)&0\cr 0&k_{j+1}^\pm(z)}
\mat{1&0\cr e_{j+1,j}^\pm(z)&1\cr}\\
&=&\mat{k_j^\pm(z)+f_{j,j+1}^\pm(z)k_{j+1}^\pm(z)e_{j+1,j}^\pm(z)&f_{j,j+1}^\pm(z)
k_{j+1}^\pm(z)\cr
k^\pm_{j+1}(z)e_{j+1,j}^\pm(z)&k^\pm_{j+1}(z)\cr}.
\en
Then comparing the relations in Proposition \ref{KEF} and those of the basic Gauss components of $L^\pm(z)$ in $E_{q,p}(\glth)$ in Sec.\ref{relBasicGC},
 we obtain the following. 
\begin{thm}\lb{gl2hom}\cite{JKOS,KK03}
For each $j$, the $L$-operators  $\cL_j^{\pm}(z)$ satisfy the same $RLL$-relations \eqref{mRLLlm}-\eqref{RLLpmlm} 
at $l=j, m=j+1$ replacing $L^{m\pm}_l(z)$ with $\cL_j^\pm(z)$.  Hence the following map gives a surjective $H$-algebra homomorphism.
\bea
&&\Phi^{(j)}\ : \ U^{(j)}_{q,p}(\glth)\to E_{q,p}(\glth),\\
&&e^+_{j+1,j}(z)\mapsto E^+_{j+1,j}(z),\quad f^+_{j,j+1}(z)\mapsto F^+_{j,j+1}(z),\quad k^+_{j}(z)\mapsto K^+_j(z)\quad k^+_{j+1}(z)\mapsto K^+_{j+1}(z).\nn
\ena
\end{thm}

Now let us consider the canonical extension of the map $\Phi^{(j)}$ to $\Phi: \cU\to \cE$ by 
\be
&&e^+_{j+1,j}(z)\mapsto E^+_{j+1,j}(z),\quad f^+_{j,j+1}(z)\mapsto F^+_{j,j+1}(z),\quad k^+_{l}(z)\mapsto K^+_l(z)\qquad (j\in I,\ 1\leq l\leq N).
\en 

\begin{thm}
 $\Phi$ gives an isomorphism 
as a topological $H$-algebra over $\FF[[p]]$.   
\end{thm}

\noindent
{\it Proof.}\ 
1) Surjectivity:  
From Theorem \ref{L_l^m} with $l=j, m=j+1$ the basic Gauss components $E^+_{j+1,j}(z)$, $F^+_{j,j+1}(z)$, $K^+_{j}(z)$ and $\hd$ generate the subalgebra $E_{q,p}(\glth)$. 
From Lemma \ref{recursion}  the $RLL$ relations allows us to  construct the other Gauss components $E^\pm_{k,j}(z)$,  $F^\pm_{j,k}(z)$  $( 3\leq j+2\leq k \leq N)$  recursively from the basic ones $E^\pm_{j+1,j}(z), F^\pm_{j,j+1}(z)$, $K^\pm_j(z), K^\pm_{j+1}(z)$.
Then the surjectivity follows from Theorem \ref{gl2hom}. 

\noindent
2) Injectivity: 
Let $(\varphi^q_{\la,k},V)$ be a highest weight representation  
of $\cU_q=U_q(\glnh)$ with the highest weight $\la$ and the level $k$.  We extend $\varphi^q_{\la,k}$ to the dynamical representation $(\varphi^{q,p}_{\la,k}=\varphi^q_{\la,k}\circ \phi_p, \hV_{\FF[[p]]})$ of $\cU$ 
as in Corollary \ref{dynamicalrep}. Then we define 
$\widetilde{\varphi}_{\lambda,k}^{q,p}$ : $\cE\ {\to}\ \End_\C \widehat{V}_{\FF[[p]]}$ by 
$\widetilde{\varphi}_{\lambda,k}^{q,p}(A)=\varphi^{q,p}_{\la,k}(a)$ for $A=\Phi(a), a\in \cU$. 
Let $\pi_p: \cU\to \cU/p\cU$ be the canonical projection. From the remark above Corollary \ref{dynamicalrep}, we  also have the 
corresponding canonical projection $\pi_{p}:\End_\C\widehat{V}_{\FF[[p]]}\to \End_\C({V_\FF}\otimes V_Q)$. 
We then consider the following diagram.  
\be
\mmatrix{
\cU&\xrightarrow{\ \Phi\ }&\cE&\xrightarrow{ \widetilde{\varphi}^{q,p}_{\la,k} }\quad &\cD_{H,\widehat{V}_{\FF[[p]]}}\\
\pi_p\downarrow\quad&\stackrel{\phi_p}{\searrow}&&\stackrel{\varphi^q_{\la,k}}{\nearrow}&\cr
\cU/p\cU& &(\FF[[p]]\otimes_\C \cU_q)\sharp\C[\cR_Q] & &\pi_{p}\downarrow\quad\\
\mbox{\rotatebox{90}{$\cong$}}&\stackrel{\pi_p}{\swarrow} &&\ &\cr
(\FF\otimes_{\C}\cU_q)\sharp \C[\cR_{\cQ}]&&\xrightarrow{\varphi^q_{\la,k}}&&\cD_{H,\widehat{V}_\FF}\cr}
\en

\begin{lem}\lb{piphi}
\

\begin{itemize}
\item[(i)] $\displaystyle{\pi_{p}\circ \phi_p=\pi_{p}}$
\item[(ii)] $\displaystyle{\pi_{p}\circ  \varphi_{\lambda,k}^{q}=
 \varphi_{\lambda,k}^{q}\circ\pi_{p}}$ on $(\FF[[p]]\otimes_\C \cU_q)\sharp\C[\cR_Q]$
\end{itemize}
\end{lem}
\noindent
{\it Proof.}\  $(i)$ follows from $u^\pm_{\vep_i}(z,p)= 1$ at $p= 0$.

$(ii)$ follows from $(i)$ and $\widehat{V}_{\FF[[p]]}=(\FF[[p]]\otimes_\C V)\otimes V_Q$. 
\qed

\begin{lem}\lb{KerinpU}
${\rm Ker}\ \Phi \subset p \cU$. 
\end{lem}
\noindent
{\it Proof.}\ 
Assume $\Ker \Phi\not\subset p \cU$. Then there exists a non zero element $a\in\Ker \Phi$ such that 
$\pi_{p}(a)\not=0$. Then from Lemma \ref{piphi} for any  level-$k$ highest weight representation $\varphi^q_{\la,k}$ of $\cU_q$, we have
\bea
0&=&\pi_{p}\circ\widetilde{\varphi}^{q,p}_{\la,k}\circ \Phi(a)=\pi_{p}\circ \varphi^q_{\la,k}\circ\phi_p(a)
=\varphi^q_{\la,k}\circ\pi_{p}(a).\lb{varphipi0}
\ena
This contradicts the fact ${\bigcap_{\la,k}\Ker \varphi^q_{\la,k}}=0$ given in \cite{DF}.
\qed

\noindent
{\it Proof of the injectivity.}\ 
Let us assume $\Ker\Phi\not=0$. Let $a\not=0\in \Ker\Phi$. Then from Lemma \ref{KerinpU} there exists $\widetilde{a}\in \cU$ such that $a=p^n\widetilde{a}$ for some positive integer $n$ and $\pi_{p}(\widetilde{a})\not=0\in \cU_q$. Then the same argument as \eqref{varphipi0} 
yields for any $\varphi^q_{\la,k}$ 
\be
0&=&\pi_{p}\circ\varphi^{q,p}_{\la,k}\circ \Phi(a)
=p^n\varphi^q_{\la,k}\circ\pi_{p}(\widetilde{a}).
\en
This again contradicts ${\bigcap_{\la,k}\Ker \varphi^q_{\la,k}}=0$.

\qed

\section*{Acknowledgements}
The author would like to thank Masatoshi Noumi for valuable discussions. He is also grateful to Ivan Cherednik, Giovanni Felder, Vassily Gorbounov, Michio Jimbo, Anatol Kirillov, Erik Koelink, Stanislav Pakulyak and Hjalmar Rosengren for useful conversations and their interests.  He is supported by the Grant-in -Aid for Scientific Research (C) 26400046 JSPS, Japan.

\appendix
\setcounter{equation}{0}
\begin{appendix}

\section{Quantum Affine Algebra $U_q(\glnhbig)$}\lb{Uqgh}

\begin{dfn}\lb{defUqglnh}
The quantum affine algebra $U_q(\glnh)$ is a 
topological algebra over $\C$ generated by 
 $k^\pm_{i,m},\ x^\pm_{j,n}, d \ (1\leq i\leq N,\ 1\leq j\leq N-1,\ 
 m\in\Z_{\geq 0},\ n\in\Z)$ and the central element $q^{\pm c/2}$. 
The defining relations are conveniently written in terms of the generating functions called 
the Drinfeld currents : 
\be
&&k^\pm_{0,i}(z)=\sum_{m\in \Z_{\geq 0}}k^\pm_{i,m}z^{\pm m},\\
&&x_j^\pm(z)=\sum_{n\in \Z}x^\pm_{j,n} z^{-n}.
\en
The relations are given by 
\be
&&[d,k^\pm_i(z)]=\pm z\frac{\partial }{\partial z} k^\pm_i(z),\qquad 
[d,x^\pm_j(z)]=- z\frac{\partial }{\partial z} x^\pm_j(z),\\
&&k^+_{i,0}k^-_{i,0}=1=k^-_{i,0}k^+_{i,0},\lb{Uqk0k0}\\
&&k^\pm_{0,j}(z_1)k^\pm_{0,l}(z_2)=k^\pm_{0.l}(z_2)k^\pm_{0,j}(z_1),\lb{Uqkk}
\\
&&k^-_{0,j}(z_1)k^+_{0,j}(z_2)=\frac{(q^{c+2}z_2/z_1,q^{2N}q^{c-2}z_2/z_1,q^{2N}q^{-c}z_2/z_1,q^{-c}z_2/z_1;q^{2N})_\infty}{
(q^{-c+2}z_2/z_1,q^{2N}q^{-c-2}z_2/z_1,q^{2N}q^{c}z_2/z_1,q^{c}z_2/z_1;q^{2N})_\infty} k^+_{0,j}(z_2)k^-_{0,j}(z_1),\\
&&k^-_{0,j}(z_1)k^+_{0,l}(z_2)=\frac{(q^{2N}q^{c+2}z_2/z_1,q^{2N}q^{c-2}z_2/z_1,q^{2N}q^{-c}z_2/z_1,q^{2N}q^{-c}z_2/z_1;q^{2N})_\infty}{
(q^{2N}q^{-c+2}z_2/z_1,q^{2N}q^{-c-2}z_2/z_1,q^{2N}q^{c}z_2/z_1,q^{2N}q^{c}z_2/z_1;q^{2N})_\infty}\nn\\
&&\hspace{8cm}\times
k^+_{0,l}(z_2)k^-_{0,j}(z_1)\quad (j<l),
\en
\be
&&k^-_{0,j}(z_1)k^+_{0,l}(z_2)=\frac{(q^{c+2}z_2/z_1,q^{c-2}z_2/z_1,q^{-c}z_2/z_1,q^{-c}z_2/z_1;q^{2N})_\infty}{
(q^{-c+2}z_2/z_1,q^{-c-2}z_2/z_1,q^{c}z_2/z_1,q^{c}z_2/z_1;q^{2N})_\infty}
k^+_{0,l}(z_2)k^-_{0,j}(z_1)\qquad (j>l),
\en
\be
&&k^+_{0,j}(z_1)x^+_j(z_2)k^+_{0,j}(z_1)^{-1}=q^{-1}\frac{1-q^{-c+j}z_1/z_2}{1-q^{-c-2+j}z_1/z_2}x^+_j(z_2),\\
&&k^+_{0,j+1}(z_1)x^+_j(z_2)k^+_{0,j+1}(z_1)^{-1}=q\frac{1-q^{-c+j}z_1/z_2}{1-q^{-c+2+j}z_1/z_2}x^+_j(z_2),\\
&&k^+_{0,j}(z_1)x^-_j(z_2)k^+_{0,j}(z_1)^{-1}=q\frac{1-q^{-2+j}z_1/z_2}{1-q^{j}z_1/z_2}x^-_j(z_2),\\
&&k^+_{0,j+1}(z_1)x^-_j(z_2)k^+_{0,j+1}(z_1)^{-1}=q^{-1}\frac{1-q^{2+j}z_1/z_2}{1-q^{j}z_1/z_2}x^-_j(z_2),\\
&&k^+_{0,l}(z_1)x^\pm_j(z_2)k^+_{0,l}(z_1)^{-1}=x^\pm_j(z_2)\qquad (l\not=j,j+1), 
\en
\be
&&k^-_{0,j}(z_1)^{-1}x^+_j(z_2)k^-_{0,j}(z_1)=q^{-1}\frac{1-q^{2-j}z_2/z_1}{1-q^{-j}z_2/z_1}x^+_j(z_2),\\
&&k^-_{0,j+1}(z_1)^{-1}x^+_j(z_2)k^-_{0,j+1}(z_1)=q\frac{1-q^{-2-j}z_2/z_1}{1-q^{-j}z_2/z_1}x^+_j(z_2),\\
&&k^-_{0,j}(z_1)^{-1}x^-_j(z_2)k^-_{0,j}(z_1)=q\frac{1-q^{c-j}z_2/z_1}{1-q^{c+2-j}z_2/z_1}x^-_j(z_2),\\
&&k^-_{0,j+1}(z_1)^{-1}x^-_j(z_2)k^-_{0,j+1}(z_1)=q^{-1}\frac{1-q^{c-j}z_2/z_1}{1-q^{c-2-j}z_2/z_1}x^-_j(z_2),\\
&&k^-_{0,l}(z_1)^{-1}x^\pm_j(z_2)k^-_{0,l}(z_1)=x^\pm_j(z_2)\qquad (l\not=j,j+1), \lb{Uqkf}
\en
\be
&&z_1(1-q^{\pm 2}z_2/z_1)
x_j^\pm(z_1)x_j^\pm(z_2)=-z_2 (1-q^{\pm 2}z_1/z_2) x_j^\pm(z_2)x_j^\pm(z_1),\\
&&z_1(1-q^{\mp 1}z_2/z_1)
x_j^\pm(z_1)x_{j+1}^\pm(z_2)=-z_2 (1-q^{\mp 1}z_1/z_2) x_{j+1}^\pm(z_2)x_j^\pm(z_1),\\
&&x_j^\pm(z_1)x_l^\pm(z_2)= x_l^\pm(z_2)x_j^\pm(z_1)\qquad (l\not=j,j+1),
\en
\be
&&[x_i^+(z_1),x_j^-(z_2)]=\frac{\delta_{i,j}}{q-q^{-1}}
\left(\delta\bigl(q^{-c}{z_1}/{z_2}\bigr)k^-_{0,i}(q^{-c/2}z_1)k^-_{0,i+1}(q^{-c/2}z_1)^{-1}\right.\nn\\
&&\hspace{5cm}\left.-\delta\bigl(q^{c}{z_1}/{z_2}\bigr)k^+_{0,i}(q^{-c/2}z_2)k^+_{0,i+1}(q^{-c/2}z_2)^{-1}
\right),\\
&&\left\{x^\pm_i(z_1)x^\pm_i(z_2)x^\pm_j(w)-(q+q^{-1})x^\pm_i(z_1)x^\pm_j(w)x^\pm_i(z_2)
+x^\pm_j(w)x^\pm_i(z_1)x^\pm_i(z_2)\right\}\nn\\
&&\hspace{5cm}+(z_1\leftrightarrow z_2)=0,\qquad |i-j|=1.  
\en
\end{dfn}

\section{Relations Among $u^\pm_{\vep_i}(z,p)$ and $u_j^\pm(z,p)$}\lb{EvHom} 

From the relations in Definition \ref{defUqglnh} the following commutation relations hold. 
\begin{lem}\lb{commuu}
Let us set $z=z_1/z_2$. 
\begin{eqnarray*}
&&\hspace{-1cm}k^+_{0,j}(z_1)u^-_{\vep_j}(q^jz_2,p)=\frac{(p^*q^{2}z,p^*q^{2N}q^{-2}z,pq^{2N}z,pz;p,q^{2N})_\infty}
{(pq^{2}z,pq^{2N}q^{-2}z,p^*q^{2N}z,p^*z;p,q^{2N})_\infty}u^-_{\vep_j}(q^jz_2,p)k^+_{0,j}(z_1),\\
&&\hspace{-1cm}k^+_{0,j}(z_1)u^-_{\vep_l}(q^lz_2,p)
=\frac{(p^*q^{2}z,p^*q^{-2}z,pz,pz;p,q^{2N})_\infty}
{(pq^{2}z,pq^{-2}z,p^*z,p^*z;p,q^{2N})_\infty}u^-_{\vep_l}(q^lz_2,p)k^+_{0,j}(z_1)\ (j<l),
\\
&&\hspace{-1cm}k^+_{0,j}(z_1)u^-_{\vep_l}(q^lz_2,p)
=\frac{(p^*q^{2N}q^{2}z,p^*q^{2N}q^{-2}z,pq^{2N}z,pq^{2N}z;p,q^{2N})_\infty}
{(pq^{2N}q^{2}z,pq^{2N}q^{-2}z,p^*q^{2N}z,p^*q^{2N}z;p,q^{2N})_\infty}\nn\\
&&\hspace{5cm}\times u^-_{\vep_l}(q^lz_2,p)k^+_{0,j}(z_1)\qquad (j>l),
\\
&&k^+_{0,j}(z_1)u_j^-(q^jz_2,p)=\frac{(pq^{-2}z;p)_\infty(p^*z;p)_\infty}{(p^*q^{-2}z;p)_\infty(pz;p)_\infty}
u_j^-(q^jz_2,p)k^+_{0,j}(z_1),\\
&&k^+_{0,j+1}(z_1)u_j^-(q^jz_2,p)=\frac{(pq^{2}z;p)_\infty(p^*z;p)_\infty}{(p^*q^{2}z;p)_\infty(pz;p)_\infty}
u_j^-(q^jz_2,p)k^+_{0,j+1}(z_1),\\
&&k^+_{0,l}(z_1)u_j^-(q^jz_2,p)=u_j^-(q^jz_2,p)k^+_{0,l}(z_1)\qquad (l\not=j,j+1),
\en
\be
&&\hspace{-1cm}u^+_{\vep_j}(q^jz_1,p)k^-_{0,j}(z_2)=\frac{(p^*q^{2}z,p^*q^{2N}q^{-2}z,pq^{2N}z,pz;p^*,q^{2N})_\infty}
{(pq^{2}z,pq^{2N}q^{-2}z,p^*q^{2N}z,p^*z;p^*,q^{2N})_\infty}k^-_{0,j}(z_2)u^+_{\vep_j}(q^jz_1,p),\\
&&\hspace{-1cm}u^+_{\vep_l}(q^lz_1,p)k^-_{0,j}(z_2)
=\frac{(p^*q^{2N}q^{2}z,p^*q^{2N}q^{-2}z,pq^{2N}z,pq^{2N}z;p^*,q^{2N})_\infty}
{(pq^{2N}q^{2}z,pq^{2N}q^{-2}z,p^*q^{2N}z,p^*q^{2N}z;p^*,q^{2N})_\infty}\nn\\
&&\hspace{5cm}\times k^-_{0,j}(z_2)u^+_{\vep_l}(q^lz_1,p)\qquad (j<l),\\
&&\hspace{-1cm}u^+_{\vep_l}(q^lz_1,p)k^-_{0,j}(z_2)
=\frac{(p^*q^{2}z,p^*q^{-2}z,pz,pz;p^*,q^{2N})_\infty}
{(pq^{2}z,pq^{-2}z,p^*z,p^*z;p^*,q^{2N})_\infty}
k^-_{0,j}(z_2)u^+_{\vep_l}(q^lz_1,p)\ (j>l),\nn\\
&&u_j^+(q^jz_1,p)k^-_{0,j}(z_2)=\frac{(p^*q^{2}z;p^*)_\infty(pz;p^*)_\infty}{(pq^{2}z;p^*)_\infty(p^*z;p^*)_\infty}
k^-_{0,j}(z_2)u_j^+(q^jz_1,p),\\
&&u_j^+(q^jz_1,p)k^-_{0,j+1}(z_2)=\frac{(p^*q^{-2}z;p^*)_\infty(pz;p^*)_\infty}{(pq^{-2}z;p^*)_\infty(p^*z;p^*)_\infty}
k^-_{0,j+1}(z_2)u_j^+(q^jz_1,p)
,\\
&&u_j^+(q^jz_1,p)k^-_{0,l}(z_2)=k^-_{0,l}(z_2)u_j^+(q^jz_1,p)\qquad (l\not=j,j+1),\\
&&u^{\pm}_{\vep_j}(z_1,p)u^{\pm}_{\vep_l}(z_2,p)=u^{\pm}_{\vep_l}(z_2,p)u^{\pm}_{\vep_j}(z_1,p)\quad (\forall j,l),\\
&&\hspace{-1cm}u^+_{\vep_j}(q^jz_1,p)u^-_{\vep_j}(q^jz_2,p)=\frac{(p^{*2}q^{c+2}z,p^{*2}q^{2N}q^{c-2}z,pp^*q^{2N}q^cz,pp^*q^cz;p,p^*,q^{2N})_\infty}
{(pp^*q^{c+2}z,pp^*q^{2N}q^{c-2}z,p^{*2}q^{2N}q^cz,p^{*2}q^cz;p,p^*,q^{2N})_\infty}\nn\\
&&\hspace{5cm}\times u^-_{\vep_j}(q^jz_2,p)(z_2)u^+_{\vep_j}(q^jz_1,p),\\
&&\hspace{-1cm}u^+_{\vep_l}(q^lz_1,p)u^-_{\vep_j}(q^jz_2,p)
=\frac{(p^{*2}q^{2N}q^{c+2}z,p^{*2}q^{2N}q^{c-2}z,pp^*q^{2N}q^cz,pp^*q^{2N}q^cz;p,p^*,q^{2N})_\infty}
{(pp^*q^{2N}q^{c+2}z,pp^*q^{2N}q^{c-2}z,p^{*2}q^{2N}q^cz,p^{*2}q^{2N}q^cz;p,p^*,q^{2N})_\infty}\nn\\
&&\hspace{5cm}\times u^-_{\vep_j}(q^jz_2,p)u^+_{\vep_l}(q^lz_1,p)\qquad (j<l),
\\
&&\hspace{-1cm}u^+_{\vep_l}(q^lz_1,p)u^-_{\vep_j}(q^jz_2,p)
=\frac{(p^{*2}q^{c+2}z,p^{*2}q^{c-2}z,pp^*q^cz,pp^*q^cz;p,p^*,q^{2N})_\infty}
{(pp^*q^{c+2}z,pp^*q^{c-2}z,p^{*2}q^cz,p^{*2}q^cz;p,p^*,q^{2N})_\infty}\nn\\
&&\hspace{5cm}\times u^-_{\vep_j}(q^jz_2,p)u^+_{\vep_l}(q^lz_1,p)\qquad (j>l),\\
&&u_j^+(q^jz_1,p)u^-_{\vep_j}(q^jz_2,p)=\frac{(p^{*}q^{c+2}z;p^*)_\infty (p^*q^cz;p)_\infty}{(p^*q^{c+2}z;p)_\infty
( p^{*}q^c z;p^*)_\infty}
u^-_{\vep_j}(q^jz_2,p)u_j^+(q^jz_1,p),
\en
\be
&&\hspace{-0.2cm}u_j^+(q^jz_1,p)u^-_{\vep_{j+1}}(q^{j+1}z_2,p)=\frac{(p^*q^{c-2}z;p^*)_\infty(p^*q^cz;p)_\infty}{(p^*q^{c-2}z;p)_\infty(p^*q^cz;p^*)_\infty}
u^-_{\vep_{j+1}}(q^{j+1}z_2,p)u_j^+(q^jz_1,p),
\\
&&u_j^+(q^jz_1,p)u^-_{\vep_l}(q^lz_2,p)=u^-_{\vep_l}(q^lz_2,p)u_j^+(q^jz_1,p)\qquad (l\not=j,j+1),\\
&&u_{\vep_j}^+(q^jz_1,p)u^-_{j}(q^jz_2,p)=\frac{(p^*q^{c-2}z;p)_\infty(p^*q^cz;p^*)_\infty}{(p^*q^{c-2}z;p^*)_\infty(p^*q^cz;p)_\infty}
u^-_{j}(q^jz_2,p)u_{\vep_j}^+(q^jz_1,p),
\\
&&\hspace{-0.2cm}u^+_{\vep_{j+1}}(q^{j+1}z_1,p)u_j^-(q^jz_2,p)=
\frac{(p^*q^{c+2}z;p)_\infty
( p^{*}q^c z;p^*)_\infty}{(p^{*}q^{c+2}z;p^*)_\infty (p^*q^cz;p)_\infty}
u_j^-(q^jz_2,p)u^+_{\vep_{j+1}}(q^{j+1}z_1,p),
\\
&&u^+_{\vep_l}(q^lz_1,p)u_j^-(q^jz_2,p)=u_j^-(q^jz_2,p)u^+_{\vep_l}(q^lz_1,p)\qquad (l\not=j,j+1),\\
&&u_j^+(z_1,p)u_j^-(z_2,p)=
\frac{(pq^{-c-2}z;p)_\infty}{(pq^{-c+2}z;p)_\infty}
\frac{(p^*q^{c+2}z;p^*)_\infty}{(p^*q^{c-2}z;p^*)_\infty}
u_j^-(z_2,p)u_j^+(z_1,p),\\
&&u_j^+(z_1,p)u_{j+1}^-(z_2,p)=
\frac{(pq^{-c+1}z;p)_\infty}{(pq^{-c-1}z;p)_\infty}
\frac{(p^*q^{c-1}z;p^*)_\infty}{(p^*q^{c+1}z;p^*)_\infty}
u_{j+1}^-(z_2,p)u_j^+(z_1,p),\\
&&u_j^+(z_1,p)u_{l}^-(z_2,p)=
u_{l}^-(z_2,p)u_j^+(z_1,p) \qquad (l\not=j,j+1),
\en
\be
&&u^+_{\vep_j}(q^{-c+j}z_1,p)x^+_j(z_2)=\frac{(p^*q^{-c+j}z;p^*)_\infty}{(p^*q^{-c-2+j}z;p^*)_\infty}x^+_j(z_2)u^+_{\vep_j}(q^{-c+j}z_1,p),\\
&&u^-_{\vep_j}(q^jz_1,p)x^+_j(z_2)=\frac{(pq^{-c-j}/z;p)_\infty}{(pq^{-c+2-j}/z;p)_\infty}x^+_j(z_2)u^-_{\vep_j}(q^jz_1,p),\\
&&u^+_{\vep_{j+1}}(q^{-c+j+1}z_1,p)x^+_j(z_2)=\frac{(p^*q^{-c+j}z;p^*)_\infty}{(p^*q^{-c+2+j}z;p^*)_\infty}x^+_j(z_2)u^+_{\vep_{j+1}}(q^{-c+j+1}z_1,p),\\
&&u^-_{\vep_{j+1}}(q^{j+1}z_1,p)x^+_j(z_2)=\frac{(pq^{-c-j}/z;p)_\infty}{(pq^{-c-2-j}/z;p)_\infty}x^+_j(z_2)u^-_{\vep_{j+1}}(q^{j+1}z_1,p),
\en
\be
&&u^\pm_{\vep_l}(z_1,p)x^+_j(z_2)=x^+_j(z_2)u^\pm_{\vep_l}(z_1,p)\qquad (l\not=j,j+1),\\
&&u_j^+(z_1,p)x_j^+(z_2)=
\frac{(p^*q^{2}z;p^*)_\infty}{(p^*q^{-2}z;p^*)_\infty}
x_j^+(z_2)u_j^+(z_1,p),\\
&&u_{j+1}^+(z_1,p)x_j^+(z_2)=
\frac{(p^*q^{-1}z;p^*)_\infty}{(p^*qz;p^*)_\infty}
x_j^+(z_2)u_{j+1}^+(z_1,p),\\
&&u_j^-(z_1,p)x_j^+(z_2)=
\frac{(pq^{-2-c}z;p)_\infty}{(pq^{2-c}z;p)_\infty}
x_j^+(z_2)u_j^-(z_1,p),\\
&&u_{j+1}^-(z_1,p)x_j^+(z_2)=
\frac{(pq^{1-c}z;p)_\infty}{(pq^{-1-c}z;p)_\infty}
x_j^+(z_2)u_{j+1}^-(z_1,p),\\
&&u_l^\pm(z_1,p)x_j^+(z_2)=x_j^+(z_2)u_l^\pm(z_1,p)\qquad (l\not=j,j+1),
\en
\be
&&u^+_{\vep_j}(q^{-c+j}z_1,p)x^-_j(z_2)=\frac{(p^*q^{-2+j}z;p^*)_\infty}{(p^*q^{j}z;p^*)_\infty}x^-_j(z_2)u^+_{\vep_j}(q^{-c+j}z_1,p),\\
&&u^-_{\vep_j}(q^jz_1,p)x^-_j(z_2)=\frac{(pq^{2-j}/z;p)_\infty}{(pq^{-j}/z;p)_\infty}x^-_j(z_2)u^-_{\vep_j}(q^jz_1,p),\\
&&u^+_{\vep_{j+1}}(q^{-c+j+1}z_1,p)x^-_j(z_2)=\frac{(p^*q^{2+j}z;p^*)_\infty}{(p^*q^{j}z;p^*)_\infty}x^-_j(z_2)u^+_{\vep_{j+1}}(q^{-c+j+1}z_1,p),\\
&&u^-_{\vep_{j+1}}(q^{j+1}z_1,p)x^-_j(z_2)=\frac{(pq^{-2-j}/z;p)_\infty}{(pq^{-j}/z;p)_\infty}x^-_j(z_2)u^-_{\vep_{j+1}}(q^{j+1}z_1,p),\\
&&u^\pm_{\vep_l}(z_1,p)x^-_j(z_2)=x^-_j(z_2)u^\pm_{\vep_l}(z_1,p)\qquad (l\not=j,j+1),
\en
\be
&&u_j^+(z_1,p)x_j^-(z_2)=
\frac{(p^*q^{-2+c}z;p^*)_\infty}{(p^*q^{2+c}z;p^*)_\infty}
x_j^-(z_2)u_j^+(z_1,p),\\
&&u_{j+1}^+(z_1,p)x_j^-(z_2)=
\frac{(p^*q^{1+c}z;p^*)_\infty}{(p^*q^{-1+c}z;p^*)_\infty}
x_j^-(z_2)u_{j+1}^+(z_1,p),\\
&&u_j^-(z_1,p)x_j^-(z_2)=
\frac{(pq^{2}z;p)_\infty}{(pq^{-2}z;p)_\infty}
x_j^-(z_2)u_j^-(z_1,p),\\
&&u_{j+1}^-(z_1,p)x_j^-(z_2)=
\frac{(pq^{-1}z;p)_\infty}{(pqz;p)_\infty}
x_j^-(z_2)u_{j+1}^-(z_1,p),\\
&&u_l^\pm(z_1,p)x_j^-(z_2)=x_j^-(z_2)u_l^\pm(z_1,p)\qquad (l\not=j,j+1).
\end{eqnarray*}
\end{lem}

\section{From the Basic Gauss Components to the Elliptic Currents 
}\lb{hc2ellc}

\subsection{Relations among the basic Gauss components }\lb{relBasicGC}
Let us consider the $m=j+1, l=j$ case of the $L^{m\pm}_j(z)$ and 
the subalgebra  $E_{q,p}(\widehat{\gl}_2)$ generated by them. 
\be
L_j^{j+1\pm}(z)&=&\mat{1&F_{j,j+1}^\pm(z)\cr 0&1\cr}\mat{K_j^\pm(z)&0\cr 0&K_{j+1}^\pm(z)}
\mat{1&0\cr E_{j+1,j}^\pm(z)&1\cr}\\
&=&\mat{K_j^\pm(z)+F_{j,j+1}^\pm(z)K_{j+1}^\pm(z)E_{j+1,j}^\pm(z)&F_{j,j+1}^\pm(z)
K_{j+1}^\pm(z)\cr
K^\pm_{j+1}(z)E_{j+1,j}^\pm(z)&K^\pm_{j+1}(z)\cr}.
\en
\begin{eqnarray}
L_j^{j+1\pm}(z)^{-1}=\left(\begin{array}{cc}
K_j^\pm(z)^{-1}&-K_j^\pm(z)^{-1}F_{j,j+1}^\pm(z)\\
-E_{j+1,j}^\pm(z)K_j^\pm(z)^{-1}&E_{j+1,j}^\pm(z)K_j^\pm(z)^{-1}
F_{j,j+1}^\pm(z)+K_{j+1}^\pm(z)^{-1}
\end{array}\right).
\end{eqnarray}

From $(j,j),(j,j)$ of \eqref{mRLLlminvinv} 
and $(j+1,j+1),(j+1,j+2)$ of \eqref{mRLLlm}
\begin{eqnarray}
K_l^\pm(z_1)K_l^\pm(z_2)=\rho(z)
K_l^\pm(z_2)K_l^\pm(z_1)\qquad~(l=j,j+1),\lb{KpKp}
\end{eqnarray}
where $z=z_1/z_2$. 
Similarly, from $(j,j),(j,j)$ of \eqref{RLLpmlminvinv}
and $(j+1,j+1),(j+1,j+1)$ of \eqref{RLLpmlm}
\begin{eqnarray}
K_l^\pm(z_1)K_l^\mp(z_2)=
\frac{\rho^{\pm*}(zq^{\mp{c}})}{\rho^{\pm}(zq^{\pm{c}})}
K_l^\mp(z_2)K_l^\pm(z_1),~(l=j,j+1).\lb{KpKm}
\end{eqnarray}
From $(j+1,j),(j+1,j)$ of \eqref{mRLLlminv}
\begin{eqnarray}
K_j^\pm(z_1)K_{j+1}^\pm(z_2)=\rho(z)
\frac{\bar{b}(z)}{\bar{b}^*(z)}
K_{j+1}^\pm(z_2)K_j^\pm(z_1).\lb{Kp1Kp2}
\end{eqnarray}
Inparticular, in the limit $z_1\to z_2$,
 we have
 \bea
 &&K_j^\pm(z_2)K_{j+1}^\pm(z_2)=\frac{(p;p)^3_\infty}{(p^*;p^*)^3_\infty}
\frac{\tes(q^2)}{\te(q^2)}
K_{j+1}^\pm(z_2)K_j^\pm(z_2).\lb{Kpm2Kpm2}
 \ena
From $(j+1,j),(j+1,j)$ of \eqref{RLLpmlminv}
\begin{eqnarray}
K_{j+1}^\pm(z_1)K_j^\mp(z_2)=
\frac{\rho^{\pm*}(zq^{\mp{c}})}
{\rho^{\pm}(zq^{\pm {c}})}\frac{\bar{b}^*(zq^{\mp c})}{\bar{b}(zq^{\pm c})}
K_j^\mp(z_2)K_{j+1}^\pm(z_1).\lb{Kp2Km1}
\end{eqnarray}
From $(j+1,j+1),(j+1,j)$ of \eqref{mRLLlm}, we obtain
\bea
K^\pm_{j+1}(z_1)^{-1}E_{j+1,j}^\pm(z_2)K_{j+1}^\pm(z_1)=E_{j+1,j}^\pm(z_2)\frac{1}{\bar{b}^*(z)}
-E_{j+1,j}^\pm(z_1)\frac{c^*(z,P_{j,j+1})}{\bar{b}^*(z)}.\lb{KpEpKp}
\ena
From $(j+1,j+1),(j+1,j)$ of \eqref{RLLpmlm}, we obtain
\bea
K^\pm_{j+1}(z_1)^{-1}E_{j+1,j}^\mp(z_2)K_{j+1}^\pm(z_1)=E_{j+1,j}^\mp(z_2)\frac{1}{\bar{b}^*(zq^{\mp c})}
-E_{j+1,j}^\pm(z_1)\frac{c^*(zq^{\mp c},P_{j,j+1})}{\bar{b}^*(zq^{\mp c})}
.\lb{KpEmKp}
\ena
Similarly, from $(j+1,j),(j+1,j+1)$ of \eqref{mRLLlm}, 
\bea
K_{j+1}^\pm(z_1)F_{j,j+1}^\pm(z_2)K_{j+1}^\pm(z_1)^{-1}=
\frac{1}{\bar{b}(z)}F_{j,j+1}^\pm(z_2)-
\frac{\bar{c}(z,(P+h)_{j,j+1})}{\bar{b}(z)}F_{j,j+1}^\pm(z_1).\lb{KpFpKp}
\ena
From $(j+1,j),(j+1,j+1)$ of \eqref{RLLpmlm},
\bea
K_{j+1}^\pm(z_1)F_{j,j+1}^\mp(z_2)K_{j+1}^\pm(z_1)^{-1}=
\frac{1}{\bar{b}(zq^{\pm c})}F_{j,j+1}^\mp(z_2)-
\frac{\bar{c}(zq^{\pm c},(P+h)_{j,j+1})}{\bar{b}(zq^{\pm c})}F_{j,j+1}^\pm(z_1).\lb{KpFmKp}
\ena
From $(j+1,j),(j,j)$ of \eqref{mRLLlminvinv} 
\begin{eqnarray}
K_j^\pm(z_2)^{-1}E_{j+1,j}^\pm(z_1)K_j^\pm(z_2)=
\frac{1}{\bar{b}^*(z)}E_{j+1,j}^\pm(z_1)-
\frac{\bar{c}^*(z,P_{j,j+1})}{\bar{b}^*(z)}E_{j+1,j}^\pm(z_2).\lb{Kp1EpKp1}
\end{eqnarray}
From $(j+1,j),(j,j)$ of \eqref{RLLpmlminvinv} 
\begin{eqnarray}
K_j^\mp(z_2)^{-1}E_{j+1,j}^\pm(z_1)K_j^\mp(z_2)=
\frac{1}{\bar{b}^*(zq^{\mp c})}E_{j+1,j}^\pm(z_1)-
\frac{\bar{c}^*(zq^{\mp c},P_{j,j+1})}{\bar{b}^*(zq^{\mp c})}E_{j+1,j}^\mp(z_2).\lb{Km1EpKm1}
\end{eqnarray}
From $(j,j),(j+1,j)$ of \eqref{mRLLlminvinv}
\begin{eqnarray}
K_j^\pm(z_2)F_{j,j+1}^\pm(z_1)K_j^\pm(z_2)^{-1}=
F_{j,j+1}^\pm(z_1)\frac{1}{\bar{b}(z)}
-F_{j,j+1}^\pm(z_2)\frac{\bar{c}(z,(P+h)_{j,j+1})}{\bar{b}(z)}
.\lb{Kp1FpKp1}
\end{eqnarray}
From $(j+1,j+1),(j,j)$ of \eqref{mRLLlm}, we obtain
\bea
K^\pm_{j+1}(z_2)^{-1}E_{j+1,j}^\pm(z_1)K_{j+1}^\pm(z_2)E_{j+1,j}^\pm(z_2)=
K^\pm_{j+1}(z_1)^{-1}E_{j+1,j}^\pm(z_2)K_{j+1}^\pm(z_1)E_{j+1,j}^\pm(z_1).\lb{KpEpKpEp}
\ena
From \eqref{KpEpKp} and \eqref{KpEpKpEp}, we obtain
\bea
&&E_{j+1,j}^\pm(z_1)E_{j+1,j}^\pm(z_2)\frac{1}{\bar{b}^*(1/z)}
-E_{j+1,j}^\pm(z_2)^2\frac{{c}^*(1/z,P_{j,j+1}-2)}{\bar{b}^*(1/z)}
\nn\\
&&=E_{j+1,j}^\pm(z_2)E_{j+1,j}^\pm(z_1)\frac{1}{\bar{b}^*(z)}
-E_{j+1,j}^\pm(z_1)^2\frac{{c}^*(z,P_{j,j+1}-2)}{\bar{b}^*(z)}
.\lb{EpEp}
\ena
From $(j+1,j+1),(j,j)$ of \eqref{RLLpmlm}, we obtain
\bea
K^\mp_{j+1}(z_2)^{-1}E_{j+1,j}^\pm(z_1)K_{j+1}^\mp(z_2)E_{j+1,j}^\mp(z_2)=
K^\pm_{j+1}(z_1)^{-1}E_{j+1,j}^\mp(z_2)K_{j+1}^\pm(z_1)E_{j+1,j}^\pm(z_1).\lb{KmEpKmEm}
\ena
From \eqref{KpEmKp} and \eqref{KmEpKmEm}, we obtain
\bea
&&E_{j+1,j}^\pm(z_1)E_{j+1,j}^\mp(z_2)\frac{1}{\bar{b}^*(1/(zq^{\mp c}))}
-E_{j+1,j}^\mp(z_2)^2\frac{{c}^*(1/(zq^{\mp c}),P_{1,2}-2)}{\bar{b}^*(1/(zq^{\mp c}))}
\nn\\
&&=E_{j+1,j}^\mp(z_2)E_{j+1,j}^\pm(z_1)\frac{1}{\bar{b}^*(zq^{\mp c})}
-E_{j+1,j}^\pm(z_1)^2\frac{{c}^*(zq^{\mp c},P_{j,j+1}-2)}{\bar{b}^*(zq^{\mp c})}
.\lb{EpEm}
\ena
In a similar way, we obtain
\bea
&&F_{j,j+1}^\pm(z_1)F_{j,j+1}^\pm(z_2)\frac{1}{\bar{b}(z)}
-F_{j,j+1}^\pm(z_1)^2\frac{\bar{c}(z,(P+h)_{j,j+1}-2)}{\bar{b}(z)}
\nn\\
&&=F_{j,j+1}^\pm(z_2)F_{j,j+1}^\pm(z_1)\frac{1}{\bar{b}(1/z)}
-F_{j,j+1}^\pm(z_2)^2\frac{\bar{c}(1/z,(P+h)_{j,j+1}-2)}{\bar{b}(1/z)}
,\lb{FpFp}\\
&&F_{j,j+1}^\pm(z_1)F_{j,j+1}^\mp(z_2)\frac{1}{\bar{b}(zq^{\pm c})}
+ F_{j,j+1}^\pm(z_1)^2\frac{\bar{c}(zq^{\pm c},(P+h)_{1,2}-2)}{\bar{b}(zq^{\pm c})}
\nn\\
&&=F_{j,j+1}^\mp(z_2)F_{j,j+1}^\pm(z_1)\frac{1}{\bar{b}(1/(zq^{\pm c}))}
-F_{j,j+1}^\pm(z_1)^2\frac{\bar{c}(1/(zq^{\pm c}),(P+h)_{j,j+1}-2)}{\bar{b}(1/(zq^{\pm c}))}
.\lb{FpFm}
\ena
From $(j+1,j), (j,j+1)$ of \eqref{mRLLlm},
\begin{eqnarray}
&&\rho^\pm(z)
\left\{
\bar{b}^*(z)K_{j+1}^\pm(z_1)E_{j+1,j}^\pm(z_1)F_{j,j+1}^\pm(z_2)K_{j+1}^\pm(z_2)
\right.\nonumber\\
&&\left.+
\bar{c}(z,(P+h)_{j,j+1})\left(
K_j^\pm(z_1)+F^\pm(z_1)K_{j+1}^\pm(z_1)E_{j+1,j}^\pm(z_1)
\right)K_{j+1}^\pm(z_2)
\right\}\nonumber\\
&=&
\rho^{*\pm}(z)
\left\{
\left(K_j^\pm(z_2)+F_{j,j+1}^\pm(z_2)K_{j+1}^\pm(z_2)E_{j+1,j}^\pm(z_2)\right)K_{j+1}^\pm(z_1)
{\bar{c}^*(z,P_{j,j+1})} \right.\nonumber\\
&&\left.+
F_{j,j+1}^\pm(z_2)K_2^\pm(z_2)K_{j+1}^\pm(z_1)E_{j+1,j}^\pm(z_1)
b^*(z,P_{j,j+1})
\right\}.
\end{eqnarray}
Using $(j+1,j), (j+1,j+1)$ and $(j+1,j+1), (j,j+1)$ of \eqref{mRLLlm}, 
we get
\bea
[E_{j+1,j}^\pm(z_1),F_{j,j+1}^\pm(z_2)]
&=&K^\pm_j(z_2)K^\pm_{j+1}(z_2)^{-1}\frac{\bar{c}^*(z,P_{j,j+1}-1)}{\bar{b}^*(z)}
\nn\\
&&\qquad-K_{j+1}^\pm(z_1)^{-1}K^\pm_j(z_1)\frac{\bar{c}(z,(P+h)_{j,j+1}-1)}{\bar{b}(z)}
.
\lb{EF}
\ena
Similarly, from $(j+1,j), (j,j+1)$ of \eqref{RLLpmlm}, 
\begin{eqnarray}
&&\rho^\pm(zq^{\pm{c}})
\left(
\bar{b}^*(zq^{\pm c})K_{j+1}^\pm(z_1)E_{j+1,j}^\pm(z_1)F_{j,j+1}^\mp(z_2)K_{j+1}^\mp(z_2)
\right.\nonumber\\
&&\left.+
\bar{c}(zq^{\pm c},(P+h)_{j,j+1})(
K_j^\pm(z_1)+F^\pm(z_1)K_{j+1}^\pm(z_1)E_{j+1,j}^\pm(z_1)
)K_{j+1}^\mp(z_2)
\right)\nonumber\\
&=&
\rho^{*\pm}(zq^{\mp c})
\left(
(K_j^\mp(z_2)+F_{j,j+1}^\mp(z_2)K_{j+1}^\mp(z_2)E_{j+1,j}^\mp(z_2))K_{j+1}^\pm(z_1)
{\bar{c}^*(zq^{\mp c},P_{j,j+1})} \right.\nonumber\\
&&\left.+
F_{j,j+1}^\mp(z_2)K_2^\mp(z_2)K_{j+1}^\pm(z_1)E_{j+1,j}^\pm(z_1)
b^*(zq^{\mp c},P_{j,j+1})
\right).
\end{eqnarray}
Using $(j+1,j), (j+1,j+1)$ and $(j+1,j+1), (j,j+1)$ of \eqref{RLLpmlm}, 
we get
\begin{eqnarray}
[E_{j+1,j}^\pm(z_1),F_{j,j+1}^\mp(z_2)]&=&
K_j^\mp(z_2)K_{j+1}^\mp(z_2)^{-1}\frac{\bar{c}^*(zq^{\mp c},P_{j,j+1}-1)}{\bar{b}^*(zq^{\mp c})} 
\nn\\
&&\qquad
-K_{j+1}^\pm(z_1)^{-1}K_j^\pm(z_1)\frac{\bar{c}(zq^{\pm c},(P+h)_{j,j+1}-1)}{\bar{b}(zq^{\pm c})}. 
\lb{EFpm}
\end{eqnarray}

\subsection{Identification with the elliptic currents}
In this section we demonstrate an identification of a combination of the  basic Gauss components with the elliptic currents of $U_{q,p}(\glnh)$. 

Let us define the total current by $E_j(zq^{j-c/2})=\mu^*(E^+_{j+1, j}(z^+)-E^-_{j+1,j}(z^-))$ with 
$\mu^*=\frac{(p^*;p^*)_\infty^3}{a^*_{j+1,j}\Theta_{p^*}(q^2)}.$
Then we obtain 
from \eqref{KpEpKp} and \eqref{KpEmKp}
\bea
&&K^+_{j+1}(z_1)^{-1}E_j(z_2q^{j-c/2})K_{j+1}^+(z_1)
=K^+_{j+1}(z_1)^{-1}\mu^*\left(E_{j+1,j}^+(z^+_2)-E_{j+1,j}^-(z_2^-)\right)K_{j+1}^+(z_1)\nn\\
&&=\mu^*\left(E_{j+1,j}^+(z^+_2)\frac{1}{\bar{b}^*(zq^{- c/2})}
-E_{j+1,j}^+(z_1)\frac{c^*(zq^{- c/2},P_{j,j+1})}{\bar{b}^*(zq^{- c/2})}\right.
\nn\\
&&\left. \qquad\qquad\qquad
 -E_{j+1,j}^-(z^-_2)\frac{1}{\bar{b}^*(zq^{- c/2})}
 +E_{j+1,j}^+(z_1)\frac{c^*(zq^{- c/2},P_{j,j+1})}{\bar{b}^*(zq^{- c/2})}
 \right)\nn\\
&&=E_j(z_2q^{j-c/2})\frac{1}{\bar{b}^*(zq^{- c/2})}
.\lb{KpEKp}
\ena
Comparing this with \eqref{kjp1eja}, we identify $E_j(z), K^+_j(z)$ with $e_j(z), k^+_j(z)$, respectively.  

Next, inserting \eqref{EpEp} and \eqref{EpEm} into $(E^+_{j+1,j}(z^+_1)-E_{j+1,j}(z_1^-))(E^+_{j+1,j}(z^+_2)-E_{j+1,j}(z_2^-))$, 
we obtain
\be
&&
E_j(z_1)E_j(z_2)=
-\frac{1}{z}\frac{\tes(zq^2)}{\tes{(q^2/z)}}E_j(z_2)E_j(z_1). 
\en
This is consistent to \eqref{eethetajj}.  

To obtain \eqref{eethetajjp1} we use \eqref{elelp1}, which is derived from \eqref{RLLpmlminv} with $m=l+2$.  

Similarly, if we define $F_j(zq^{j-c/2})=\mu(F_{j,j+1}^+(z_-)-F_{j,j+1}^-(z_+))$ with 
$\mu=\frac{(p;p)_\infty^3}{a_{j,j+1}\Theta_p(q^2)}$, we recover the relations \eqref{ffthetajj}-\eqref{ffthetajjp1} 
from \eqref{KpFpKp}-\eqref{KpFmKp}, \eqref{FpFp}-\eqref{FpFm} and the $F^+$ counterpart of \eqref{elelp1}. Hence we identify $F_j(z)$ with $f_j(z)$. 

Finally let us check the relation \eqref{eifj}. 
From \eqref{EF} and \eqref{EFpm}, we have 
\bea
&&(\mu\mu^*)^{-1}[E(z_1),F(z_2)]\nn\\
&&=[E_{j+1,j}^+(z_1^+)-E_{j+1,j}^-(z_1^-), F_{j,j+1}^+(z_2^-)-F_{j,j+1}^-(z_2^+)]\nn\\
&&=[E_{j+1,j}^+(z_1^+), F_{j,j+1}^+(z_2^-)]+[E_{j+1,j}^-(z_1^-), F_{j,j+1}^-(z_2^+)]\nn\\
&&\qquad-[E_{j+1,j}^+(z_1^+),F_{j,j+1}^-(z_2^+)]-[E_{j+1,j}^-(z_1^-),F_{j,j+1}^+(z_2^-)].\lb{bracketEF}
\ena
Let us substitute \eqref{EF} and \eqref{EFpm} into this.  Noting the remark below Proposition \ref{derRLL},  
one finds that the terms containing $K^-_j(z_2^+)K^{-}_{j+1}(z_2^+)^{-1}$ from the 2nd and 3rd terms in \eqref{bracketEF} cancel out each other and the same is true 
for the terms containing $K^+_{j+1}(z_1^+)^{-1}K^{+}_j(z_1^+)$
 from the 1st and the 3rd terms in \eqref{bracketEF}.  We obtain
 \bea
&&\hspace{-1cm}(\mu\mu^*)^{-1}[E_j(z_1),F_j(z_2)]\nn\\
&&\hspace{-1cm}=K^+_j(z_2^-)K^{+}_{j+1}(z_2^-)^{-1}q^{-1}\tes(q^2)
\left(\left. \frac{\tes(q^{-2(P_{j,j+1}-1)}q^{c}z)}{\tes(q^{-2(P_{j,j+1}-1)})\tes(zq^{c})}\right|_{+}
-\left. \frac{\tes(q^{-2(P_{j,j+1}-1)}q^{c}z)}{\tes(q^{-2(P_{j,j+1}-1)})\tes(zq^{c})}\right|_{-}\right)\nn\\
&&\hspace{-1cm} - K^-_{j+1}(z_1^-)^{-1}K^{-}_j(z_1^-)q^{-1}\te(q^2)\left(\left. 
\frac{\te(q^{-2((P+h)_{j,j+1}-1)}q^{-c}z)}{\te(q^{-2((P+h)_{j,j+1}-1)})\te(zq^{-c})}
\right|_+-\left. \frac{\te(q^{-2((P+h)_{j,j+1}-1)}q^{-c}z)}{\te(q^{-2((P+h)_{j,j+1}-1)})\te(zq^{-c})}
\right|_-\right).\nn
 \ena
Here 
\bea
\left.\frac{\te(q^{2s}z)(p;p)^3_\infty}{\te(q^{2s})\te(z)}\right|_+&=&
\frac{\te(q^{2s}z)(p;p)^3_\infty}{\te(q^{2s})\te(z)}\nn\\
&=&\sum_{n\in \Z}\frac{1}{1-q^{2s}p^n}z^{n}\nn\\
&=&\sum_{l\in \Z_{\geq 0}}\left(\frac{q^{2sl}}{1-{p^l}{z}}-
\frac{q^{-2s(l+1)}p^{l+1}/z}{1-p^{l+1}/z}\right)\lb{zpositive}
\ena
for $|p|<|z|<1$ and  
\bea
\left.\frac{\te(q^{2s}z)(p;p)^3_\infty}{\te(q^{2s})\te(z)}\right|_-
&=&-\frac{\te(q^{-2s}/z)(p;p)^3_\infty}{\te(q^{-2s})\te(1/z)}\nn\\
&=&-\sum_{n\in \Z}\frac{1}{1-q^{-2s}p^n}z^{-n}\nn\\
&=&-\sum_{l\in \Z_{\geq 0}}\left(\frac{q^{-2sl}{p^l}/{z}}{1-{p^l}/{z}}-
\frac{q^{2s(l+1)}}{1-p^{l+1}z}\right)\lb{znegative}
\ena
for $1<|z|<|p^{-1}|$. Then the difference between these two expansions turns out to be the formal delta function $\delta(z)=\sum_{n\in \Z}z^n$.  
We thus obtain
\bea
&&[E_j(z_1),F_j(z_2)]
\nn\\
&=&\mu\mu^*q^{-1}\left\{K_j^+(z_2^-)K_{j+1}^+(z_2^-)^{-1}\frac{\tes(q^2)}{(p^*;p^*)_\infty^3}
\delta\left(\frac{z_1}{z_2}q^c\right)
 -K_{j+1}^-(z_1^-)^{-1}K_j^-(z_1^-)\frac{\te(q^2)}{(p;p)_\infty}
\delta\left(\frac{z_1}{z_2}q^{-c}\right)\right\}\nn\\
&=&-\frac{\kappa}{q-q^{-1}}\left\{
K_j^+(z_2^-)K_{j+1}^+(z_2^-)^{-1}
\delta\left(\frac{z_1}{z_2}q^c\right)
-K_j^-(z_1^-)K_{j+1}^-(z_1^-)^{-1}
\delta\left(\frac{z_1}{z_2}q^{-c}\right)\right\}.
 \ena
In the last line we used \eqref{Kpm2Kpm2} and 
\be
&&q^{-1}\mu\mu^*\frac{\tes(q^2)}{(p^*;p^*)_\infty^3}=-\frac{\kappa}{q-q^{-1}}.
\en
This is consistent to \eqref{eifj}.

\section{Relation to Jimbo-Miwa-Okado\rq{}s Face Weight}\lb{JMOsW}
In this Appendix we give a relationship between our 
$R$-matrix \eqref{tRmat} and Jimbo-Miwa-Okado\rq{}s $A^{(1)}_{n-1}$ type face weight in \cite{JMO}.

\subsection{Fractional powers in $z$}
So far we have defined the elliptic algebras $U_{q,p}(\glnh)$ and $E_{q,p}(\glnh)$ 
by using  the generating functions $e_j(z), f_j(z), k_j(z)$ and $L_{ij}(z)$, respectively. The coefficients of their relations are given in terms of the theta function 
$\Theta_{p}(z)$, which enables us to expand every things to power series in $p$. 

However for a practical use it is  convenient  to introduce operators such as 
$z^{\pm \frac{P+h-1}{r}}$ and $z^{\pm \frac{P-1}{r^*}}$ into the algebras. 
Here $r$ and $r^*$ are introduced by $p=q^{2r}$ and $p^*=pq^{-2c}=q^{2r^*}$ with $r^*=r-c$. 
The  main reason for this can be seen in the following example.  
Let us consider Jacobi's theta function
\be
&&\vartheta_1(u,\tau)=i\sum_{n=-\infty}^\infty(-1)^ne^{\pi i\tau(n-1/2)^2}e^{2\pi i u(n-1/2)}.
\en
Identifying $z=q^{2u}$, $p=e^{-\frac{2\pi i}{\tau} }$ we have 
\be
&&\vartheta_1(u,\tau)=e^{\frac{\pi i}{4}}\tau^{-\frac{1}{2}}p^{\frac{1}{8}}q^{\frac{u^2}{r}-u}\Theta_p(z).
\en
Then we have 
\bea
&&z^{\frac{s-1}{r}}c(z,s)=z^{\frac{s-1}{r}}\frac{\Theta_p(q^2)\Theta_p(q^{2s}z)}
{\Theta_p(q^{2s})\Theta_p(q^{2}z)}=\frac{\vartheta_1(\frac{1}{r},\tau)\vartheta_1(\frac{s+u}{r},\tau)}
{\vartheta_1(\frac{s}{r},\tau)\vartheta_1(\frac{1+u}{r},\tau)}.\lb{cu}
\ena
This is invariant under the shift $z \mapsto z p$ i.e. $u \mapsto u+r$. 
Compare this with \eqref{shiftR}.

Motivated by this, let us consider the transformation
\bea
\hR^\pm(z,s)=\left(\Ad z^{-\frac{1}{r}\theta_V(s)}\otimes \id\right)
\left(z^{\frac{1}{r}T_{V,V}}
R^\pm(z,s)\right),
\lb{hatR}
\ena
where $\theta_V(s)$ is given in \eqref{thetaV} and
\be
&&T_{V,V}=\sum_{j=1}^{N-1}\pi_V(h_j)\otimes \pi_V(h^{j})
\en
Then one finds 
\bea
\hR^{\pm}(u,s)&=&\hrho^{\pm}(u)\widehat{\bar{R}}(u,s),
\lb{Rp}\\
\widehat{\bar{R}}(u,s)
&=&
\sum_{j=1}^{N}
E_{jj}\otimes E_{jj}+
\sum_{1 \leq j<l \leq N}
\left(\hb(z,s_{j,l})
E_{jj}
\otimes E_{ll}+
\widehat{\bar{b}}(z)
E_{ll}\otimes E_{jj}
\right)\nonumber\\
&&+
\sum_{1 \leq j<l \leq N}
\left(\hc_{}(z,s_{j,l})
E_{jl}\otimes E_{lj}+
{\hc}
(z,-s_{j,l})E_{lj}\otimes E_{jl}
\right),\lb{rmat}
\end{eqnarray}
with
\begin{eqnarray}
&&\hrho^+(u)=q^{-\frac{N-1}{N}}
z^{\frac{N-1}{rN}}
\frac{\{q^2z\}
\{q^{2N-2}z\}
\{p/z\} \{pq^{2N}/z\}
}{
\{z\}
\{q^{2N}z\}
\{pq^2/z\} \{pq^{2N-2}/z\}
},\label{def:rhop}\\
&&\hrho^-(u)=q^{\frac{N-1}{N}}
z^{\frac{N-1}{rN}}
\frac{\{pq^2z\}
\{pq^{2N-2}z\}
\{1/z\} \{q^{2N}/z\}
}{
\{pz\}
\{pq^{2N}z\}
\{q^2/z\} \{q^{2N-2}/z\}
},\label{def:rhom}\\
&&\hb_{}(u,s)=
\frac{[s+1]
[s-1][u]
}{
[s]^{2} [u+1]
},
~~~\widehat{\bar{b}}_{}(u)=\frac{[u]}{[u+1]},\\
&&\hc_{}(u,\pm s)=
\frac{[1][s+u]}{
[s][u+1]},
\lb{rmatcomp}
\end{eqnarray}
where we set 
\be
[u]=\vartheta_1\left(\frac{u}{r},\tau\right).
\en
We also need for $p^*=e^{-\frac{2\pi i}{\tau^*}}$
\be
[u]^*=\vartheta_1\left(\frac{u}{r^*},\tau^*\right).
\en  
We have 
\bea
&&\rho^+(u+r)=\rho^-(u),\qquad 
\frac{\rho^+(u)}{\rho^{+*}(u)}=\frac{\rho^-(u)}{\rho^{-*}(u)},\qquad 
\hrho^+(u)\hrho^-(-u)=1\lb{rhoprhom},\\
&&\hR^-(u,s)^{-1}=P\hR^+(-u,s)P\lb{RpRm}
\ena
for $\hrho^\pm(u)=\hrho^\pm(u)|_{r\to r^*, p\to p^*}$. Hence we obtain
\begin{prop}
\bea
&&\hR^+(u+r,s)=\hR^-(u,s).
\ena
\end{prop} 
Let us further define $\hL^\pm(u)$ from $L^+(z)$ in Sec.\ref{secEqp} by
\bea
&&\hL^+(u)= \left(z^{-\frac{1}{r}\theta_V(P)}\otimes \id\right)
z^{\frac{1}{r}T_{V}}
L^+(z)\left(z^{\frac{1}{r^*}\theta_V(P)}\otimes \id\right)
,\lb{defhL}\\
&&\hL^-(u)=\hL^+(u+r^*+{c}/{2}),\nn
\ena
where $z=q^{2u}$, and $T_V$ is given in \eqref{TV}.
Note that we do not need the extra $\Ad\ q^{-2\theta_V(P)}$ action to define
$\hL^-(u)$ unlike \eqref{LmfromLp}. 

\begin{prop}
 $\hL^\pm(u)$ satisfy the same relations as \eqref{RLL}, \eqref{mRLL} and 
\eqref{RLLpm} with replacing $R^\pm(z,s)$ by $\hR^\pm(u,s)$. Namely, 
\bea
&&\hR^\pm(u,P+h)\hL^\pm(u_1)\hL^\pm(u_2)=\hL^\pm(u_2)\hL^\pm(u_1)\hR^{\pm*}(u,P),\lb{hRLL}\\
&&\hR^\pm(u\pm{c}/{2},P+h)\hL^\pm(u_1)\hL^\mp(u_2)=\hL^\mp(u_2)\hL^\pm(u_1)\hR^{\pm*}(u\mp{c}/{2},P),\lb{hRLLpm}
\ena
where $u=u_1-u_2$.
\end{prop}

\noindent
{\it Remark.}\ The $R$-matrices and $L$-operators in the previous works \cite{Konno,JKOS,KK03,Konno09} are $\hR\pm(u,s)$ and $\hL^\pm(u)$ 
given in this section except for the prefactor: $\rho^\pm(u)$ in the previous works are $\hrho^\mp(u)$ in \eqref{def:rhom} and \eqref{def:rhop}.

\subsection{Gauge transformation from Jimbo-Miwa-Okado's $W$}

For $1\leq j\leq N$ let us define $\displaystyle{F(P,P+\hj)
}$ by
\be
F(s,s+\hj)&=&\left(\prod_{m=j+1}^{N}\frac{[s_{j,m}+1]}{[s_{j,m}]}\right)^{\frac{1}{2}},
\en
where $s_{j,m}=s_{\bep_j}-s_{\bep_m}$ as in Sec.2. 

Our $R$-matrix is related to Jimbo-Miwa-Okado's $W$ denoted by $W_{JMO}$ as follows.
Let us set 
\be
&&\hR^+(u,s)^{ij}_{kl}=W\left(\left.\mmatrix{s&s+\hi\cr s+\hl&a+\hi+\hj}\right|\ u\right)
\qquad (i+j=k+l).
\en
Then we have 
\be
&&W\left(\left.\mmatrix{s&s+\hi\cr s+\hl&s+\hi+\hj}\right|\ u\right)
=\hrho^+(u)\frac{[1]}{[u+1]}\frac{F(s,s+\hi)F(s+\hi,s+\hi+\hj)}{
F(,s+\hl)F(s+\hl,s+\hl+\hk)}\\
&&\qquad\qquad\qquad\qquad \qquad\qquad\qquad \qquad \times W_{JMO}\left(\left.\mmatrix{s&s+\hi\cr s+\hl&s+\hi+\hj}\right|\ u\right).
\en
Here $(s+\hi)_{\bep_j}=s_{\bep_i}+\delta_{i,j}$ etc. and 
\be
&&W_{JMO}\left(\left.\mmatrix{s&s+\hi\cr s+\hi&s+2\hi}\right|\ u\right)=\frac{[1+u]}{[u]},\\
&&W_{JMO}\left(\left.\mmatrix{s&s+\hi\cr s+\hi&s+\hi+\hj}\right|\ u\right)=\frac{[s_{i,j}-u]}{[s_{i,j}]},\\
&&W_{JMO}\left(\left.\mmatrix{s&s+\hj\cr s+\hi&s+\hi+\hj}\right|\ u\right)=\frac{[u]}{[1]}\left(\frac{[s_{i,j}+1][s_{i,j}-1]}{[s_{i,j}]^2}\right)^{1/2}.
\en

\section{Elliptic Quantum Determinant}\lb{ellqdet}

\subsection{Jimbo-Kuniba-Miwa-Okado's projection}
Quantum determinant of the $L$-operator depends on a choice of the gauge for the $R$-matrix. 
Let  $\cV=\FF\otimes_{\C}V,\  V=\oplus_{a=1}^N \C v_a$ and $E_{i,j}v_a=\delta_{j,a}v_i$. 
Let us consider the following $R$ matrix $R'(z,s)\in \End(\cV\otimes \cV)$  given by 
\bea
R'(u,s)=\rho^{+'}(u)\left[\sum_{j=1}^{N}\alpha(u)
E_{jj}\otimes E_{jj}+
\sum_{j\not= l}
\biggl(
\beta_{}(u,s_{l,j})
E_{jj}
\otimes E_{ll}+
\gamma_{}(u,s_{l,j})
E_{jl}\otimes E_{lj}\biggr)\right],\lb{Rprime}
\ena
where 
\be
&&\rho^{+'}(u)=\hrho^+(u)\frac{[1]}{[u+1]},\label{rhop}\\
&&\alpha(u)=\frac{[u+1]}{[1]},\quad 
\beta(u,s)=\frac{[u][s+1]}{[1][s]},\quad
\gamma(u,s)=\frac{[s-u]}{[s]}.
\en
Note that $R'(u,s)^{ij}_{kl}=W\left(\left.\mmatrix{s&s+\widehat{i}\cr s+\widehat{l}&s+\widehat{i}+\widehat{j}}\right|u\right)$ is the face weight in \cite{JKMO}, which is gauge equivalent to $\hR^+(u,s)$ in \eqref{Rp}. 

Instead of giving the gauge transformation of the $R$ matrix  we give  the gauge transformation of the $L$-operator $\hL^+(u)$ in \eqref{defhL} to the one satisfying the $RLL$ relation with the new $R$ matrix \eqref{Rprime}. Namely we define $L(u)=\sum_{1\leq i,j\leq N}E_{i,j}L_{ij}(u)$  by
\bea
L_{ij}(u)=\prod_{m=i+1}^N\frac{[(P+h)_{im}+1]}{[1]}\ \hL^+_{ij}(u)\ 
\prod_{n=1}^{j-1}\frac{[1]^*}{[P_{nj}+1]}.\lb{gaugeL}
\ena
Then  from \eqref{hRLL}  we obtain
\begin{prop}
\bea
&&{R\rq{}}^{(12)}(u_1-u_2,P+h)L^{(1)}(u_1)L^{(2)}(u_2)=L^{(2)}(u_2)L^{(1)}(u_1){R\rq{}}^{*(12)}(u_1-u_2,P).\lb{RpLL}
\ena
\end{prop}

One has 
\be
&&R\rq{}(-1,s)=\rho_0\sum_{j\not= l}\frac{[s_{l,j}+1]}{[s_{l,j}]}\biggl(
E_{jj}\otimes E_{ll}-
E_{jl}\otimes E_{lj}\biggr),
\en
where 
\be
&&\rho_0=-\lim_{u\to -1\ (z\to q^{-2})}\hrho^+(u).
\en
Hence 
\be
&&R\rq{}(-1,s)v_a\otimes v_b=\rho_0\left(\frac{[s_{b,a}+1]}{[s_{b,a}]}v_a\otimes v_b-\frac{[s_{a,b}+1]}{[s_{a,b}]}v_b\otimes v_a\right)\in \cV\wedge \cV. 
\en

In order to generalize this it is convenient to consider the `transposition'  of $R'(u,s)$:
\bea
&&R(u,s)={{}^{t_1t_2}R\rq{}}^{(21)}(u,s), \qquad R^*(u,s)=R(u,s)|_{r\to r^*, p\to p^*}.
\ena
In fact this yields
\be
&&R(-1,s)=\rho_0\sum_{j\not= l}\frac{[s_{l,j}+1]}{[s_{l,j}]}\biggl(
E_{ll}\otimes E_{jj}-
E_{jl}\otimes E_{lj}\biggr).
\en
Hence
\be
&&R(-1,s)v_a\otimes v_b=\rho_0\frac{[s_{a,b}+1]}{[s_{a,b}]}\left(v_a\otimes v_b-v_b\otimes v_a\right)\in \cV\wedge \cV . 
\en
Accordingly taking the transpositions ${}^{t_1}$ and ${}^{t_2}$ of \eqref{RpLL},  flipping the two tensor components and exchanging $u_1$ and $u_2$, we obtain
\bea
&&{R}^{*(12)}(u_2-u_1,P)\ {}^tL^{(1)}(u_1){}^tL^{(2)}(u_2)={}^tL^{(2)}(u_2){}^tL^{(1)}(u_1){R}^{(12)}(u_2-u_1,P+h).\lb{tRpLL}
\ena

Let us generalize these formulas as follows\cite{JKMO,Hasegawa,FV97}. 
For $2\leq k\leq N$, define\\
 $R(u_1,\cdots,u_{k-1};u_k,s)^{1\cdots k-1;k}$ and $\Pi_{k}(u_1,\cdots,u_{k-1};u_k,s)$ $\in \End_\FF(\cV^{\otimes k})$ 
by
\bea
&&\hspace{-1cm}R(u_1,\cdots,u_{k-1};u_{k},s)^{1 \cdots k-1;k}=R^{k-1k}(u_k-u_{k-1},s)
R^{k-2k}(u_k-u_{k-2},s+h^{(k-1)})\nn\\
&&\qquad\qquad\qquad\qquad\qquad\qquad\cdots R^{1k}(u_k-u_1,s+\sum_{j=2}^{k-1}h^{(j)}),\\
&&\hspace{-1cm}\Pi_{k}(u_1,\cdots,u_{k-1};u_k,s)=\frac{1}{k!}
R(u_1,\cdots,u_{k-1};u_{k},s)^{1\cdots k-1;k}
R(u_1,\cdots,u_{k-2};u_{k-1},s)^{1\cdots k-2;k-1}\nn\\
&&\qquad\qquad\qquad \qquad \qquad\qquad\qquad \qquad\cdots  R(u_1,u_2;u_3,s)^{12;3}
R(u_1;u_2,s)^{1;2}.\lb{defPi}
\ena  
We also need $\Pi^*_k(s)$ defined by the same  formula as \eqref{defPi} with replacing $R(u,s)$ by $R^*(u,s)$. 
By using the DYBE \eqref{dybe} repeatedly, one obtains another expression of 
$\Pi_{k}(u_1,\cdots,u_{k-1};u_k,s)$
\begin{prop}
\be
&&\Pi_{k}(u_1,\cdots,u_{k-1};u_k,s)=\frac{1}{k!}R(u_1;u_{2},\cdots,u_k,s)^{1;2 \cdots k}
R(u_2;u_{3},\cdots,u_k,s)^{2;3 \cdots k}\cdots 
R(u_{k-1};u_k,s)^{k-1k}.
\en
where for $j<k$ 
\be
R(u_j;u_{j+1},\cdots,u_k,s)^{j;j+1 \cdots k}&=&R^{jj+1}(u_{j+1}-u_j,s+\sum_{\stackrel{i=1}{\mbox{\tiny $\not=j,j+1$}}}^{k}h^{(i)})
R^{jj+2}(u_{j+2}-u_j,s+\sum_{\stackrel{i=1}{\mbox{\tiny $\not=j,j+1,j+2$}}}^{k}h^{(i)})\\
&&\qquad\qquad \cdots  R^{jk}(u_k-u_{j},s).
\en
\end{prop}

\begin{prop}
Let $L(z)$ be the $L$ operator defined by \eqref{gaugeL}. Then we have 
\bea
&&\Pi^*_{k}(u_1,\cdots,u_{k-1};u_k,P)\ \trL^{(1)}(u_1)\trL^{(2)}(u_2)\cdots \trL^{(k)}(u_k)\nn\\
&&\qquad =\trL^{(k)}(u_k)\trL^{(k-1)}(u_{k-1})\cdots \trL^{(1)}(u_1)\Pi_{ k}(u_1,\cdots,u_{k-1};u_k,P+h-\sum_{j=1}^kh^{(j)}).\lb{PiLLPi}
\ena
\end{prop}
\noindent
Note that 
\be
&&R^{ij}(u,s+h^{(i)}+h^{(j)})=R^{ij}(u,s).
\en

Now let us consider the operators $\Pi_k(u_1,\cdots,u_{k-1};u_k,P+h)$  and $\Pi^*_{k}(u_1,\cdots,u_{k-1};u_k,P)$ with 
the specialization of the spectrum parameters $(u_1,u_2,\cdots,u_k)=(u,u-1, \cdots, u-(k-1))$. We denote the resultant operators by     
$\Pi_k(P+h)$ and $ \Pi^*_k(P)$, respectively. 
Let us set $[1,N]=\{1,2,\cdots,N\}$, $I=\{i_1,i_2,\cdots, i_k\}\subseteq [1,N]$ with $ i_1<i_2<\cdots < i_k$ and define 
\be
&&v_I=\cC_I\ v_{i_1}\wedge v_{i_2}\wedge \cdots \wedge v_{i_k},\\
&&v_{i_1}\wedge v_{i_2}\wedge \cdots \wedge v_{i_k}=\frac{1}{k!}\sum_{\s\in \gS_k}\sgn \s 
 \ v_{i_{\s(1)}}\tot v_{i_{\s(2)}}\tot \cdots \tot v_{i_{\s(k)}},\\
&&\cC_I=\prod_{1\leq a<b\leq k}\sqrt{\frac{\rho_0^*[a]^*}{[1]^*}\frac{\rho_0[a]}{[1]}}\frac{[(P+h)_{i_a,i_b}+1]}{[P_{i_a,i_b}]^*}.
\en
\begin{prop}\lb{actPi}
For $2\leq k\leq N$ and $s=P, P+h\in H$, 
\bea
{\rm Im}\ \Pi_{k}(s)={\wedge}^k\cV. 
\ena
In particular, ${\rm Im}\ \Pi_{N}(s)$ is the one dimensional subspace 
of $\cV^{\tot N}$ spanned by 
\be
&&v_{[1,N]}= \cC_{[1,N]} \ v_{1}
\wedge v_{2}\wedge \cdots\wedge v_{N}. 
\en
\end{prop}
\noindent
{\it Proof.}\ 
By induction one has  
\be
\Pi_{k}^*(P)v_{i_1}\tot v_{i_2}\tot \cdots \tot v_{i_k}&=&\prod_{1\leq a<b\leq k}\rho^*_0\frac{[a]^*[P_{i_a,i_b}+1]^*}{[1]^*[P_{i_a,i_b}]^*}
v_{i_1}\wedge v_{i_2}\wedge \cdots \wedge v_{i_k}\quad  \in \wedge^k \cV.
\en
Then noting the identity
\bea
&&\prod_{1\leq a<b\leq k}\rho^*_0\frac{[a]^*[P_{i_a,i_b}+1]^*}{[1]^*[P_{i_a,i_b}]^*}=\cN_I \cC_I
\ena
where 
\bea
&&\cN_I=\prod_{1\leq a<b\leq k}\sqrt{\frac{\rho^*_0[a]^*[1]}{\rho_0[a][1]^*}}\frac{[P_{i_a,i_b}+1]^*}{[(P+h)_{i_a,i_b}+1]},
\ena
 one obtains
\be
&&\Pi_{k}^*(P)v_{i_1}\tot v_{i_2}\tot \cdots \tot v_{i_k}=\cN_I\ v_I.
\en
Similarly using the identity
\bea
&&\prod_{1\leq a<b\leq k}\rho_0\frac{[a][(P+h)_{i_a,i_b}+1]}{[1][(P+h)_{i_a,i_b}]}=\cN_I\rq{} \cC_I,
\ena
with 
\bea
&&\cN_I\rq{}=\prod_{1\leq a<b\leq k}\sqrt{\frac{\rho_0[a][1]^*}{\rho^*_0[a]^*[1]}}\frac{[P_{i_a,i_b}]^*}{[(P+h)_{i_a,i_b}]},
\ena
one obtains
\be
&&\Pi_{k}(P+h)v_{i_1}\tot v_{i_2}\tot \cdots \tot v_{i_k}=\cN_I\rq{}\ v_I.
\en
\qed

Consider the projection operator $A_k: \cV^{\otimes k}\to \wedge^k \cV$ 
\be
&&A_k=\frac{1}{k!}\sum_{\stackrel{1\leq j_1,\cdots,j_k\leq N}{\mbox{\tiny\ $ j_a\not=j_b (a\not= b)$}}}\sum_{\s\in \gS_k}\sgn\s E_{j_{\s(1)},j_1}\tot \cdots \tot E_{j_{\s(k)},j_k}.
\en 
Note that $A_k\Pi_k^*(P)=\Pi_k^*(P)$.
\begin{dfn}
Let $I=\{i_1,i_2,\cdots,i_k\}, J=\{j_1,j_2,\cdots,j_k\}\subset [1,N]$ with $i_a<i_b, j_a<j_b$ for $1\leq a<b\leq k$. We define the quantum minor determinant $l(z)^J_I$ of $L(z)$ by
\be
&&\Pi^*_{k}(P)\ \trL^{(1)}(u)\trL^{(2)}(u-1)\cdots \trL^{(k)}(u-(k-1))v_{i_1}\tot \cdots \tot v_{i_k}\nn\\
&&=A_k\ \trL^{(k)}(u-(k-1))\trL^{(k-1)}(u-{(k-2)})\cdots \trL^{(1)}(u)\Pi_{k}(P+h)v_{i_1}\tot \cdots \tot v_{i_k}\nn\\
&&=\sum_{1\leq j_1<\cdots<j_k\leq N}l(z)_I^J\ v_J.
\en
\end{dfn} 
For $\tau\in \gS_k$ we set $\tau(I)=\{i_{\tau(1)},\cdots,i_{\tau(k)}\}$ and define\cite{Hartwig} 
\be
&&\sgn_I(\tau,P+h)=\prod_{1\leq a<b\leq k\atop \tau(a)>\tau(b)}\frac{[(P+h)_{i_{\tau(a)},i_{\tau(b)}}+1]}{[(P+h)_{i_{\tau(b)},i_{\tau(a)}}+1]}
\qquad \sgn^*_I(\tau,P)=\prod_{1\leq a<b\leq k\atop \tau(a)>\tau(b)}\frac{[P_{i_{\tau(a)},i_{\tau(b)}}+1]^*}{[P_{i_{\tau(b)},i_{\tau(a)}}+1]^*}.
\en
Then we have
\begin{prop}\lb{qdet}
\be
l(u)_I^J&=&\cN_{J}\sum_{\s\in \gS_k}\sgn^*_{J}(\s,P)L_{i_1j_{\s(1)}}(u)L_{i_2j_{\s(2)}}(u-1)\cdots L_{i_kj_{\s(k)}}(u-(k-1))\nn\\
&=&\cN_{J}\rq{}\frac{F_I(P+h)}{F_J(P+h)}\sum_{\s\in \gS_k}\sgn \s\ L_{i_{\s(k)}j_k}(u-(k-1))L_{i_{\s(k-1)}j_{k-1}}(u-(k-2))\cdots L_{i_{\s(1)}j_1}(u),
\lb{qminor}
\en
where
\be
&&F_I(P+h)=\prod_{1\leq a<b\leq k}\frac{[(P+h)_{i_a,i_b}+1]}{[(P+h)_{i_a,i_b}]}.
\en
In particular, in the case $I=J=[1,N]$ we obtain the quantum determinant of $L(z)$:
\be  
\edet L(u)&=&\cN_{[1,N]}\sum_{\s\in \gS_N}\sgn^*_{[1,N]}(\s,P)L_{1\s(1)}(u)L_{2\s(2)}(u-1)\cdots L_{N\s(N)}(u-(N-1))\nn\\
&=&\cN_{[1,N]}\rq{}\sum_{\s\in \gS_N}\sgn \s\ L_{\s(N)N}(u-(N-1))L_{\s(N-1)N-1}(u-(N-2))\cdots L_{\s(1)1}(u).
\lb{qdetcomp}
\en
\end{prop}
\noindent
{\it Proof.} 
The statement follows from a standard calculation given for example in \cite{Molev} and the formulas
\be
&&\sgn^*_{[1,N]}(\s,P)=\prod_{1\leq a<b\leq N\atop \s(a)>\s(b)}\frac{[P_{i_{\s(a)},i_{\s(b)}}+1]^*}{[P_{i_{\s(b)},i_{\s(a)}}+1]^*}= \prod_{1\leq a<b\leq N}\frac{[P_{i_{\s(a)},i_{\s(b)}}+1]^*}{[P_{i_a,i_b}+1]^*},\\
&&\sgn\s =\prod_{1\leq a<b\leq N}\frac{[P_{i_a,i_b}]^*}{[P_{i_{\s(a)},i_{\s(b)}}]^*}=\prod_{1\leq a<b\leq N}\frac{[(P+h)_{i_a,i_b}]}{[(P+h)_{i_{\s(a)},i_{\s(b)}}]}.
\en
\qed

Note that the formulas in Proposition \ref{qdet} are consistent to the ones obtained by Hartwig using the co-module algebras\cite{Hartwig}. 

Now using the identity $\sgn\tau \ \cC_{\tau(I)}/\cC_I=\sgn_I(\tau,P+h)$ we have 
\begin{prop}
\bea
&&v_{\tau(I)}=\sgn_I(\tau,P+h)v_I.
\ena
\end{prop}
Then we obtain
\begin{prop}\lb{ltau}
For $\tau\in \gS_k$,
\bea
&&l(u)_{\tau(I)}^J=\sgn_I(\tau,P+h)\ l(u)_{I}^J,\\
&&l(u)_I^{\tau(J)}=\frac{1}{\sgn_J(\tau,P+h)}\ l(u)_{I}^J.
\ena
\end{prop}
Then noting \eqref{coprof} and Proposition \ref{Coprofsoverf}
we obtain
\begin{prop}
\bea
&&\Delta(l(u)_I^J)=\sum_{1\leq l_1<\cdots<l_k\leq N}\frac{\cN_J}{\cN_L}\ l(u)_I^L\ \tot\  l(u)_L^J.
\ena
In particular,
\bea
&&\Delta(\qdet L(u))=\qdet L(u)\tot \qdet L(u).
\ena
\end{prop}

Next for $\leq l\leq k$ let us set $\widehat{i_l}=I\backslash \{i_l\}$ and define 
\be
&&\cN^{(l)}_{I}=\frac{\cN_I}{\cN_{\widehat{i_l}}},\qquad F^{(l)}_{I}(P+h)=\frac{F_I(P+h)}{F_{\widehat{i_l}}(P+h)}
\en
etc.. Note that 
\be
&&F^{(l)}_{I}(P+h)=\prod_{1\leq a<l}\frac{[(P+h)_{i_a,i_l}+1]}{[(P+h)_{i_a,i_l}]}\prod_{l< a\leq k}\frac{[(P+h)_{i_l,i_a}+1]}{[(P+h)_{i_l,i_a}]}.
\en
Then by a direct calculation using the expressions of $l(u)_I^J$ in Proposition \ref{qdet} we obtain 
\begin{prop}\lb{laplacel}
\be
l(u)_I^J&=&\sum_{l=1}^k\cN^{(l)}_{J}\prod_{l<a\leq k}\frac{[P_{j_a,j_l}+1]^*}{[P_{j_l,j_a}+1]^*}\ l(u)_{\widehat{i_k}}^{\widehat{j_l}}\ L_{i_kj_l}(u-(k-1)),\\
&=&\sum_{l=1}^k L_{i_1j_l}(u)\ l(u-1)_{\widehat{i_1}}^{\widehat{j_l}}\ \cN^{(l)}_{J}\prod_{1\leq a<l}\frac{[P_{j_l,j_a}+1]^*}{[P_{j_a,j_l}+1]^*},\\
&=& \sum_{l=1}^k L_{i_lj_k}(u-(k-1))\ l(u)_{\widehat{i_l}}^{\widehat{j_k}}\ \cN^{(k)}_{J}\frac{F^{(l)}_{I}(P+h)}{F^{(k)}_{J}(P+h)}(-1)^{k-l},\\
&=&\sum_{l=1}^k\cN^{(1)}_{J}\frac{F^{(l)}_{I}(P+h)}{F^{(1)}_{J}(P+h)}(-1)^{k-l}\   l(u-1)_{\widehat{i_l}}^{\widehat{j_1}}\ L_{i_lj_1}(u)
\en
\end{prop}
\noindent

\begin{prop}
For $1\leq i \leq N$, 
\bea
\qdet L(u)&=& \sum_{l=1}^N\cN^{(l)}_{[1,N]}\prod_{i<a\leq N}\frac{[(P+h)_{i,a}+1]}{[(P+h)_{a,i}+1]}\prod_{l<a\leq N}\frac{[P_{a,l}+1]^*}{[P_{l,a}+1]^*}\ l(u)_{\widehat{i}}^{\widehat{l}}\ L_{il}(u-(N-1)),\nn\\
&=&\sum_{l=1}^N L_{il}(u)\ l(u-1)_{\widehat{i}}^{\widehat{l}}\ \cN^{(l)}_{[1,N]}\prod_{1\leq a<i}\frac{[(P+h)_{a,i}+1]}{[(P+h)_{i,a}+1]}\prod_{1\leq a<l}\frac{[P_{l,a}+1]^*}{[P_{a,l}+1]^*},\nn\\
&=& \sum_{l=1}^N L_{li}(u-(N-1))\ l(u)_{\widehat{l}}^{\widehat{i}}\ (-1)^{N-l} \cN^{(i)}_{[1,N]}\frac{F^{(l)}_{[1,N]}(P+h)}{F^{(i)}_{[1,N]}(P+h)}\prod_{i<a\leq N}\frac{[(P+h)_{i,a}+1]}{[(P+h)_{a,i}+1]},\nn\\
&=&\sum_{l=1}^N(-1)^{N-l}\cN^{(i)}_{[1,N]}\frac{F^{(l)}_{[1,N]}(P+h)}{F^{(i)}_{[1,N]}(P+h)}\prod_{1\leq a <j}\frac{[(P+h)_{a,i}+1]}{[(P+h)_{i,a}+1]}\   l(u-1)_{\widehat{l}}^{\widehat{j}}\ L_{li}(u)\nn
\ena
\end{prop}
\noindent
{\it Proof. }\ Consider the case $I=J=[1,N]$ in Proposition \ref{laplacel} and use Proposition \ref{ltau} for 
the cyclic permutations $\tau=(i,i+1,\cdots,N)$ in the 1st and the 4th lines,  whereas for 
$\tau=(i,i-1, \cdots,2,1)$ in the 2nd and the 3rd lines.

\subsection{Gauge transformation}
Inserting \eqref{gaugeL} into the expressions of  $l(u)_I^J$ in Proposition \ref{qdet} we define the quantum minor determinant $\hl^+(u)_I^J$ of $\hL^+(u)$  by 
\bea
l(u)_I^J&=&\left(\prod_{1\leq a<b\leq k}\frac{[1]^*}{[1]}\frac{[(P+h)_{i_a,i_b}+1]}{[(P+h)_{j_a,j_b}+1]}\right)
\ \hl^+(u)_I^J.
\ena
For $\s\in \gS$ we set 
\bea
&&\hsgn_{I}(\s,P+h)=
\prod_{1\leq a<b\leq k\atop \s(a)>\s(b)}\frac{[(P+h)_{i_{\s(a)},i_{\s(b)}}+1]}{[(P+h)_{i_{\s(b)},i_{\s(a)}}]},\qquad  
\hsgn_I^*(\s,P)=\prod_{1\leq a<b\leq k\atop \s(a)>\s(b)}\frac{[P_{i_{\s(a)},i_{\s(b)}}+1]^*}{[P_{i_{\s(b)},i_{\s(a)}}]^*}.\nn
\ena

\begin{prop}\lb{qdethL}
\be
&& \hl^+(u)_I^J\nn\\
&&=\cN_k\sum_{\s\in \gS_k} \hsgn^*_{J}(\s,P)\hL^+_{i_1j_{\s(1)}}(u)\hL^+_{i_2j_{\s(2)}}(u-1)\cdots \hL^+_{i_kj_{\s(k)}}(u-(k-1)), \nn\\
&&=\cN_k^{-1}\sum_{\s\in \gS_k} \hsgn_{I}(\s,P+h)
\hL^+_{i_{\s(k)}j_k}(u-(k-1))\hL^+_{i_{\s(k-1)}j_{k-1}}(u-(k-2))\cdots \hL^+_{i_{\s(1)}j_1}(u),
\en
where 
\be
&&\cN_k=\prod_{1\leq a<b\leq k}\sqrt{\frac{\rho^*_0[a]^*[1]}{\rho_0[a][1]^*}}.
\en
\end{prop}
\noindent
{\it Proof.}\ Inserting \eqref{gaugeL} and using \eqref{lgr} and \eqref{rgr}, we have
\be
&&\hspace{-1cm}L_{i_1j_{\s(1)}}(u)L_{i_2j_{\s(2)}}(u-1)\cdots L_{i_kj_{\s(k)}}(u-(k-1))\\
&&=\prod_{1\leq a<b\leq k}\frac{[(P+h)_{i_a,i_b}+1]}{[1]}\ 
\hL^+_{i_1j_{\s(1)}}(u)\hL^+_{i_2j_{\s(2)}}(u-1)\cdots \hL^+_{i_kj_{\s(k)}}(u-(k-1))\\
&&\qquad\qquad\qquad\qquad\qquad \times\prod_{1\leq a<b\leq k}\frac{[1]^*}{[P_{j_a,j_b}+1]^*}\prod_{1\leq a<b\leq k\atop\s(a)>\s(b)}\frac{[P_{j_{\s(b)},j_{\s(a)}}+1]^*}{[P_{j_{\s(b)},j_{\s(a)}}]^*}\\
&&\hspace{-1cm}L_{i_{\s(k)}j_k}(u-(k-1))L_{i_{\s(k-1)}j_{k-1}}(u-(k-2))\cdots L_{i_{\s(1)}j_1}(u)\\
&&=\prod_{1\leq a<b\leq k}\frac{[(P+h)_{i_a,i_b}]}{[1]}\prod_{1\leq a<b\leq k\atop\s(a)>\s(b)}\frac{[(P+h)_{i_{\s(a)},i_{\s(b)}}+1]}{[(P+h)_{i_{\s(a)},i_{\s(b)}}]}\\
&&\qquad\times \hL^+_{i_{\s(k)}j_k}(u-(k-1))\hL^+_{i_{\s(k-1)}j_{k-1}}(u-(k-2))\cdots \hL^+_{i_{\s(1)}j_1}(u) \prod_{1\leq a<b\leq k}\frac{[1]^*}{[P_{j_a,j_b}]^*}.
\en
\qed
\begin{cor}\lb{qdethLcomp}
\bea
&&\qdet \hL^+(u)\nn\\
&&=\cN_N\sum_{\s\in \gS_N} \hsgn^*_{[1,N]}(\s,P)\hL^+_{1\s(1)}(u)\hL^+_{2\s(2)}(u-1)\cdots \hL^+_{N\s(N)}(u-(N-1)), \nn\\
&&=\cN_N^{-1}\sum_{\s\in \gS_N} \hsgn_{[1,N]}(\s,P+h)
\hL^+_{\s(N)N}(u-(N-1))\hL^+_{\s(N-1)N-1}(u-(N-2))\cdots \hL^+_{\s(1)1}(u).\nn
\ena
\end{cor}

\begin{prop}
\bea
&&\hl^+(u)_{\tau(I)}^J=\hsgn_I(\tau,P+h)\  \hl^+(u)_{I}^J,\\
&&\hl^+(u)_I^{\tau( J)}=\hsgn_J^*(\tau,P)\ \hl^+(u)_{I}^J.
\ena
\end{prop}

\begin{prop}
\bea
&&\Delta(\hl^+(u)_I^J)=\sum_{1\leq l_1<\cdots<l_k\leq k}\ \hl^+(u)_I^L\ \tot\ \hl^+(u)_L^J.
\ena
\end{prop}
\noindent
{\it Proof.}\ Note $\Delta(\cN_k)=\cN_k\tot \cN_k$ and 
\be
&&\Delta(\hsgn_I^*(\s,P))=1\tot\  \hsgn_I^*(\s,P), \qquad \Delta(\hsgn_I(\s,P+h))=\hsgn_I(\s,P+h)\tot 1. 
\en
 \qed

\begin{prop}\lb{LaplaceExp}
\be
\hl^+(u)_I^J&=&\sum_{l=1}^k \cN\rq{}_{k}\prod_{l<a\leq k}\frac{[P_{j_a,j_l}+1]^*}{[P_{j_l,j_a}]^*}\ \hl^+(u)_{\widehat{i_k}}^{\widehat{j_l}}\ \hL^+_{i_kj_l}(u-(k-1)),\\
&=&\sum_{l=1}^k \hL^+_{i_1j_l}(u)\ \hl^+(u-1)_{\widehat{i_1}}^{\widehat{j_l}}\ \cN\rq{}_{k}\prod_{1\leq a<l}\frac{[P_{j_l,j_a}+1]^*}{[P_{j_a,j_l}]^*},\\
&=& \sum_{l=1}^k \hL^+_{i_lj_k}(u-(k-1))\ \hl^+(u)_{\widehat{i_l}}^{\widehat{j_k}}\ {\cN\rq{}_{k}}^{-1}\prod_{l<a\leq k}\frac{[(P+h)_{i_a,i_l}+1]}{[(P+h)_{i_l,i_a}]},\\
&=&\sum_{l=1}^k{\cN\rq{}_{k}}^{-1}\prod_{1\leq a <l}\frac{[(P+h)_{i_l,i_a}+1]}{[(P+h)_{i_a,i_l}]}\   \hl^+(u-1)_{\widehat{i_l}}^{\widehat{j_1}}\ \hL^+_{i_lj_1}(u),
\en
where 
\be
&&\cN_k\rq{}=\frac{\cN_k}{\cN_{k-1}}=\prod_{1\leq a\leq k-1}\sqrt{\frac{\rho^*_0[a]^*[1]}{\rho_0[a][1]^*}}.
\en
\end{prop}
\begin{prop}
For $1\leq i\leq N$,
\be
\qdet \hL^+(u)&=& \sum_{l=1}^N\cN\rq{}_{N}\prod_{i<a\leq N}\frac{[(P+h)_{i,a}]}{[(P+h)_{a,i}+1]}\prod_{l<a\leq N}\frac{[P_{a,l}+1]^*}{[P_{l,a}]^*}\ \hl^+(u)_{\widehat{i}}^{\widehat{l}}\ \hL^+_{il}(u-(N-1)),\nn\\
&=&\sum_{l=1}^N \hL^+_{il}(u)\ \hl^+(u-1)_{\widehat{i}}^{\widehat{l}}\ \cN\rq{}_{N}\prod_{1\leq a<i}\frac{[(P+h)_{a,i}]}{[(P+h)_{i,a}+1]}\prod_{1\leq a<l}\frac{[P_{l,a}+1]^*}{[P_{a,l}]^*},\nn\\
&=& \sum_{l=1}^N \hL^+_{li}(u-(N-1))\ \hl^+(u)_{\widehat{l}}^{\widehat{i}}\   {\cN\rq{}_{N}}^{-1}\prod_{l<a\leq N}\frac{[(P+h)_{a,l}+1]}{[(P+h)_{l,a}]}\prod_{i<a\leq N}\frac{[P_{a,i}+1]^*}{[P_{i,a}]^*},\nn\\
&=&\sum_{l=1}^N{\cN\rq{}_{N}}^{-1}\prod_{1\leq a <l}\frac{[(P+h)_{l,a}+1]}{[(P+h)_{a,l}]}\prod_{1\leq a  <i}\frac{[P_{a,i}]^*}{[P_{i,a}+1]^*}\   \hl^+(u-1)_{\widehat{l}}^{\widehat{i}}\ \hL^+_{li}(u)
\en
\end{prop}
Comparing  this and the axiom for the antipode $S$, we determine the action of $S$ on $\hL^+_{il}(u)$ and $\hl^+(u)_{\widehat{l}}^{\widehat{i}}$. For each there are four different expressions. For example,
\begin{prop}\lb{SL}
\bea
&&S(\hL^+_{il}(u)) =\hl^+(u-1)_{\widehat{l}}^{\widehat{i}}\ \cN\rq{}_{N}\prod_{1\leq a<l}\frac{[(P+h)_{a,l}]}{[(P+h)_{l,a}+1]}\prod_{1\leq a<i}\frac{[P_{i,a}+1]^*}{[P_{a,i}]^*}(\edet \hL^+(u))^{-1},\nn\\
&&S(\hl^+(u)_{\widehat{l}}^{\widehat{i}})=\cN\rq{}_{N}\prod_{i<a\leq N}\frac{[(P+h)_{i,a}]}{[(P+h)_{a,i}+1]}\prod_{l<a\leq N}\frac{[P_{a,l}+1]^*}{[P_{l,a}]^*}(\edet \hL^+(u))^{-1}.\nn
\ena
\end{prop}
Combining these we obtain
\begin{prop}
\be
&&S^2(\hL^+_{il}(u))=\prod_{a\in \widehat{i}}\frac{[(P+h)_{i,a}+1]}{[(P+h)_{i,a}]}\ \hL^+_{il}(u-N)\prod_{a\in \widehat{l}}\frac{[P_{a,l}+1]^*}{[P_{a,l}]^*}. 
\en
\end{prop}
Proposition \ref{LaplaceExp} also yields
\begin{prop}
\be
&&(\hL^+(u)^{-1})_{ij}=S(\hL^+_{ij}(u)). 
\en
\end{prop} 
 
\subsection{Formulas for the half currents}
In this section we follow the idea in \cite{Iohara}. 

For $1\leq a, b\leq N$ let us define $\hL^+(u)_{a,a}=(\hL^+_{i,j}(u))_{a\leq i,j\leq N}$ and
\bea
\hL^+(u)_{a,b}&=&\mat{\hL^+_{ab}(u)&\hL^+_{aa+1}(u)&\cdots&\hL^+_{aN}(u)\cr
\hL^+_{a+1b}(u)&\hL^+_{a+1a+1}(u)&\cdots&\hL^+_{a+1N}(u)\cr
\vdots&\vdots&&\vdots\cr
\hL^+_{Nb}(u)&\hL^+_{Na+1}(u)&\cdots&\hL^+_{NN}(u)\cr} \qquad \mbox{for}\ a>b\\
&=&\mat{\hL^+_{ab}(u)&\hL^+_{ab+1}(u)&\cdots&\hL^+_{aN}(u)\cr
\hL^+_{b+1b}(u)&\hL^+_{b+1b+1}(u)&\cdots&\hL^+_{b+1N}(u)\cr
\vdots&\vdots&&\vdots\cr
\hL^+_{Nb}(u)&\hL^+_{Nb+1}(u)&\cdots&\hL^+_{NN}(u)\cr} \qquad \mbox{for}\ a<b.
\ena
Then we have the following Gauss decompositions.
\begin{lem}
\noindent
For $a>b$
\bea
&&\hspace{-1.5cm}\hL^+(u)_{a,b}=
\left(\begin{array}{ccccc}
1&F_{a,a+1}^+(z)&F_{a,a+2}^+(z)&\cdots&F_{a,N}^+(z)\\
0&1&F_{a+1,a+2}^+(z)&\cdots&F_{a+1,N}^+(z)\\
\vdots&\ddots&\ddots&\ddots&\vdots\\
\vdots&&\ddots&1&F_{N-1,N}^+(z)\\
0&\cdots&\cdots&0&1
\end{array}\right)\left(
\begin{array}{cccc}
K^+_a(z)E^+_{a,b}(u)&0&\cdots&0\\
0&K^+_{a+1}(z)&&\vdots\\
\vdots&&\ddots&0\\
0&\cdots&0&K^+_{N}(z)
\end{array}
\right)\nn\\
&&\qquad\qquad\qquad\qquad\qquad\times
\left(
\begin{array}{ccccc}
1&0&\cdots&\cdots&0\\
E^+_{a+1,b}(z)&1&\ddots&&\vdots\\
E^+_{a+2,b}(z)&
E^+_{a+2,a+1}(z)&\ddots&\ddots&\vdots\\
\vdots&\vdots&\ddots&1&0\\
E^+_{N,b}(z)&E^+_{N,a+1}(z)
&\cdots&E^+_{N,N-1}(z)&1
\end{array}
\right).\lb{Labls}
\ena
For $a<b$
\bea
&&\hspace{-1.5cm}\hL^+(u)_{a,b}=
\left(\begin{array}{ccccc}
1&F_{a,b+1}^+(z)&F_{a,b+2}^+(z)&\cdots&F_{a,N}^+(z)\\
0&1&F_{b+1,b+2}^+(z)&\cdots&F_{b+1,N}^+(z)\\
\vdots&\ddots&\ddots&\ddots&\vdots\\
\vdots&&\ddots&1&F_{N-1,N}^+(z)\\
0&\cdots&\cdots&0&1
\end{array}\right)\left(
\begin{array}{cccc}
F^+_{a,b}(u)K^+_b(z)&0&\cdots&0\\
0&K^+_{b+1}(z)&&\vdots\\
\vdots&&\ddots&0\\
0&\cdots&0&K^+_{N}(z)
\end{array}
\right)\nn\\
&&\qquad\qquad\qquad\qquad\qquad\times
\left(
\begin{array}{ccccc}
1&0&\cdots&\cdots&0\\
E^+_{b+1,b}(z)&1&\ddots&&\vdots\\
E^+_{b+2,b}(z)&
E^+_{b+2,b+1}(z)&\ddots&\ddots&\vdots\\
\vdots&\vdots&\ddots&1&0\\
E^+_{N,b}(z)&E^+_{N,b+1}(z)
&\cdots&E^+_{N,N-1}(z)&1
\end{array}
\right).\lb{Labsl}
\ena
The formula for $\hL^+(u)_{a,a}$ is the same as  \eqref{Laa} with l=a. 
\end{lem}
Let us write the Gauss decomposition of $\hL^+(u)_{a,b}$ in the above Lemma as
\be
&&\hL^+(u)_{a,b}=\cF_{a,b}(u)\cK_{a,b}(u)\cE_{a,b}(u).
\en
Then we have
\be
&&\cF_{a,b}(u)^{-1}=\cK_{a,b}(u)\cE_{a,b}(u)\hL^+(u)_{a,b}^{-1}.
\en
Comparing the $(1,1)$ component of the both sides
 we obtain the following.  
\begin{lem}
\bea
&&K^+_a(u)=\frac{1}{(\hL^+(u)_{a,a}^{-1})_{11}}\qquad \mbox{for}\ a=b, \\
&&E^+_{a,b}(u)=(\hL^+(u)_{a,a}^{-1})_{11}\frac{1}{(\hL^+(u)_{a,b}^{-1})_{11}}\qquad \mbox{for}\ a>b,\\
&&F^+_{a,b}(u)=\frac{1}{(\hL^+(u)_{a,b}^{-1})_{11}}(\hL^+_{b}(u)^{-1})_{11}\qquad \mbox{for}\ a<b.
\ena
\end{lem}
Noting 
\be
&&(\hL^+(u)_{a,b}^{-1})_{11}=(\hl^+(u-1)_{a,b})_{\widehat{1}}^{\widehat{1}}
(\qdet\hL^+(u)_{a,b})^{-1}\cN\rq{}_{N-a+1},\\
&&(\hl^+(u-1)_{a,b})_{\widehat{1}}^{\widehat{1}}=\qdet\hL^+(u-1)_{a+1,a+1},
\en
we have
\begin{thm}
\bea
&&K^+_a(u)=\frac{\qdet\hL^+(u)_{a,a}}{\cN\rq{}_{N-a-1}}\frac{1}{\qdet\hL^+(u-1)_{a+1a+1}},\\
&&E^+_{a,b}(u)=\frac{1}{\qdet\hL^+(u)_{a,a}}{\qdet\hL^+(u)_{a,b}}\qquad \mbox{for}\ a>b,\\
&&F^+_{a,b}(u)={\qdet\hL^+(u)_{a,b}}\frac{1}{\qdet\hL^+(u)_{b,b}}\qquad \mbox{for}\ a<b.
\ena
\end{thm}

\begin{cor}
Let us define
\bea
&&K(u)=K^+_1(u)K^+_2(u-1)\cdots K^+_N(u-(N-1)).
\ena
Then
\bea
&&\qdet\hL^+(u)=\cN_N\ K(u).
\ena
Moreover the $q$-determinant $\qdet \hL^+(u)$ belongs to the center of $E_{q,p}(\glnh)$.
\end{cor}
\noindent
{\it Proof.}\ 
Since $K^+_l(v), E_i(v)=const.(E^+_{i+1,i}(v+c/4)-E^-_{i+1,i}(v-c/4)), F_i(v)=const'.(F^+_{i,i+1}(v-c/4)-F^-_{i,i+1}(v+c/4))\ (1\leq l\leq N, 1\leq i\leq N-1)$ satisfy the same commutation relations as the elliptic currents 
of $U_{q,p}(\glnh)$, $K(u)$ commutes with $K^+_l(v), E_i(v), F_i(v), \FF$ due to Proposition \ref{CenterK}. Hence  by Definition \ref{modeGC} 
$K(u)$ commutes with $K^\pm_l(v), E^\pm_{i+1,i}(v), F^\pm_{i,i+1}(v)$ so that 
$K(u)$ commutes with $\hL^+(v)$. 
\qed

\end{appendix}

\renewcommand{\baselinestretch}{0.7}

\end{document}